\documentclass{article}

\usepackage{arxiv}

\usepackage[utf8]{inputenc} 
\usepackage[T1]{fontenc}    
\usepackage{hyperref}       
\usepackage{url}            
\usepackage{booktabs}       
\usepackage{amsmath}
\usepackage{amsfonts}       
\usepackage{amssymb}
\usepackage{amsthm}
\usepackage{algorithmic}
\usepackage{algorithm}
\usepackage{nicefrac}       
\usepackage{microtype}      
\usepackage{lipsum}		
\usepackage{graphicx,xcolor}
\usepackage[numbers]{natbib}
\usepackage{doi}
\usepackage{cleveref}
\usepackage{subcaption}
\usepackage{svg}
\usepackage{csquotes}
\usepackage{enumerate}
\usepackage{enumitem}
\usepackage{cancel}
\usepackage{todonotes}
\usepackage{lineno}

\usepackage{localmacros}

\title{Pulling back symmetric Riemannian geometry for data analysis}

\author{ 
Willem Diepeveen
\\
	Faculty of Mathematics\\
	University of Cambridge\\
	Cambridge, UK \\
	\texttt{wd292@cam.ac.uk} \\
}

\hypersetup{
pdfauthor={Willem Diepeveen},
pdfkeywords={manifold-valued data, Riemannian manifold, interpolation, dimension reduction, deep learning},
}

\begin{document}
\maketitle

\begin{abstract}
    Data sets tend to live in low-dimensional non-linear subspaces. Ideal data analysis tools for such data sets should therefore account for such non-linear geometry. The symmetric Riemannian geometry setting can be suitable for a variety of reasons. First, it comes with a rich mathematical structure to account for a wide range of non-linear geometries that has been shown to be able to capture the data geometry through empirical evidence from classical non-linear embedding. Second, many standard data analysis tools initially developed for data in Euclidean space can also be generalised efficiently to data on a symmetric Riemannian manifold. A conceptual challenge comes from the lack of guidelines for constructing a symmetric Riemannian structure on the data space itself and the lack of guidelines for modifying successful algorithms on symmetric Riemannian manifolds for data analysis to this setting. This work considers these challenges in the setting of pullback Riemannian geometry through a diffeomorphism. The first part of the paper characterises diffeomorphisms that result in proper, stable and efficient data analysis. The second part then uses these best practices to guide construction of such diffeomorphisms through deep learning. As a proof of concept, different types of pullback geometries -- among which the proposed construction -- are tested on several data analysis tasks and on several toy data sets. The numerical experiments confirm the predictions from theory, i.e., that the diffeomorphisms generating the pullback geometry need to map the data manifold into a geodesic subspace of the pulled back Riemannian manifold while preserving local isometry around the data manifold for proper, stable and efficient data analysis, and that pulling back positive curvature can be problematic in terms of stability.
\end{abstract}




\blfootnote{Our code for the main algorithms is available at \href{https://github.com/wdiepeveen/Pulling-back-symmetric-Riemannian-geometry-for-data-analysis}{https://github.com/wdiepeveen/Pulling-back-symmetric-Riemannian-geometry-for-data-analysis}. }

\section{Introduction }
\label{sec:intro-pull-back-geometry}

\begin{figure}[h!]
    \centering
    \begin{subfigure}{0.48\linewidth}
    \centering
        \includegraphics[width=\linewidth]{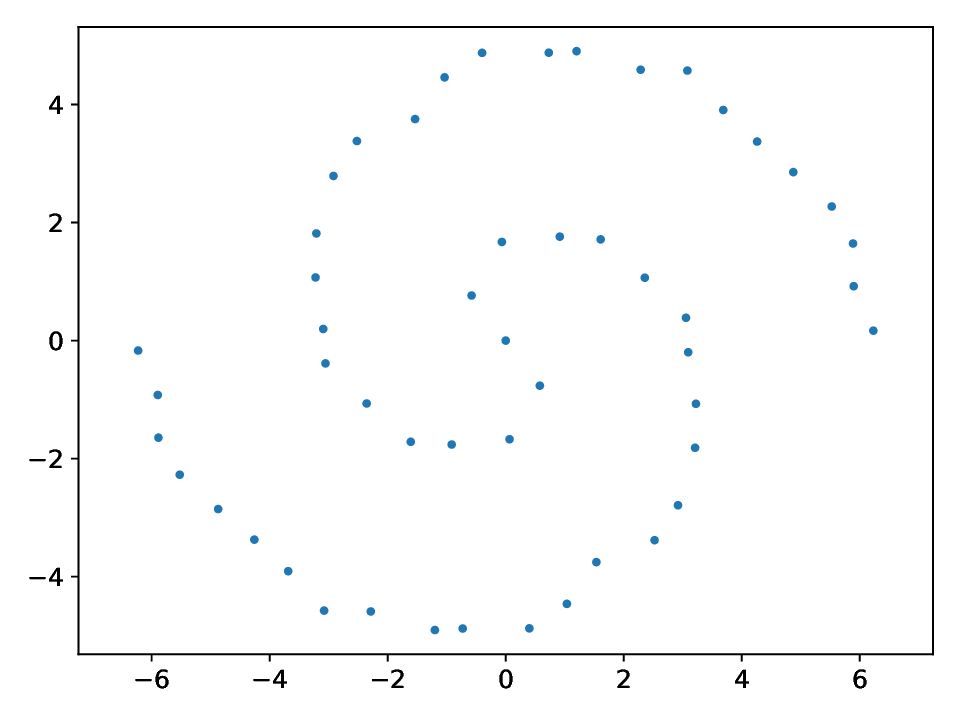}
    \caption{A toy spiral data set.}
    \label{fig:swirl-example-data}
    \end{subfigure}
    \hfill
    \begin{subfigure}{0.48\linewidth}
    \centering
        \includegraphics[width=\linewidth]{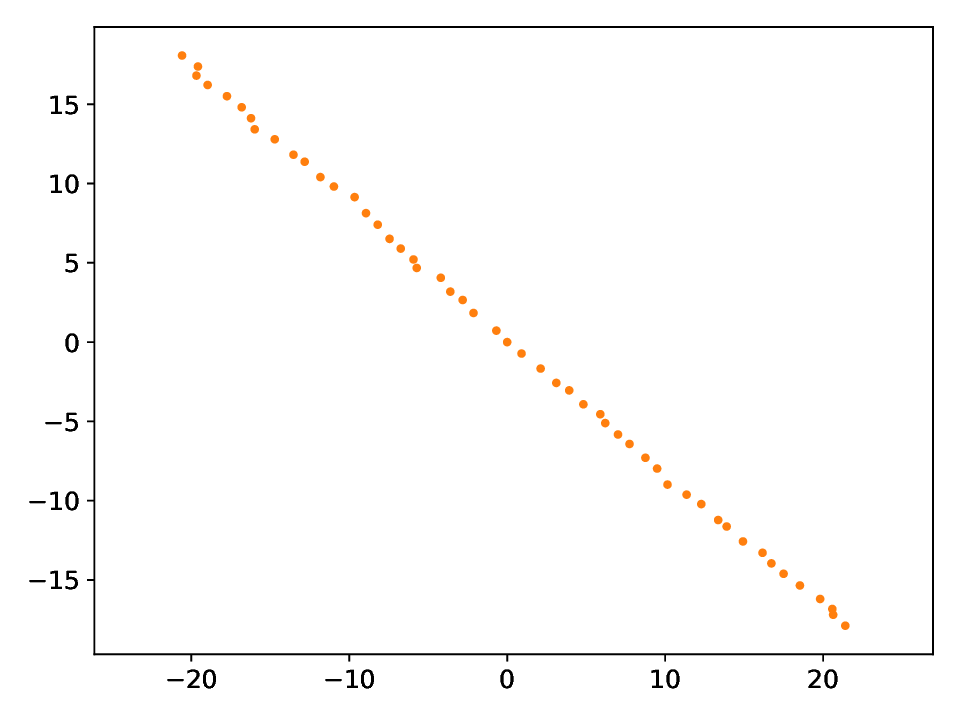}
    \caption{A spiral data set seen from the tangent space $\tangent_{(0,0)} \Real^2$.}
    \label{fig:swirl-example-logs}
    \end{subfigure}
    \caption{The non-linear $\Real^2$-valued data set in (a) looks linear and 1-dimensional from $\Vector:=(0,0)$ under a learned non-standard Riemannian structure on $\Real^2$ (b).}
    \label{fig:swirl-example}
\end{figure}

An increasingly common viewpoint is that data sets in $\Real^\dimInd$ typically reside in non-linear and often low-dimensional subspaces, e.g., as in \cref{fig:swirl-example-data}. To deal with such non-linearity in the data analysis, a suitable notion of distance can already be key for clustering and classification \cite{ghojogh2022spectral,kaya2019deep}, and even hints of the global geometry of the data can be used to improve upon algorithms for further understanding of data sets through geometry-regularized decomposition methods \cite{cai2008non,guan2011manifold,qin2020blind,shu2017multiple,zhang2012low}. Having said that, access to just a distance can be restrictive for tasks beyond clustering and classification and not having a rich enough notion of the global data geometry can be restrictive for fully understanding the data. In particular, a suitable framework on which we do our data analysis ideally at least enables \emph{interpolation and extrapolation over non-linear paths through the data, computing non-linear means on such non-linear paths, and low-rank approximation over curved subspaces spanned by such non-linear paths} as argued in \cite{diepeveen2023riemannian}. Then, from these basic tasks \emph{more advanced signal processing and recovery problems such as non-linear data decomposition and inverse problems can be constructed and solved in a more natural way}. Focusing on the basic tasks for now, we would first like to \emph{have access to the non-linear data space} and \emph{have ways of doing data analysis on there}.


Regarding access to the data space, there is a large amount of work that uses the assumption that the data space is a non-linear \emph{data manifold} and that subsequently applies a non-linear dimension reduction to construct a (Euclidean) embedding manifold \cite{belkin2001laplacian,bishop1998gtm,coifman2006diffusion,demers1992non,hastie1989principal,kingma2013auto,lawrence2005probabilistic,mcinnes2018umap,roweis2000nonlinear,sammon1969nonlinear,scholkopf1998nonlinear,tenenbaum2000global,van2008visualizing,zhang2004principal}. However, there are methodological challenges that prevent us from taking next steps, i.e., data analysis on data manifolds. To see this, consider that in these non-linear embedding approaches there either is a mapping between the embedding manifold and the data manifold available or there is no such mapping. For the latter case, such embeddings do not directly help us setting up a framework capable of the tasks described above since there is no way of mapping interpolants, barycentres and other structural information back to the data manifold. For the former case, there is still some ambiguity over what curves we want to do our data analysis and there is the potential challenge how to communicate this between the embedding manifold and the data manifold. 






The curves we could base our data analysis on can be length-minimising curves -- or \emph{geodesics} -- in the framework of Riemannian geometry \cite{carmo1992riemannian,sakai1996riemannian}, which is a suitable choice for a myriad of reasons. First, all of the aforementioned data analysis tasks have a Riemannian interpretation. Indeed, interpolation can be performed over non-linear geodesics or related higher order interpolation schemes \cite{bergmann2018bezier}, extrapolation can be done using the Riemannian exponential mapping, the data mean is naturally generalized to the Riemannian barycentre \cite{karcher1977riemannian}, and low rank approximation knows several extensions that find the most important geodesics through a data set using the Riemannian logarithmic mapping \cite{diepeveen2023curvature,fletcher2004principal}. In addition, Riemannian geometry also allows to specify the type of non-linearity, i.e., there is a choice which curves are length-minimising – which then also naturally implies the distance between any two points. In particular, there is the possibility to take data geometry into account when constructing a custom Riemannian manifold, which can help in the interpretability \footnote{\Cref{fig:swirl-example-logs} is an example hereof as the data looks low-dimensional and linear from any point in the data set.}. Finally, Riemannian geometry allows to pass geometric information between the embedding manifold and the data manifold as explored in \cite{arvanitidis2017latent,chen2018metrics,louis2019riemannian,rumpf2015variational,shao2018riemannian,tosi2014metrics} -- all of which are in line with the non-linear dimension reduction setting --, but also allows for other approaches through fitting a submanifold and reconstructing Riemannian structure \cite{fefferman2018fitting,fefferman2020reconstructionB,fefferman2019fitting,Yao2023}, or through remetrizing all of the ambient space around the data manifold
\cite{diepeveen2023riemannian,gruffaz2021learning,Lebanon2006,lou2021learning,peltonen2004improved}. 





In practice, we arguably want to impose additional assumptions on the general Riemannian geometry setting. In particular, we want cheap access to basic manifold mappings -- geodesics, exponential mapping, logarithmic mapping, distance -- and we want that the imposed Riemannian structure yields a \emph{symmetric Riemannian manifold}. Whereas the reason for the former requirement on cheap access to mappings is straightforward, the latter requirement needs some motivation. We argue that there are several potential upshots without severe loss of generality. 
Regarding the upshots, there are more efficient tools available for basic data processing such as efficient higher order interpolation \cite{bergmann2018bezier}, efficient low rank approximation \cite{diepeveen2023curvature} and several efficient optimisation schemes 
\cite{bacak2016second,bergmann2021fenchel,bergmann2016parallel,diepeveen2021inexact} that we can consider for even more complicated downstream tasks. In other words, we can piggy back off of a richer literature from Riemannian data processing.
Regarding the loss of generality, it has been observed several times \cite{donoho2005image,lopez2021symmetric,sala2018representation,sonthalia2021can} that symmetric spaces in the setting of data embeddings are a very natural, versatile and rich class of spaces for capturing the geometry of real-world data. So unsurprisingly, there is a whole zoo of non-linear embedding schemes into symmetric Riemannian manifolds such as spherical embeddings \cite{liu2017sphereface,wilson2014spherical} ($\Sphere^\dimInd$), hyperbolic embeddings \cite{keller2020hydra,nickel2017poincare,sonthalia2020tree,walter2004h} ($\Hyperboloid^\dimInd$), combinations of these through product manifold embedding \cite{gu2019learning,skopek2020mixed}, and embedding into other non-constant curvature spaces such as the Grassmann manifold \cite{cruceru2021computationally} ($\mathrm{Gr}(\dimInd,k)$) and the space of symmetric positive definite matrices \cite{cruceru2021computationally} ($\mathcal{P}(\dimInd)$) -- all under their standard Riemannian metric. In other words, if we can harness this empirical observation in the construction of our Riemannian structure, there might be a lot to be gained.








Summarized, we need access to the non-linear data manifold and need ways of doing data analysis on there. From the above discussion we have seen that this can be done, if we can construct a suitable Riemannian manifold that comes with cheap access to manifold mappings and a symmetric Riemannian structure. Our first question is straightforward: \emph{how should we go about constructing such a symmetric Riemannian manifold?} Then, our second question is more subtle: \emph{how should we use or modify successful algorithms on symmetric Riemannian manifolds for data analysis in this setting?} 


\subsection{Related work}
\label{sec:intro-pullback-geometry-related-work}

Out of the above two main questions there is a sizeable literature addressing the first one, whereas the second question has -- to the best of our knowledge -- been left unaddressed. So in the following, we will mainly consider work related to the question of constructing and learning Riemannian manifolds from data. As hinted above, there have been three main ways of approaching this: (i) fitting a submanifold, (ii) constructing a chart, and (iii) remetrizing the ambient space. We will briefly survey these three approaches and the subcategories thereof. In particular, we will focus on the extent to which we can expect to have access to manifold mappings and the possibility to get a symmetric Riemannian structure.


\paragraph{Fitting a submanifold}
is an approach that assumes a data set and a metric structure on this data set, so that we have access to all pairwise distances. The goal is to globally fit a smooth submanifold to the data and to construct a Riemannian structure that approximately realizes the pairwise distances. Even though the developments of this idea are recent, there is solid theoretical understanding of testing such a manifold hypothesis \cite{Fefferman2016}, fitting the manifold \cite{fefferman2018fitting,fefferman2019fitting,Yao2023}, reconstructing the Riemannian structure \cite{fefferman2020reconstructionB}, and necessary data complexity under noise on the pairwise distances \cite{fefferman2020reconstruction}. However, in the end these methods do not come with cheap access to manifold mappings, nor with a Riemannian structure that is symmetric by construction. So for our practical purposes, this approach is not ideal.



\paragraph{Constructing a chart}
is a local version of the above setting, but comes with an explicit embedding mapping between a low-dimensional Euclidean space and the data manifold. In particular, it uses that a manifold is by definition locally homeomorphic to Euclidean space through \emph{charts}. Constructions focus typically on a single chart and have already been proposed in the early 2000s \cite{brand2002charting,lin2008riemannian,roweis2001global,zhang2004principal}. Unlike the above approach that assumes a pairwise metric structure, the geometry is either taken from the the ambient space of the data manifold and then imposed onto the Euclidean embedding manifold, or the Euclidean structure of the embedding manifold is used and then imposed onto the data manifold. The data analysis itself is typically done on the Euclidean embedding manifold rather than on the data manifold.

Both ways of imposing Riemannian geometry have their pros and cons.
Using Riemannian geometry based on the ambient space of the data manifold was pioneered in \cite{tosi2014metrics} and was followed up by works such as \cite{arvanitidis2017latent,chen2018metrics,shao2018riemannian}, but comes with challenges. Realistically, the imposed geometry on the embedding manifold needs to be augmented to actually capture the data geometry properly \cite{arvanitidis2017latent,hauberg2018only} or to capture prior knowledge \cite{arvanitidis2020geometrically}, computational feasibility of geodesics needs to be addressed \cite{arvanitidis2019fast} and so for the general stability of geodesics \cite{arvanitidis2021prior}. In addition, in general it is too optimistic to hope that there is a chart that can cover the whole manifold (a global chart), which can be addressed by either using multiple charts \cite{schonsheck2019chart} or by choosing a non-Euclidean embedding manifold \cite{lopez2022gd}. Besides aforementioned challenges, symmetry of the imposed Riemannian structure is too much to hope for. So overall, despite the large amount of work out, the setting in this line of work is also not ideal for our purposes. Having said that, the second way, which has received much less attention, does come with closed-form geodesics and other manifold mappings, and is symmetric. To the best of our knowledge, the only work that addresses this setting is \cite{louis2019riemannian}, but it also comes with challenges. First, it requires prior knowledge on what distance in the Euclidean embedding space should mean. For example, in \cite{louis2019riemannian} the authors interpret distance as time that has passed between data points. Such an interpretation makes sense in \cite{louis2019riemannian} because the data set of interest contains time series, but in general it is not clear how distances in the Euclidean embedding space should relate to the geometry of the data. In addition, the geometry imposed onto the data manifold is isometric to Euclidean, which limits the types of data geometries that can be modelled with it, e.g., it is impossible to have circular geodesics. Finally, this approach also comes with the single chart issue from before, which is harder to solve than for the previous case as all charts will have their own geometry that needs to be made compatible, compared to using the same geometry from the ambient space of the data manifold. So overall, despite the closed-form manifold mappings and symmetry of the Riemannian structure, this setting is also not ideal for our purposes either.







\paragraph{Remetrizing the ambient space} aims to construct a Riemannian structure on all of the ambient space and while doing that sidesteps the problem of having to find a data manifold first and does not come with topological challenges chart-based approaches suffer from. The basic setting of constructing a data-driven metric tensor field from scratch has been considered in several ways \cite{arvanitidis2016locally,hauberg2012geometric,peltonen2004improved}, but even approximations for the basic manifold mappings can realistically only be expected when having additional special structure from the specific choice of metric \cite{diepeveen2023riemannian,rumpf2015variational}. In these settings, requiring symmetry of the Riemannian structure on top of that is too much to hope for in general. Instead, pulling back or pushing forward geometry through a diffeomorphism that maps between the ambient space and a symmetric Riemannian manifold might be a better way of constructing Riemannian geometry \cite{cuzzolin2008learning,gruffaz2021learning,hauser2017principles,Lebanon2006}. In this setting we do get closed-form manifold mappings and inherit symmetry of the Riemannian structure. In addition, we can impose more general geometries than ones that are isometric to Euclidean \cite{cuzzolin2008learning,Lebanon2006}. Even though this will not necessarily capture every data geometry out there\footnote{For that we can try to embed the ambient space in an even higher dimensional space and pull the geometry back from there \cite{lou2021learning}, but we will lose closed-form mappings and symmetry of the remetrized ambient space.} the empirical evidence mentioned in the above section suggests that this does not have to be a problem in practice. In conclusion, the pullback (or pushforward) approach can be a framework that satisfies our criteria, but similarly to the second chart-based approach it requires prior knowledge on what distance in the symmetric Riemannian manifold should mean. In other words, this boils down to the question: \emph{what are good diffeomorphisms and how do we construct them?}









\subsection{Contributions}
\label{sec:intro-pullback-geometry-contributions}
In the discussion above, we have argued that remetrizing all of ambient space and modeling geometry through pulling back (or pushing forward) from a symmetric Riemannian manifold using a diffeomorphism can be a very suitable setting, in which we want to model data-driven geometry, do data analysis and solve other downstream problems. 


Without loss of generality, we will talk about pullback Riemannian geometry in the paper. In particular, we will consider pullback geometry on $\Real^\dimInd$, i.e., on real-valued data.
Under the pullback geometry framework the two main questions can be specified somewhat more: \emph{how should we go about constructing diffeomorphisms into a symmetric Riemannian manifold?} and \emph{how should we use or modify successful algorithms on symmetric Riemannian manifolds for data analysis in the pullback geometry framework?} As it turns out, it makes a lot of sense to answer these questions jointly, which has not been attempted in prior work as pointed out before. In this work, we make the following three contributions towards answering the above rephrased main questions.

\paragraph{Characterisation of diffeomorphisms for proper and stable data analysis.}
We will see that diffeomorphisms should map the data manifold into low-dimensional geodesic subspaces of the pulled back Riemannian manifold for proper interpolation (\cref{thm:pull-back-mappings}) through and barycentres (\cref{thm:well-posedness-barycentre}) within the data set, but should also do it in such a way that it is a local isometry on the data in order to get $\ell^2$-stability (\cref{thm:stability-geodesics-symmetric,thm:stability-barycentre}), which is a feature we ideally inherit from data analysis in the Euclidean setting in case our data is somewhat noisy. In addition, \cref{thm:stability-geodesics-symmetric,thm:stability-barycentre} show and quantify that the only other instabilities are due to curvature of the pulled back Riemannian manifold. The interplay between curvature effects and diffeomorphism effects on geodesics and barycentres is considered carefully in several numerical experiments with 2-dimensional toy data sets.

\paragraph{Characterisation of diffeomorphisms for efficient data analysis.}

Next, we will see how we can piggyback off of existing theory for data compression on symmetric Riemannian manifolds and how a recent efficient algorithm for low rank approximation motivates construction of the \emph{Riemannian autoencoder} (RAE) and the \emph{curvature corrected Riemannian autoencoder} (CC-RAE) mappings, which are non-linear compression mappings that have additional nice mathematical properties such as an interpretable latent space that traditional neural network-based autoencoders do not have. Although the focus is on algorithm design, we will also observe several times that the above-mentioned best practices for diffeomorpisms are necessary for useful and efficient algorithms. The developed ideas are once again tested in several numerical experiments with 2-dimensional toy data sets.

\paragraph{Construction of diffeomorphisms for proper, stable and efficient data analysis.}

Finally, having learned that diffeomorphisms need to map the data manifold into a geodesic subspace of the pulled back Riemannian manifold while preserving local isometry around the data manifold for proper, stable and efficient data analysis under the pullback geometry, we address how to construct such diffeomorphisms in a general setting using insights from several empirical observations. In particular, we propose a learning problem \cref{eq:pwd-learning-problem} to train invertible neural networks, which is then tested as a proof of concept through the above data analysis tasks on 2-dimensional toy data, see \cref{fig:swirl-example}.

\vspace{0.3cm}
Within the broader scope of data processing, we believe that this mathematical framework has important implications on how to construct Riemannian geometry and has important implications for handling classically Euclidean data in general.

\subsection{Outline}
\label{sec:intro-pullback-geometry-outline} 
This article is structured as follows. \Cref{sec:prelim-pullback-geometry} covers basic notation from differential and Riemannian geometry and states several known basic results for $\Real^\dimInd$ under a pullback Riemannian geometry. In \cref{sec:basic-processing-pullback} we consider well-posedness and stability of interpolation and Riemannian barycentres under pullback geometry. \Cref{sec:non-linear-compression-pullback} focuses on leveraging ideas from data compression on symmetric Riemannian manifolds to the pullback geometry setting and subsequently on constructing the (non-linear) RAE and CC-RAE mappings. \Cref{sec:learning-diffeo} proposes a deep learning-based approach to constructing diffeomorphisms, which covers all the best practices from the previous sections. In \cref{sec:numerics-pullback-manifolds} we visualize the predicted instabilities we get from curvature and unsuitable diffeomorphisms and test the practical ideas developed in this work. Finally, we summarize our findings in \cref{sec:conclusions-pullback-manifolds}.



\section{Preliminaries}
\label{sec:prelim-pullback-geometry}
\subsection{Notation}

Here we present some basic notations from differential and Riemannian geometry, see \cite{boothby2003introduction,carmo1992riemannian,lee2013smooth,sakai1996riemannian} for details. 

Let $\manifold$ be a smooth manifold. We write $C^\infty(\manifold)$ for the space of smooth functions over $\manifold$. The \emph{tangent space} at $\mPoint \in \manifold$, which is defined as the space of all \emph{derivations} at $\mPoint$, is denoted by $\tangent_\mPoint \manifold$ and for \emph{tangent vectors} we write $\tangentVector_\mPoint \in \tangent_\mPoint \manifold$. For the \emph{tangent bundle} we write $\tangent\manifold$ and smooth vector fields, which are defined as \emph{smooth sections} of the tangent bundle, are written as $\vectorfield(\manifold) \subset \tangent\manifold$.

A smooth manifold $\manifold$ becomes a \emph{Riemannian manifold} if it is equipped with a smoothly varying \emph{metric tensor field} $(\,\cdot\,, \,\cdot\,) \colon \vectorfield(\manifold) \times \vectorfield(\manifold) \to C^\infty(\manifold)$. This tensor field induces a \emph{(Riemannian) metric} $\distance_{\manifold} \colon \manifold\times\manifold\to\Real$. The metric tensor can also be used to construct a unique affine connection, the \emph{Levi-Civita connection}, that is denoted by $\nabla_{(\,\cdot\,)}(\,\cdot\,) : \vectorfield(\manifold) \times \vectorfield(\manifold) \to \vectorfield(\manifold)$. 
This connection is in turn the cornerstone of a myriad of manifold mappings.
One is the notion of a \emph{geodesic}, which for two points $\mPoint,\mPointB \in \manifold$ is defined as a curve $\geodesic_{\mPoint,\mPointB} \colon [0,1] \to \manifold$ with minimal length that connects $\mPoint$ with $\mPointB$. A subset $\mathcal{C}\subset \manifold$ is called \emph{geodesically convex} if,  given any two points $\mPoint, \mPointB \in\mathcal{C}$, there is a unique geodesic that connects $\mPoint$ with $\mPointB$ contained in $\mathcal{C}$. Somewhat related are geodesically convex functions $f:\manifold\to \Real$, which are functions such that $t \mapsto f \circ \geodesic_{\mPoint,\mPointB} (t)$ is convex for every geodesic $\geodesic_{\mPoint,\mPointB}$.
Another closely related notion to geodesics is the curve $t \mapsto \geodesic_{\mPoint,\tangentVector_\mPoint}(t)$  for a geodesic starting from $\mPoint\in\manifold$ with velocity $\dot{\geodesic}_{\mPoint,\tangentVector_\mPoint} (0) = \tangentVector_\mPoint \in \tangent_\mPoint\manifold$. This can be used to define the \emph{exponential map} $\exp_\mPoint \colon \mathcal{D}_\mPoint \to \manifold$ as 
\begin{equation}
\exp_\mPoint(\tangentVector_\mPoint) := \geodesic_{\mPoint,\tangentVector_\mPoint}(1)
\quad\text{where $\mathcal{D}_\mPoint \subset \tangent_\mPoint\manifold$ is the set on which $\geodesic_{\mPoint,\tangentVector_\mPoint}(1)$ is defined.} 
\end{equation}
The manifold $\manifold$ is said to be \emph{(geodesically) complete} whenever $\mathcal{D}_\mPoint = \mathcal{T}_{\mPoint}\manifold$ for all $\mPoint \in \manifold$.
Furthermore, the \emph{logarithmic map} $\log_\mPoint \colon \exp(\mathcal{D}'_\mPoint ) \to \mathcal{D}'_\mPoint$ is defined as the inverse of $\exp_\mPoint$, so it is well-defined on  $\mathcal{D}'_{\mPoint} \subset \mathcal{D}_{\mPoint}$ where $\exp_\mPoint$ is a diffeomorphism. 
Moreover, the \emph{Riemannian gradient} of a smooth function $f \colon \manifold\to \Real$ denotes the unique vector field $\Grad f \in \vectorfield(\manifold)$ such that 
\begin{equation}
    (\Grad f,
      \tangentVector )_\mPoint :=  \tangentVector_\mPoint f := D_{\mPoint} f(\cdot) [\tangentVector_\mPoint] ,
    \quad\text{holds for any $\tangentVector \in \vectorfield(\manifold)$ and $\mPoint\in\manifold$,}
\end{equation}
where $D_{\mPoint}f(\cdot): \tangent_\mPoint \manifold \to \Real$ denotes the differential of $f$. 
Finally, we write $\curvature (\,\cdot\,, \,\cdot\,)(\,\cdot\,) \colon \vectorfield(\manifold) \times \vectorfield(\manifold) \times \vectorfield(\manifold) \to \vectorfield(\manifold)$ for the \emph{curvature operator}, which can be used to define the sectional curvature $\sectionalcurvature: \vectorfield(\manifold) \times \vectorfield(\manifold) \to C^\infty(\manifold)$ as
\begin{equation}
    \sectionalcurvature( \tangentVector,\tangentVectorB) :=\left\{\begin{matrix}
    \frac{(\curvature(\tangentVector, \tangentVectorB) \tangentVectorB, \tangentVector)_{(\cdot)}}{\|\tangentVector\|_{(\cdot)}^2 \|\tangentVectorB\|_{(\cdot)}^2 - ( \tangentVector, \tangentVectorB)_{(\cdot)}^2} & \text{if  $\tangentVector, \tangentVectorB \in \vectorfield(\manifold)$ are linearly independent},  \\
    0 &  \text{otherwise}.\\
    \end{matrix}\right.
\end{equation}

Beyond basic concepts from Riemannian geometry, we will need several additional notions. First, let $(\manifold, (\cdot,\cdot))$ and $(\manifoldB, (\cdot,\cdot)')$ be Riemannian manifolds. A mapping $F:\manifold \to \manifoldB$ has \emph{local Lipschitz constant} $\operatorname{Lip}_{\mPoint}(F)\geq 0$ on an open neighbourhood $\mathcal{U}(\mPoint)\subset\manifold$ of $\mPoint\in \manifold$ if
\begin{equation}
    \distance_{\manifoldB}(F(\mPointB), F(\mPoint)) \leq \operatorname{Lip}_{\mPoint}(F) \; \distance_{\manifold}(\mPointB, \mPoint), \quad \text{for all } \mPointB \in \mathcal{U}(\mPoint).
\end{equation}
As a special case, consider that if $F$ is a local isometry on $\mathcal{U}(\mPoint)$, i.e., 
\begin{equation}
    \distance_{\manifoldB}(F(\mPointB), F(\mPoint)) = \distance_{\manifold}(\mPointB, \mPoint), \quad \text{for all } \mPointB \in \mathcal{U}(\mPoint),
\end{equation}
it has unit local Lipschitz constant at $\mPoint$.

Next, let $\mathcal{U}(\mPoint)\subset \manifold$ again be an open neighbourhood of $\mPoint\in \manifold$. A mapping $\reflection_{\mPoint}: \mathcal{U}(\mPoint) \to \manifold$ is called a \emph{geodesic reflection} at $\mPoint$ if
\begin{equation}
    \reflection_{\mPoint}(\mPoint)=\mPoint \text { and } D_{\mPoint} \reflection_{\mPoint}(\cdot) =-\operatorname{id}_{\mPoint}.
    \label{eq:geodesic-reflection-def-properties}
\end{equation}
A Riemannian manifold $\manifold$ is called \emph{(locally) symmetric}, if there exists a geodesic reflection at any point $\mPoint \in \manifold$ that is an isometry on a local neighbourhood of $\mPoint$, i.e., there exists a geodesic reflection $\reflection_{\mPoint}: \mathcal{U}(\mPoint) \to \manifold$ and a neighbourhood $\mathcal{V}(\mPoint)\subset \mathcal{U}(\mPoint)$ such that for all $\mPointB, \mPointC \in \mathcal{V}(\mPoint)$ we have
\begin{equation}
    \distance_{\manifold}(\reflection_{\mPoint}(\mPointB), \reflection_{\mPoint}(\mPointC))=\distance_{\manifold}(\mPointB, \mPointC) .
    \label{eq:symmetry-properties}
\end{equation}
A symmetric space is said to be a \emph{globally symmetric} space if in addition its geodesic symmetries can be extended to isometries on all of $\manifold$.

Finally, if $(\manifold, (\cdot,\cdot))$ is a $\dimInd$-dimensional Riemannian manifold, $\manifoldB$ is a $\dimInd$-dimensional smooth manifold and $\diffeo:\manifoldB \to \manifold$ is a diffeomorphism, the \emph{pullback metric}
\begin{equation}
    (\tangentVector, \tangentVectorB)^\diffeo := (\diffeo_*[\tangentVector], \diffeo_*[\tangentVectorB]) := (D_{(\cdot)}\diffeo[\tangentVector_{(\cdot)}], D_{(\cdot)}\diffeo[\tangentVectorB_{(\cdot)}])_{\diffeo(\cdot)}
\end{equation}
defines a Riemannian structure on $\manifoldB$, which we denote by $(\manifoldB, (\cdot,\cdot)^\diffeo)$. The mapping $\diffeo_* : \vectorfield(\manifoldB) \to \vectorfield(\manifold)$ denotes the \emph{pushforward} of $\diffeo$.



\subsection{Riemannian geometry on $(\Real^\dimInd, (\cdot,\cdot)^\diffeo)$}
\label{sec:prelim-pullback-geometry-theoretical-results}
In the following, we will state several \emph{known results} for Riemannian geometry on $(\Real^\dimInd, (\cdot,\cdot)^\diffeo)$, where the mapping $\diffeo:\Real^\dimInd \to \manifold$ is a diffeomorphism and $(\manifold, (\cdot,\cdot))$ is a $\dimInd$-dimensional Riemannian manifold. In particular, we will remind the reader that we can pull back basic mappings and additional structure, e.g., symmetry. For the sake of completeness, the proofs are also provided in \cref{app:proofs-prelim-pullback-geometry-theoretical-results}.


We outlined in \cref{sec:intro-pull-back-geometry} that we want to be able to interpolate, extrapolate, compute non-linear means, and do data decomposition on the Riemannian manifold $(\Real^\dimInd, (\cdot,\cdot)^\diffeo)$. For that, access to basic manifold mappings -- geodesics, exponential mapping, logarithmic mapping, distance, parallel transport -- is required to carry out these tasks. \Cref{sec:intro-pullback-geometry-related-work} hinted that these mappings are accessible in closed-form under pullback geometry. More concretely, the result below states that we can map points in $\Real^\dimInd$ and tangent vectors at these points to the manifold $\manifold$ and tangent spaces thereon, use the corresponding mapping of interest there and map the result back. In particular, if we have closed-form expressions for manifold mappings on $(\manifold, (\cdot,\cdot))$, we get closed-form expressions on $(\Real^\dimInd, (\cdot,\cdot)^\diffeo)$.

\begin{proposition}
\label{thm:pull-back-mappings}
    Let $(\manifold, (\cdot,\cdot))$ be a $\dimInd$-dimensional Riemannian manifold and let $\diffeo:\Real^\dimInd \to \manifold$ be a smooth diffeomorphism such that $\diffeo(\Real^\dimInd) \subset \manifold$ is a geodesically convex set.

    Then, 
    \begin{enumerate}[label=(\roman*)]
        \item length-minimising geodesics $\geodesic^\diffeo_{\Vector, \VectorB}:[0,1] \to \Real^\dimInd$ on $(\Real^\dimInd, (\cdot,\cdot)^\diffeo)$ are given by 
        \begin{equation}
            \geodesic^\diffeo_{\Vector, \VectorB}(t) = \diffeo^{-1}(\geodesic_{\diffeo(\Vector), \diffeo(\VectorB)}(t))
            \label{eq:thm-geodesic-remetrized}
        \end{equation}
        where $\geodesic_{\mPoint, \mPointB}:[0,1]\to \manifold$ are length-minimising geodesics on the manifold $\manifold$ generated by $(\cdot,\cdot)$.
        \item the logarithmic map $\log^\diffeo_{\Vector} (\cdot):\Real^\dimInd \to \tangent_\Vector \Real^\dimInd$  on $(\Real^\dimInd, (\cdot,\cdot)^\diffeo)$ is given by 
        \begin{equation}
            \log^\diffeo_{\Vector} \VectorB = \diffeo^{-1}_{*}[\log_{\diffeo(\Vector)} \diffeo(\VectorB)]
            \label{eq:thm-log-remetrized}
        \end{equation}
        where $\log_{\mPoint} (\cdot):\manifold \to \tangent_\mPoint \manifold$ is the logarithmic map on the manifold $\manifold$ generated by $(\cdot,\cdot)$.
        \item the exponential map $\exp^\diffeo_{\Vector} (\cdot):\mathcal{G}_{\Vector}\to \Real^\dimInd$ for $\mathcal{G}_{\Vector} := \log^\diffeo_{\Vector} (\Real^\dimInd) \subset \tangent_\Vector \Real^\dimInd$ on $(\Real^\dimInd, (\cdot,\cdot)^\diffeo)$ is given by 
        \begin{equation}
             \exp^\diffeo_\Vector (\tangentVector_\Vector) = \diffeo^{-1}(\exp_{\diffeo(\Vector)}(\diffeo_*[\tangentVector_\Vector]))
             \label{eq:thm-exp-remetrized}
        \end{equation}
        where $\exp_{\mPoint} (\cdot):\tangent_\mPoint \manifold \to \manifold$ is the exponential map on the manifold $\manifold$ generated by $(\cdot,\cdot)$.
        \item the distance $\distance^\diffeo_{\Real^\dimInd}:\Real^\dimInd \times\Real^\dimInd \to \Real$ on $(\Real^\dimInd, (\cdot,\cdot)^\diffeo)$ is given by 
        \begin{equation}
            \distance^\diffeo_{\Real^\dimInd}(\Vector, \VectorB) = \distance_{\manifold}(\diffeo(\Vector), \diffeo(\VectorB)),
            \label{eq:thm-distance-remetrized}
        \end{equation}
        where $\distance_{\manifold}:\manifold\times \manifold \to \Real$ is the distance on the manifold $\manifold$ generated by $(\cdot,\cdot)$.
        \item parallel transport along geodesics $\mathcal{P}^\diffeo_{\VectorB \leftarrow \Vector} :\tangent_\Vector \Real^\dimInd \to \tangent_\VectorB \Real^\dimInd$ on $(\Real^\dimInd, (\cdot,\cdot)^\diffeo)$ is given by 
        \begin{equation}
            \mathcal{P}^\diffeo_{\VectorB \leftarrow \Vector} \tangentVector_\Vector = \diffeo^{-1}_* [\mathcal{P}_{\diffeo(\VectorB) \leftarrow \diffeo(\Vector)}  \diffeo_*[\tangentVector_\Vector]]
            \label{eq:thm-parallel-transport-remetrized}
        \end{equation}
        where $\mathcal{P}_{\mPointB \leftarrow \mPoint} :\tangent_\mPoint \manifold \to \tangent_\mPointB \manifold$ is parallel transport on the manifold $\manifold$ generated by $(\cdot,\cdot)$.
    \end{enumerate}
\end{proposition}

The above result gives well-posedness of almost all mappings directly, except for the exponential mapping. For that we need an an extra condition on both the manifold $\manifold$ and the diffeomorpism $\diffeo$.

\begin{proposition}
\label{thm:completeness}
    Let $(\manifold, (\cdot,\cdot))$ be a $\dimInd$-dimensional geodesically convex and complete Riemannian manifold and let $\diffeo:\Real^\dimInd \to \manifold$ be a smooth global diffeomorphism, i.e., $\diffeo(\Real^\dimInd) = \manifold$. 
    
    Then, $(\Real^\dimInd,  (\cdot,\cdot)^\diffeo)$ is geodesically complete and $\exp^\diffeo_{\Vector}$ is well-defined on $\tangent_{\Vector}\Real^\dimInd$ for any $\Vector\in \Real^\dimInd$.
\end{proposition}




Next, we also indicated in \cref{sec:intro-pullback-geometry-related-work} that symmetry of the Riemannian structure is inherited under pullback geometry. This is made more concrete in the following result.

\begin{proposition}
\label{thm:local-symmetry}
Let $(\manifold, (\cdot,\cdot))$ be a $\dimInd$-dimensional (locally) symmetric Riemannian manifold and let $\diffeo:\Real^\dimInd \to \manifold$ be a smooth diffeomorphism such that $\diffeo(\Real^\dimInd) \subset \manifold$ is a geodesically convex set.

Then, $(\Real^\dimInd, (\cdot,\cdot)^\diffeo)$ is a symmetric Riemannian manifold.
\end{proposition}

The extension to global symmetry also holds under slightly stronger assumptions.

\begin{corollary}
    Let $(\manifold, (\cdot,\cdot))$ be a $\dimInd$-dimensional globally symmetric and geodesically convex Riemannian manifold and let $\diffeo:\Real^\dimInd \to \manifold$ be a smooth global diffeomorphism, i.e., $\diffeo(\Real^\dimInd) = \manifold$.
    
    
    Then $(\Real^\dimInd, (\cdot,\cdot)^\diffeo)$ is a globally symmetric Riemannian manifold.
\end{corollary}

Finally, although somewhat beyond the scope of this work, it is good to mention to the interested reader that if the Riemannian manifold $(\manifold, (\cdot,\cdot))$ is \emph{Hadamard} --  complete, simply connected with non-positive sectional curvature --, this is also inherited under pullback geometry.
\begin{proposition}
\label{thm:hadamard}
    Let $(\manifold, (\cdot,\cdot))$ be a $\dimInd$-dimensional Hadamard manifold and let $\diffeo:\Real^\dimInd \to \manifold$ be a smooth global diffeomorphism, i.e., $\diffeo(\Real^\dimInd) = \manifold$. 
    
    Then, the Riemannian manifold $(\Real^\dimInd,  (\cdot,\cdot)^\diffeo)$ is a Hadamard manifold.
\end{proposition}



\section{Basic data processing on $(\Real^\dimInd, (\cdot,\cdot)^\diffeo)$ }
\label{sec:basic-processing-pullback}

Throughout the rest of the paper we focus on how to construct suitable diffeomorphisms $\diffeo:\Real^\dimInd \to \manifold$ into some Riemannian manifold $(\manifold, (\cdot, \cdot))$ and how to modify existing algorithms for symmetric Riemannian manifolds to $(\Real^\dimInd, (\cdot,\cdot)^\diffeo)$. In this section we only consider basic tasks such as interpolation and computing barycentres. As there are standard definitions of these tasks, we will not have to modify anything algorithmically once we have chosen a Riemannian manifold and a diffeomorphism. Regarding construction, there are arguably better and worse choices for $\diffeo$ and this section will address some basic requirements the construction should follow.

In particular, we will see that diffeomorphisms should map the data manifold into low-dimensional geodesic subspaces of $(\manifold, (\cdot, \cdot))$ for proper interpolation through and barycentres within the data set, but should also do it in such a way that it is a local isometry -- with respect to the original ambient Riemannian structure $(\Real^\dimInd, (\cdot,\cdot)_2)$ -- in order to get $\ell^2$-stability, which is a feature we ideally inherit from data analysis in the Euclidean setting in case our data is somewhat noisy.

\subsection{Diffeomorphisms for proper and stable interpolation}
In the following we will consider geodesic interpolation $t \mapsto \geodesic^\diffeo_{\Vector,\VectorB}(t)$ between points $\Vector,\VectorB\in \Real^\dimInd$ on $(\Real^\dimInd, (\cdot,\cdot)^\diffeo)$. The above-mentioned requirement that diffeomorphisms should map the data manifold into a low-dimensional geodesic subspace for proper interpolation can be understood from \cref{thm:pull-back-mappings} (i) directly. One might wonder next where the rest of the ambient space should be mapped to. As indicated above, this is a particularly important question if we have slightly noisy data that does not lie perfectly on the low-dimensional data manifold. This question boils down to stability, which can -- to a certain extent -- be guaranteed through local isometry as we will see in the following. 


Our goal in this section is to characterize stability of geodesics in $\Vector$ and $\VectorB$ and we will see that $\ell^2$-stability is fully determined by the diffeomorphism's deviation from isometry around the data manifold and the curvature of $(\manifold, (\cdot, \cdot))$. Our strategy for characterizing stability follows four main steps. For the first three steps we consider $\distance_\manifold$-stability on $\manifold$ and only at the final step, we pull the results to $\Real^\dimInd$, which gives us our first main result (\cref{thm:stability-geodesics-symmetric}).

For the first step we consider $\distance_\manifold$-stability of \emph{geodesic variations}, i.e., smooth mappings $\Gamma: [0,1]\times[-\epsilon, \epsilon] \to \manifold$ for $\epsilon>0$ such that the curve $t \mapsto \Gamma(t,0)$ is a geodesic. The following result -- that is an adaptation of \cite[Lemma 3.1]{diepeveen2023curvature} -- tells us that stability is governed by first-order behaviour of $\Gamma$ in the direction of the variation.

\begin{lemma}
\label{eq:geodesic-variation-taylor-second-order}
    Let $(\manifold, (\cdot, \cdot))$ be a Riemannian manifold and $\Gamma: [0,1]\times[-\epsilon, \epsilon] \to \manifold$ be a geodesic variation.

    Then,
    \begin{equation}
    \distance_{\manifold}(\Gamma(t,\epsilon), \Gamma(t,0))^2 =  \epsilon^2 \Bigl\|\frac{\partial}{\partial s} \Gamma(t,s)\Bigr\|^2_{\Gamma(t,0)} + \mathcal{O}(\epsilon^3).
    \label{eq:error-geodesic-variation-approx}
\end{equation}
    
\end{lemma}

\begin{proof}
    We will show that \cref{eq:error-geodesic-variation-approx} holds through Taylor expansion of $s \mapsto \distance_{\manifold}(\Gamma(t,s), \Gamma(t,0))^2$ for a fixed $t$ around $s=0$. Expanding up to second order gives
\begin{multline}
    \distance_{\manifold}(\gamma_{\mPoint,\mPointC} (t),\gamma_{\mPoint,\mPointB}(t))^2 = \distance_{\manifold}(\Gamma(t,\epsilon), \Gamma(t,0))^2 
    = \distance_{\manifold}(\Gamma(t,0), \Gamma(t,0))^2 +  \epsilon \frac{\mathrm{d}}{\mathrm{d} s} \distance_{\manifold}(\Gamma(t,s), \Gamma(t,0))^2\mid_{s=0} \\
    + \frac{\epsilon^2}{2} \frac{\mathrm{d}^2}{\mathrm{d} s^2} \distance_{\manifold}(\Gamma(t,s), \Gamma(t,0))^2\mid_{s=0} + \frac{1}{2}\int_0^\epsilon \frac{\mathrm{d}^3}{\mathrm{d} s'^3} \distance_{\manifold}(\Gamma(t,s'), \Gamma(t,0))^2 (s - s')^2\mathrm{d}s'.
\end{multline}
For proving the result \cref{eq:error-geodesic-variation-approx}, we must show that the zeroth and first order terms are zero, i.e.,
\begin{equation}
    \distance_{\manifold}(\Gamma(t,0), \Gamma(t,0))^2=0,
    \label{eq:jacobi-bound-low-rank-mfld-tensor-thm-taylor0}
\end{equation}
and
\begin{equation}
    \frac{\mathrm{d}}{\mathrm{d} s} \distance_{\manifold}(\Gamma(t,s), \Gamma(t,0))^2\mid_{s=0} = 0,
    \label{eq:jacobi-bound-low-rank-mfld-tensor-thm-taylor1}
\end{equation}
and that the second order term satisfies
\begin{equation}
    \frac{1}{2} \frac{\mathrm{d}^2}{\mathrm{d} s^2} \distance_{\manifold}(\Gamma(t,s), \Gamma(t,0))^2\mid_{s=0} = \Bigl\|\frac{\partial}{\partial s} \Gamma(t,s)\Bigr\|^2_{\Gamma(t,0)}.
    \label{eq:jacobi-bound-low-rank-mfld-tensor-thm-taylor2}
\end{equation}
It is clear that the remainder term is $\mathcal{O}(\epsilon^3)$.

Trivially, \cref{eq:jacobi-bound-low-rank-mfld-tensor-thm-taylor0} holds, because $\distance_{\manifold}$ is a distance. For showing \cref{eq:jacobi-bound-low-rank-mfld-tensor-thm-taylor1} notice that
\begin{multline}
    \frac{\mathrm{d}}{\mathrm{d} s} \distance_{\manifold}(\Gamma(t,s), \Gamma(t,0))^2 = D_{\Gamma(t,s)} \distance_{\manifold}(\cdot, \Gamma(t,0))^2 \Bigl[ \frac{\partial}{\partial s} \Gamma(t,s) \Bigr] = \Bigl( \operatorname{grad} \distance_{\manifold}(\cdot, \Gamma(t,0))^2 \mid_{\Gamma(t,s)}, \frac{\partial}{\partial s} \Gamma(t,s)\Bigr)_{\Gamma(t,s)} \\
    = \Bigl(- 2 \log_{\Gamma(t,s)} \Gamma(t,0), \frac{\partial}{\partial s} \Gamma(t,s)\Bigr)_{\Gamma(t,s)},
    \label{eq:towards-jacobi-bound-low-rank-mfld-tensor-thm-taylor1}
\end{multline}
where we used the chain rule in the first equality, the definition of the Riemannian gradient in the second equality and the fact that $\operatorname{grad} \distance_{\manifold}(\cdot, \mathbf{a})^2\mid_{\mathbf{b}}= -2 \log_\mathbf{b} \mathbf{a} \in \tangent_\mathbf{b} \manifold$ for any $\mathbf{a},\mathbf{b}\in \manifold$. Then, evaluating \cref{eq:towards-jacobi-bound-low-rank-mfld-tensor-thm-taylor1} at $s=0$ yields \cref{eq:jacobi-bound-low-rank-mfld-tensor-thm-taylor1} as $\log_{\Gamma(t,0)} \Gamma(t,0) = 0_{\Gamma(t,0)}\in \tangent_{\Gamma(t,0)}\manifold$.

Finally, it remains to show that \cref{eq:jacobi-bound-low-rank-mfld-tensor-thm-taylor2} holds. From \cref{eq:towards-jacobi-bound-low-rank-mfld-tensor-thm-taylor1} we find that
\begin{multline}
    \frac{\mathrm{d}^2}{\mathrm{d} s^2} \distance_{\manifold}(\Gamma(t,s), \Gamma(t,0))^2 = \frac{\mathrm{d} }{\mathrm{d} s} \Bigl(- 2 \log_{\Gamma(t,s)} \Gamma(t,0), \frac{\partial }{\partial s} \Gamma(t,s)\Bigr)_{\Gamma(t,s)} \\
    = \frac{\mathrm{d} }{\mathrm{d} s} \Bigl(- 2 \log_{(\cdot)} \Gamma(t,0), \frac{\partial }{\partial s'} \Gamma(t,s') \mid_{s' = (\Gamma(1,\star))^{-1}(\cdot) }\Bigr)_{\Gamma(t,s)} \\
    = \frac{\partial }{\partial s} \Gamma(t,s) \Bigl(- 2 \log_{(\cdot)} \Gamma(t,0), \frac{\partial }{\partial s'} \Gamma(t,s') \mid_{s' = (\Gamma(1,\star))^{-1}(\cdot) } \Bigr)_{(\cdot)} \\
    \overset{\text{metric compatibility}}{=} \Bigl(- 2 \nabla_{\frac{\partial }{\partial s} \Gamma(t,s)}\log_{(\cdot)} \Gamma(t,0), \frac{\partial }{\partial s} \Gamma(t,s) \Bigr)_{\Gamma(t,s)} \\ + \Bigl(- 2 \log_{\Gamma(t,s)} \Gamma(t,0), \nabla_{\frac{\partial }{\partial s} \Gamma(t,s)} \frac{\partial }{\partial s'} \Gamma(t,s') \mid_{s' = (\Gamma(1,\star))^{-1}(\cdot) }\Bigr)_{\Gamma(t,s)}.
    \label{eq:towards-jacobi-bound-low-rank-mfld-tensor-thm-taylor2}
\end{multline}

Evaluating \cref{eq:towards-jacobi-bound-low-rank-mfld-tensor-thm-taylor2} at $s=0$ gives once again zero for the second term because $\log_{\Gamma(t,0)} \Gamma(t,0) = 0_{\Gamma(t,0)}\in \tangent_{\Gamma(t,0)}\manifold$. However, using that $(-\nabla_{X} \log_{(\cdot)}\mathbf{a}, Y)_{\mathbf{a}} = (X,Y)_{\mathbf{a}}$ for any $\mathbf{a}\in \manifold$, the first term in the final line of \cref{eq:towards-jacobi-bound-low-rank-mfld-tensor-thm-taylor2} reduces to
\begin{multline}
    \Bigl(- 2 \nabla_{\frac{\partial }{\partial s} \Gamma(t,s)}\log_{(\cdot)} \Gamma(t,0), \frac{\partial }{\partial s} \Gamma(t,s) \Bigr)_{\Gamma(t,s)}\mid_{s=0} = 2 \Bigl(\frac{\partial }{\partial s} \Gamma(t,s)\mid_{s=0} , \frac{\partial }{\partial s} \Gamma(t,s) \mid_{s=0}\Bigr)_{\Gamma(t,0)} \\
    = 2 \Bigl\|\frac{\partial}{\partial s} \Gamma(t,s)\Bigr\|^2_{\Gamma(t,0)}.
\end{multline}
When the above is substituted back into \cref{eq:towards-jacobi-bound-low-rank-mfld-tensor-thm-taylor2}, \cref{eq:jacobi-bound-low-rank-mfld-tensor-thm-taylor2} holds, which proves the claim.
\end{proof}

Next, we can use \cref{eq:geodesic-variation-taylor-second-order} to get expressions for stability of geodesics at their beginning and end points.

\begin{lemma}
\label{lem:evaluated-variation-geodesics}
    Let $(\manifold, (\cdot, \cdot))$ be a Riemannian manifold and $\mPoint,\mPointB \in \manifold$ be distinct points on the manifold.

Then,
\begin{enumerate}[label=(\roman*)]
    \item as $\mPointC\to \mPoint$ the mapping $t \mapsto \distance_{\manifold}(\gamma_{\mPointC,\mPointB} (t),\gamma_{\mPoint,\mPointB}(t))^2$ behaves as
    \begin{equation}
    \distance_{\manifold}(\gamma_{\mPointC,\mPointB} (t),\gamma_{\mPoint,\mPointB}(t))^2 =  \distance_{\manifold}(\mPointC,\mPoint)^2 \Bigl\|D_{\mPoint} \gamma_{(\cdot),\mPointB}(t) \Bigl[\frac{\log_{\mPoint} \mPointC}{\|\log_{\mPoint} \mPointC\|_{\mPoint}}\Bigr] \Bigr\|^2_{\gamma_{\mPoint,\mPointB}(t)} + \mathcal{O}(\distance_{\manifold}(\mPointC,\mPoint)^3).
    \label{eq:jacobi-bound-variation-i}
\end{equation}
\item as $\mPointC\to \mPointB$ the mapping $t \mapsto \distance_{\manifold}(\gamma_{\mPoint,\mPointC} (t),\gamma_{\mPoint,\mPointB}(t))^2$ behaves as
\begin{equation}
    \distance_{\manifold}(\gamma_{\mPoint,\mPointC} (t),\gamma_{\mPoint,\mPointB}(t))^2 =  \distance_{\manifold}(\mPointC,\mPointB)^2 \Bigl\|D_{\mPointB} \gamma_{\mPoint,(\cdot)}(t) \Bigl[\frac{\log_{\mPointB} \mPointC}{\|\log_{\mPointB} \mPointC\|_{\mPointB}}\Bigr] \Bigr\|^2_{\gamma_{\mPoint,\mPointB}(t)} + \mathcal{O}(\distance_{\manifold}(\mPointC,\mPointB)^3).
    \label{eq:jacobi-bound-variation-ii}
\end{equation}
\end{enumerate}

\end{lemma}

\begin{proof}
(i) For showing \cref{eq:jacobi-bound-variation-i} define the geodesic variation $\Gamma_0:[0,1]\times[-\epsilon, \epsilon] \to \manifold$ with $\epsilon := \distance_{\manifold}(\mPointC,\mPoint)$ by
\begin{equation}
    \Gamma_0(t,s):= \exp_{\exp_{\mPoint}(\frac{s}{\epsilon}\log_{\mPoint} \mPointC)}\Bigl(t \log_{\exp_{\mPoint}(\frac{s}{\epsilon}\log_{\mPoint} \mPointC)} \mPointB\Bigr).
\end{equation}
Note that $\Gamma_0(t,0) = \gamma_{\mPoint,\mPointB}(t)$ and $\Gamma_0(t,\epsilon) = \gamma_{\mPointC,\mPointB} (t)$, and that
\begin{equation}
    \frac{\partial}{\partial s}\Gamma_0(t,s)\mid_{s=0} = D_{\mPoint} \gamma_{(\cdot),\mPointB}(t) \Bigl[\frac{1}{\epsilon} \log_{\mPoint} \mPointC \Bigr] = D_{\mPoint} \gamma_{(\cdot),\mPointB}(t) \Bigl[\frac{\log_{\mPoint} \mPointC}{\distance_{\manifold}(\mPointC,\mPoint)}\Bigr]  = D_{\mPointB} \gamma_{(\cdot),\mPointB}(t) \Bigl[\frac{\log_{\mPoint} \mPointC}{\|\log_{\mPoint} \mPointC\|_{\mPoint}}\Bigr].
\end{equation}
Invoking \cref{eq:error-geodesic-variation-approx} from \cref{eq:geodesic-variation-taylor-second-order} gives the desired result.

(ii) For showing \cref{eq:jacobi-bound-variation-ii} define the geodesic variation $\Gamma_1:[0,1]\times[-\epsilon, \epsilon] \to \manifold$ with $\epsilon := \distance_{\manifold}(\mPointC,\mPointB)$ by
\begin{equation}
    \Gamma_1(t,s):= \exp_{\exp_{\mPointB}(\frac{s}{\epsilon}\log_{\mPointB} \mPointC)}\Bigl((1-t) \log_{\exp_{\mPointB}(\frac{s}{\epsilon}\log_{\mPointB} \mPointC)} \mPoint\Bigr).
\end{equation}
Note that $\Gamma_1(t,0) = \gamma_{\mPoint,\mPointB}(t)$ and $\Gamma_1(t,\epsilon) = \gamma_{\mPoint,\mPointC} (t)$, and that
\begin{equation}
    \frac{\partial}{\partial s}\Gamma_1(t,s)\mid_{s=0} = D_{\mPointB} \gamma_{\mPoint,(\cdot)}(t) \Bigl[\frac{1}{\epsilon} \log_{\mPointB} \mPointC \Bigr] = D_{\mPointB} \gamma_{\mPoint,(\cdot)}(t) \Bigl[\frac{\log_{\mPointB} \mPointC}{\distance_{\manifold}(\mPointC,\mPointB)}\Bigr]  = D_{\mPointB} \gamma_{\mPoint,(\cdot)}(t) \Bigl[\frac{\log_{\mPointB} \mPointC}{\|\log_{\mPointB} \mPointC\|_{\mPointB}}\Bigr] 
\end{equation}
Invoking \cref{eq:error-geodesic-variation-approx} from \cref{eq:geodesic-variation-taylor-second-order} gives the desired result.

\end{proof}

In general it will be hard to evaluate the differentials in the expressions \cref{eq:jacobi-bound-variation-i,eq:jacobi-bound-variation-ii} in \cref{lem:evaluated-variation-geodesics}. However, if our Riemannian manifold is symmetric, we can get closed-form expressions.

\begin{lemma}
\label{lem:evaluated-variation-geodesics-symmetric}
    Let $(\manifold, (\cdot, \cdot))$ be a $\dimInd$-dimensional symmetric Riemannian manifold and $\mPoint,\mPointB \in \manifold$ be distinct points on the manifold. Furthermore, let $\{\Theta^{\sumIndB}_{\mPoint}\}_{\sumIndB=1}^{\dimInd}\subset \tangent_{\mPoint} \manifold$ and $\{\Theta^{\sumIndB}_{\mPointB}\}_{\sumIndB=1}^{\dimInd}\subset \tangent_{\mPointB} \manifold$ be orthonormal frames that diagonalize the operator 
\begin{equation}
    \Theta_{\mPoint} \mapsto \curvature_{\mPoint}(\Theta_\mPoint, \log_{\mPoint}\mPointB) \log_{\mPoint}\mPointB \quad \text{and} \quad \Theta_{\mPointB} \mapsto \curvature_{\mPointB}(\Theta_\mPointB, \log_{\mPointB}\mPoint) \log_{\mPointB}\mPoint,
\end{equation}
with respective eigenvalues $\lambda_{\sumIndB}$ and $\mu_{\sumIndB}$ for $\sumIndB=1, \ldots, \dimInd$ and define $\betaGeoEnd:\Real \times [0,1]\to\Real$ as
\begin{equation}
    \betaGeoEnd(\kappa, t) := \left\{\begin{matrix}
\frac{\sinh(\sqrt{-\kappa} t)}{\sinh(\sqrt{-\kappa})}, & \kappa <0, \\
t, & \kappa = 0, \\
\frac{\sin(\sqrt{\kappa} t)}{\sin(\sqrt{\kappa})}, & \kappa >0.
\label{eq:beta-geodesic}
\end{matrix}\right.
\end{equation}

Then, 
\begin{enumerate}[label=(\roman*)]
    \item as $\mPointC\to \mPoint$ the mapping $t \mapsto \distance_{\manifold}(\gamma_{\mPointC,\mPointB} (t),\gamma_{\mPoint,\mPointB}(t))^2$ behaves as
    \begin{equation}
    \distance_{\manifold}(\gamma_{\mPointC,\mPointB} (t),\gamma_{\mPoint,\mPointB}(t))^2 = \distance_{\manifold}(\mPointC,\mPoint)^2 \sum_{\sumIndB=1}^{\dimInd}  \betaGeoEnd(\lambda_{\sumIndB}, 1 - t)^2 \Bigl( \frac{\log_{\mPoint} \mPointC}{\|\log_{\mPoint} \mPointC\|_\mPoint},  \Theta^{\sumIndB}_{\mPoint} \Bigr)_{\mPoint}^2 + \mathcal{O}(\distance_{\manifold}(\mPointC,\mPoint)^3).
    \label{eq:jacobi-bound-variation-symmetric-space-i}
\end{equation}
\item as $\mPointC\to \mPointB$ the mapping $t \mapsto \distance_{\manifold}(\gamma_{\mPoint,\mPointC} (t),\gamma_{\mPoint,\mPointB}(t))^2$ behaves as
\begin{equation}
    \distance_{\manifold}(\gamma_{\mPoint,\mPointC} (t),\gamma_{\mPoint,\mPointB}(t))^2 = \distance_{\manifold}(\mPointC,\mPointB)^2 \sum_{\sumIndB=1}^{\dimInd}  \betaGeoEnd(\mu_{\sumIndB}, t)^2 \Bigl( \frac{\log_{\mPointB} \mPointC}{\|\log_{\mPointB} \mPointC\|_\mPointB},  \Theta^{\sumIndB}_{\mPointB} \Bigr)_{\mPoint}^2 + \mathcal{O}(\distance_{\manifold}(\mPointC,\mPointB)^3).
    \label{eq:jacobi-bound-variation-symmetric-space-ii}
\end{equation}
\end{enumerate}

\end{lemma}

\begin{proof}

(i) The equality \cref{eq:jacobi-bound-variation-symmetric-space-i} follows directly from \cref{lem:evaluated-variation-geodesics} by evaluating $D_{\mPoint} \gamma_{(\cdot),\mPointB}(t) \Bigl[\frac{\log_{\mPoint} \mPointC}{\|\log_{\mPoint} \mPointC\|_{\mPoint}}\Bigr]$ in \cref{eq:jacobi-bound-variation-i} using \cite[Lemma~1]{bergmann2019recent}
    
(ii) The equality \cref{eq:jacobi-bound-variation-symmetric-space-ii} follows directly from \cref{lem:evaluated-variation-geodesics} by evaluating $D_{\mPointB} \gamma_{\mPoint,(\cdot)}(t) \Bigl[\frac{\log_{\mPointB} \mPointC}{\|\log_{\mPointB} \mPointC\|_{\mPointB}}\Bigr]$ in \cref{eq:jacobi-bound-variation-ii} using \cite[Lemma~1]{bergmann2019recent}
\end{proof}

Finally, we can use the expressions \cref{eq:jacobi-bound-variation-symmetric-space-i,eq:jacobi-bound-variation-symmetric-space-ii} in \cref{lem:evaluated-variation-geodesics-symmetric} for our original goal of quantifying geodesic stability. The following theorem tells us that $\ell^2$-stability is fully determined by the $\diffeo$'s deviation from isometry around the data manifold and the curvature of $(\manifold, (\cdot, \cdot))$. In particular, positive curvature will inherently cause instabilities.

\begin{theorem}[Stability of geodesics]
\label{thm:stability-geodesics-symmetric}
    Let $(\manifold, (\cdot,\cdot))$ be a $\dimInd$-dimensional symmetric Riemannian manifold and let $\diffeo:\Real^\dimInd \to \manifold$ be a smooth diffeomorphism such that $\diffeo(\Real^\dimInd) \subset \manifold$ is a geodesically convex set. Furthermore, 
    let $\Vector,\VectorB \in \Real^\dimInd$ be distinct points, let $\{\lambda_{\sumIndB}\}_{\sumIndB=1}^\dimInd,\{\mu_{\sumIndB}\}_{\sumIndB=1}^\dimInd \subset \Real$ be the eigenvalues of the operators 
\begin{equation}
    \Theta_{\diffeo(\Vector)} \mapsto \curvature_{\diffeo(\Vector)}(\Theta_{\diffeo(\Vector)}, \log_{\diffeo(\Vector)}\diffeo(\VectorB)) \log_{\diffeo(\Vector)}\diffeo(\VectorB) \quad \text{and} \quad \Theta_{\diffeo(\VectorB)} \mapsto \curvature_{\diffeo(\VectorB)}(\Theta_\mPointB, \log_{\diffeo(\VectorB)}\diffeo(\Vector)) \log_{\diffeo(\VectorB)}\diffeo(\Vector),
    \label{eq:thm-geo-stab-curvature-op}
\end{equation}
and define $\betaGeoEnd:\Real \times [0,1]\to\Real$ as in \cref{eq:beta-geodesic}. Finally, consider open neighbourhoods $\mathcal{U}(\Vector), \mathcal{U}(\VectorB)\subset \Real^\dimInd$ 
on which $\diffeo$ has local Lipschitz constants $\operatorname{Lip}_\Vector(\diffeo), \operatorname{Lip}_\VectorB(\diffeo)$ and consider an open neighbourhood and $\mathcal{V}(\geodesic_{\diffeo(\Vector), \diffeo(\VectorB)}(t)) \subset \manifold$ on which $\diffeo^{-1}$ has local Lipschitz constant $\operatorname{Lip}_{\geodesic_{\diffeo(\Vector), \diffeo(\VectorB)}(t)}(\diffeo^{-1})$.
    
Then,
\begin{enumerate}[label=(\roman*)]
    \item as $\VectorC \to \Vector$
    \begin{equation}
    \|\geodesic^\diffeo_{\VectorC, \VectorB}(t) - \geodesic^\diffeo_{\Vector, \VectorB}(t)\|_2 \leq \operatorname{Lip}_{\geodesic_{\diffeo(\Vector), \diffeo(\VectorB)}(t)}(\diffeo^{-1}) \betaGeoEnd(\lambda_{\max}, 1-t) \operatorname{Lip}_\Vector(\diffeo)  \|\VectorC - \Vector\|_2 + o(\|\VectorC - \Vector\|_2).
    \label{eq:thm-stability-geodesics-symmetric-i}
\end{equation}
    \item as $\VectorC \to \VectorB$
    \begin{equation}
    \|\geodesic^\diffeo_{\Vector, \VectorC}(t) - \geodesic^\diffeo_{\Vector, \VectorB}(t)\|_2 \leq \operatorname{Lip}_{\geodesic_{\diffeo(\Vector), \diffeo(\VectorB)}(t)}(\diffeo^{-1})  \betaGeoEnd(\mu_{\max}, t) \operatorname{Lip}_\VectorB(\diffeo) \|\VectorC - \VectorB\|_2 + o(\|\VectorC - \VectorB\|_2).
    \label{eq:thm-stability-geodesics-symmetric-ii}
\end{equation}
\end{enumerate}

\end{theorem}

\begin{proof}
Let  $\{\Theta^{\sumIndB}_{\diffeo(\Vector)}\}_{\sumIndB=1}^{\dimInd}\subset \tangent_{\diffeo(\Vector)} \manifold$ and $\{\Theta^{\sumIndB}_{\diffeo(\VectorB)}\}_{\sumIndB=1}^{\dimInd}\subset \tangent_{\diffeo(\VectorB)} \manifold$ be orthonormal frames that diagonalize the operators in \cref{eq:thm-geo-stab-curvature-op} corresponding to the eigenvalues $\{\lambda_{\sumIndB}\}_{\sumIndB=1}^\dimInd$ and $\{\mu_{\sumIndB}\}_{\sumIndB=1}^\dimInd$.

(i) As $\VectorC \to \Vector$ we may assume that $\VectorC \in \mathcal{U}(\Vector)$ and $\geodesic_{\diffeo(\VectorC), \diffeo(\VectorB)}(t) \in \mathcal{V}(\geodesic_{\diffeo(\Vector), \diffeo(\VectorB)}(t))$. Subsequently,
\begin{multline}
    \|\geodesic^\diffeo_{\VectorC, \VectorB}(t) - \geodesic^\diffeo_{\Vector, \VectorB}(t)\|_2^2 \overset{\text{\cref{thm:pull-back-mappings} (i)}}{=} \|\diffeo^{-1}(\geodesic_{\diffeo(\VectorC), \diffeo(\VectorB)}(t)) - \diffeo^{-1}(\geodesic_{\diffeo(\Vector), \diffeo(\VectorB)}(t) )\|_2^2 \\
    \leq \operatorname{Lip}_{\geodesic_{\diffeo(\Vector), \diffeo(\VectorB)}(t)}(\diffeo^{-1})^2 \distance_{\manifold}(\geodesic_{\diffeo(\VectorC), \diffeo(\VectorB)}(t) ,  \geodesic_{\diffeo(\Vector), \diffeo(\VectorB)}(t) )^2\\
    \overset{\text{\cref{lem:evaluated-variation-geodesics-symmetric} (ii)}}{=} \operatorname{Lip}_{\geodesic_{\diffeo(\Vector), \diffeo(\VectorB)}(t)}(\diffeo^{-1})^2 \sum_{\sumIndB=1}^{\dimInd} \beta(\lambda_{\sumIndB}, t)^2 \Bigl( \frac{\log_{\diffeo(\Vector)} \diffeo(\VectorC)}{\|\log_{\diffeo(\Vector)} \diffeo(\VectorC)\|_{\diffeo(\Vector)}},  \Theta^{\sumIndB}_{\diffeo(\Vector)} \Bigr)_{\diffeo(\Vector)}^2 \distance_{\manifold}(\diffeo(\VectorC),\diffeo(\Vector))^2 \\
    \hfill + \mathcal{O}(\distance_{\manifold}(\diffeo(\VectorC),\diffeo(\Vector))^3)\\
    \leq \operatorname{Lip}_{\geodesic_{\diffeo(\Vector), \diffeo(\VectorB)}(t)}(\diffeo^{-1})^2 \beta(\lambda_{\max}, t)^2 \sum_{\sumIndB=1}^{\dimInd}  \Bigl( \frac{\log_{\diffeo(\Vector)} \diffeo(\VectorC)}{\|\log_{\diffeo(\Vector)} \diffeo(\VectorC)\|_{\diffeo(\Vector)}},  \Theta^{\sumIndB}_{\diffeo(\Vector)} \Bigr)_{\diffeo(\Vector)}^2 \distance_{\manifold}(\diffeo(\VectorC),\diffeo(\Vector))^2 + \mathcal{O}(\distance_{\manifold}(\diffeo(\VectorC),\diffeo(\Vector))^3)\\
    = \operatorname{Lip}_{\geodesic_{\diffeo(\Vector), \diffeo(\VectorB)}(t)}(\diffeo^{-1})^2 \beta(\lambda_{\max}, t)^2  \distance_{\manifold}(\diffeo(\VectorC),\diffeo(\Vector))^2 + \mathcal{O}(\distance_{\manifold}(\diffeo(\VectorC),\diffeo(\Vector))^3) \\
    \leq \operatorname{Lip}_{\geodesic_{\diffeo(\Vector), \diffeo(\VectorB)}(t)}(\diffeo^{-1})^2 \beta(\mu_{\max}, t)^2 \operatorname{Lip}_\Vector(\diffeo) \|\VectorC - \Vector\|_2^2 + \mathcal{O}(\|\VectorC - \Vector\|_2^3),
\end{multline}
which implies \cref{eq:thm-stability-geodesics-symmetric-i}.

(ii) The proof of \cref{eq:thm-stability-geodesics-symmetric-ii} is analogous to the above.

\end{proof}

\subsection{Diffeomorphisms for proper and stable barycentres}

In the following we will consider Riemannian barycentres
\begin{equation}
        \Vector^* \in \argmin_{\Vector\in \Real^\dimInd} \Bigl\{ \frac{1}{2\dataPointNum}\sum_{\sumIndA = 1}^\dataPointNum \distance_{\Real^\dimInd}^\diffeo(\Vector, \Vector^\sumIndA)^2 \Bigr\}
        \label{eq:diffeo-barycentre-formal-set}
    \end{equation}
of a data set $\{\Vector^\sumIndA\}_{\sumIndA=1}^\dataPointNum \subset \Real^\dimInd$ on $(\Real^\dimInd, (\cdot,\cdot)^\diffeo)$. The above-mentioned requirement that diffeomorphisms should map the data manifold into a low-dimensional geodesic subspace for proper barycentres is less obvious than for geodesics. Before we can make any claim, we need to be somewhat careful about \emph{well-posedness of barycentres}, i.e., existence and uniqueness, in the first place.


As a first step towards well-posedness (\cref{thm:well-posedness-barycentre}), the lemma below will help us by characterizing a class of real-valued functions on $\Real^\dimInd$ that are (strongly) geodesically convex on $(\Real^\dimInd,  (\cdot,\cdot)^\diffeo)$.

\begin{lemma}
\label{lem:strong-convexity-composition}
    Let $(\manifold, (\cdot,\cdot))$ be a $\dimInd$-dimensional Riemannian manifold and let $\diffeo:\Real^\dimInd \to \manifold$ be a smooth diffeomorphism such that $\diffeo(\Real^\dimInd) \subset \manifold$ is a geodesically convex set. 
    
    If a function $f:\diffeo(\Real^\dimInd)\to \Real$ is (strongly) geodesically convex on the Riemannian manifold $(\diffeo(\Real^\dimInd), (\cdot,\cdot))$, then $f\circ \diffeo: \Real^\dimInd \to \Real$ is (strongly) geodesically convex on $(\Real^\dimInd,  (\cdot,\cdot)^\diffeo)$.
\end{lemma}

\begin{proof}
    We need to check that for every geodesic $\geodesic^\diffeo_{\Vector, \VectorB}(t)\in\Real^\dimInd$ we have
    \begin{equation}
        (f\circ \diffeo) (\geodesic^\diffeo_{\Vector, \VectorB}(t)) \leq (1-t) (f\circ \diffeo) (\Vector) + t (f\circ \diffeo) (\VectorB)
        \label{eq:lem-pull-back-convexity}
    \end{equation}
    and that we have a strict inequality in the strongly convex case. We can verify \cref{eq:lem-pull-back-convexity} directly:
    \begin{multline}
        (f\circ \diffeo) (\geodesic^\diffeo_{\Vector, \VectorB}(t)) = f (\diffeo(\geodesic^\diffeo_{\Vector, \VectorB}(t)))  \overset{\text{\cref{thm:pull-back-mappings} (i)}}{=} f (\diffeo(\diffeo^{-1}(\geodesic_{\diffeo(\Vector), \diffeo(\VectorB)}(t)))) = f (\geodesic_{\diffeo(\Vector), \diffeo(\VectorB)}(t)) \\
        \overset{\text{convexity}}{\leq} (1-t) f (\diffeo(\Vector)) + t f (\diffeo(\VectorB)) = (1-t) (f\circ \diffeo) (\Vector) + t (f\circ \diffeo) (\VectorB),
    \end{multline}
    and find the claim for strong geodesic convexity analogously by replacing the inequality with a strict inequality.
\end{proof}

Next, \cref{lem:strong-convexity-composition} will help showing existence and uniqueness of barycentres, which will then allow us to find the barycentre in a more convenient way. That is, the following theorem tells us that we can compute the barycentre of a data set $\{\diffeo(\Vector^\sumIndA)\}_{\sumIndA=1}^\dataPointNum \subset \manifold$ on $(\manifold, (\cdot,\cdot))$ first and then map it back. From this observation it follows that diffeomorphisms should map the data manifold into a low-dimensional geodesic subspace for a proper barycentre that will live on the data manifold as well.


\begin{theorem}[Well-posedness of the Riemannian barycentre]
\label{thm:well-posedness-barycentre}
    Let $(\manifold, (\cdot,\cdot))$ be a $\dimInd$-dimensional Riemannian manifold and let $\diffeo:\Real^\dimInd \to \manifold$ be a smooth diffeomorphism such that $\diffeo(\Real^\dimInd) \subset \manifold$ is a geodesically convex set. Furthermore,  let $\{\Vector^\sumIndA\}_{\sumIndA=1}^\dataPointNum, \subset \Real^\dimInd$ be a data set.
    
    Then, the Riemannian barycentre of the data set $\{\Vector^\sumIndA\}_{\sumIndA=1}^\dataPointNum$
    \begin{equation}
        \Vector^* := \argmin_{\Vector\in \Real^\dimInd} \Bigl\{ \frac{1}{2\dataPointNum}\sum_{\sumIndA = 1}^\dataPointNum \distance_{\Real^\dimInd}^\diffeo(\Vector, \Vector^\sumIndA)^2 \Bigr\}
        \label{eq:diffeo-barycentre-formal}
    \end{equation}
    is well-defined and also satisfies
    \begin{equation}
    \Vector^* = \diffeo^{-1}(\mPoint^*), \quad \text{where}\quad \mPoint^* := \argmin_{\mPoint\in \diffeo(\Real^\dimInd)} \Bigl\{ \frac{1}{2\dataPointNum} \sum_{\sumIndA = 1}^\dataPointNum \distance_{\manifold}(\mPoint, \diffeo(\Vector^\sumIndA))^2 \Bigr\}.
    \label{eq:diffeo-barycentre}
\end{equation}
\end{theorem}

\begin{proof}
    First, we note that the mapping $\Vector \mapsto \distance_{\Real^\dimInd}^\diffeo(\Vector, \VectorB)^2$ is strongly geodesically convex for every $\VectorB\in \Real^\dimInd$. Indeed, this directly follows from \cref{lem:strong-convexity-composition} after rewriting the function into the composition form, i.e.,  $\distance_{\Real^\dimInd}^\diffeo(\cdot, \VectorB)^2 \overset{\text{\cref{thm:pull-back-mappings} (iv)}}{=} \distance_{\manifold}(\diffeo(\cdot), \diffeo(\VectorB))^2$, and realizing that $\mPoint\mapsto \distance_{\manifold}(\mPoint, \diffeo(\VectorB))^2$ is strongly geodesically convex on the geodesically convex set $\diffeo(\Real^\dimInd)$. Since strong geodesic convexity is closed under addition and multiplication with positive scalars, we conclude that the mapping
    \begin{equation}
        \Vector\mapsto \frac{1}{2\dataPointNum}\sum_{\sumIndA = 1}^\dataPointNum \distance_{\Real^\dimInd}^\diffeo(\Vector, \Vector^\sumIndA)^2
        \label{eq:barycentre-loss}
    \end{equation}
    is strongly geodesically convex.

    Well-posedness of $\Vector^*$ in \cref{eq:diffeo-barycentre-formal} follows from the standard argument by Karcher \cite{karcher1977riemannian} that any geodesically convex set in $\Real^\dimInd$ containing $\{\Vector^\sumIndA\}_{\sumIndA=1}^\dataPointNum$ must have non-zero and outward pointing Riemannian gradients of the function in \cref{eq:barycentre-loss}. So there must be a minimiser in the closure of this geodesically convex set, which gives us existence. Uniqueness follows from strong convexity on $\Real^\dimInd$ we showed before.
    
    Similarly, existence and uniqueness follows for $\mPoint^*$ in \cref{eq:diffeo-barycentre}. It is then easily checked that $\diffeo^{-1}(\mPoint^*)$ satisfies the first-order optimality conditions:
    \begin{multline}
        \operatorname{grad} \frac{1}{2\dataPointNum}\sum_{\sumIndA = 1}^\dataPointNum \distance_{\Real^\dimInd}^\diffeo(\cdot, \Vector^\sumIndA)^2 \mid_{\diffeo^{-1}(\mPoint^*)} = - \frac{1}{\dataPointNum} \sum_{\sumIndA = 1}^\dataPointNum \log^\diffeo_{\diffeo^{-1}(\mPoint^*)} \Vector^\sumIndA
        \overset{\text{\cref{thm:pull-back-mappings} (ii)}}{=} - \frac{1}{\dataPointNum} \sum_{\sumIndA = 1}^\dataPointNum \diffeo_*^{-1}[\log_{\mPoint^*} \diffeo(\Vector^\sumIndA) ] \\
        =  \diffeo_*^{-1}[ - \frac{1}{\dataPointNum} \sum_{\sumIndA = 1}^\dataPointNum \log_{\mPoint^*} \diffeo(\Vector^\sumIndA) ] = \diffeo_*^{-1}[0_{\mPoint^*}] = 0_{\diffeo^{-1}(\mPoint^*)}.
    \end{multline}
    So we conclude that $\Vector^* = \diffeo^{-1}(\mPoint^*)$.
\end{proof}

\begin{corollary}
\label{cor:euclidean-diffeo-mean}
    If $(\manifold, (\cdot,\cdot)):= (\Real^\dimInd, (\cdot,\cdot)_2)$, the barycentre problem has a closed-form solution given by 
    \begin{equation}
        \Vector^* = \diffeo^{-1}\Bigl(\frac{1}{\dataPointNum} \sum_{\sumIndA=1}^\dataPointNum \diffeo(\Vector^\sumIndA) \Bigr)
        \label{eq:euclidean-diffeo-mean}
    \end{equation}
\end{corollary}

Our second goal in this section is to characterize stability of barycentres with respect to the data set $\{\Vector^\sumIndA\}_{\sumIndA=1}^\dataPointNum$ and we will see that $\ell^2$-stability is once again fully determined by the diffeomorphism's deviation from isometry around the data manifold and the curvature of $(\manifold, (\cdot, \cdot))$. Our strategy for characterizing stability requires less effort than for geodesics, but we will need to make an approximation. In particular, given the barycentre $\Vector^*$ of the original data set and a new data set $\{\VectorB^{\sumIndA}\}_{\sumIndA=1}^\dataPointNum$, we will consider the \emph{approximate barycentre} of $\{\VectorB^{\sumIndA}\}_{\sumIndA=1}^\dataPointNum \subset \Real^\dimInd$ from $\Vector^*$ given by
\begin{equation}
    \tilde{\VectorB}_{\Vector^*}:= \exp_{\Vector^*}^\diffeo \Bigl(\frac{1}{\dataPointNum} \sum_{\sumIndA=1}^\dataPointNum \log_{\Vector^*}^\diffeo (\VectorB^\sumIndA) \Bigr) = \diffeo^{-1}\Bigl(\exp_{\diffeo(\Vector^*)} \Bigl(\frac{1}{\dataPointNum} \sum_{\sumIndA=1}^\dataPointNum \log_{\diffeo(\Vector^*)} \diffeo(\VectorB^\sumIndA) \Bigr) \Bigr),
    \label{eq:approximate-barycentre}
\end{equation}
for analyzing stability (\cref{thm:stability-barycentre}).

In general, the approximate barycentre $\tilde{\VectorB}_{\Vector^*}$ and the actual barycentre $\VectorB^*$ cannot be expected to be the same. As a first motivation for our choice, consider that we do get consistency for a special case.


\begin{proposition}
    If $(\manifold, (\cdot,\cdot)):= (\Real^\dimInd, (\cdot,\cdot)_2)$, the approximate barycentre $\tilde{\VectorB}_{\Vector^*}$ is the barycentre on $(\Real^\dimInd, (\cdot,\cdot)^\diffeo)$.
\end{proposition}
 
\begin{proof}
    \begin{multline}
    \tilde{\VectorB}_{\Vector^*} 
    = \diffeo^{-1}\Bigl(\exp_{\diffeo(\Vector^*)} \Bigl(\frac{1}{\dataPointNum} \sum_{\sumIndA=1}^\dataPointNum \log_{\diffeo(\Vector^*)} \diffeo(\VectorB^\sumIndA) \Bigr) \Bigr) 
    = \diffeo^{-1}\Bigl(\diffeo(\Vector^*) + \Bigl(\frac{1}{\dataPointNum} \sum_{\sumIndA=1}^\dataPointNum  \diffeo(\VectorB^\sumIndA) - \diffeo(\Vector^*) \Bigr) \Bigr)  
    \\
    = \diffeo^{-1}\Bigl(\frac{1}{\dataPointNum} \sum_{\sumIndA=1}^\dataPointNum  \diffeo(\VectorB^\sumIndA) \Bigr) 
    \overset{\text{\cref{cor:euclidean-diffeo-mean}}}{=} \VectorB^*
\end{multline}
\end{proof}

More generally, we do expect that the approximate barycentre is a good first-order approximation of the actual barycentre, because it is just one Riemannian gradient descent step with properly chosen step size from initialisation $\Vector^*$ for solving the barycentre problem. So if the $\{\VectorB^{\sumIndA}\}_{\sumIndA=1}^\dataPointNum$ are small variations of the $\{\Vector^{\sumIndA}\}_{\sumIndA=1}^\dataPointNum$, we expect that the approximate barycentre is a good object to analyze for stability. 

Finally, we are ready to quantify stability of barycentres. The following theorem tells us that $\ell^2$-stability is fully determined by the $\diffeo$'s deviation from isometry around the data manifold and the curvature of $(\manifold, (\cdot, \cdot))$. In particular, positive curvature will inherently cause instabilities.

\begin{theorem}[Stability of the approximate Riemannian barycentre]
\label{thm:stability-barycentre}
    Let $(\manifold, (\cdot,\cdot))$ be a $\dimInd$-dimensional symmetric Riemannian manifold and let $\diffeo:\Real^\dimInd \to \manifold$ be a smooth diffeomorphism such that $\diffeo(\Real^\dimInd) \subset \manifold$ is a geodesically convex set. Furthermore, let $\{\Vector^\sumIndA\}_{\sumIndA=1}^\dataPointNum \subset \Real^\dimInd$ be a data set with Riemannian barycentre \cref{eq:diffeo-barycentre-formal} $\Vector^* \in \Real^\dimInd$, let $\{\kappa^\sumIndA_{\sumIndB}\}_{\sumIndB=1}^{\dimInd} \subset \Real$ be the eigenvalues of the operators
    \begin{equation}
    \Theta_{\diffeo(\Vector^\sumIndA)} \mapsto \curvature_{\diffeo(\Vector^\sumIndA)}(\Theta_{\diffeo(\Vector^\sumIndA)}, \log_{\diffeo(\Vector^\sumIndA)}\diffeo(\Vector^*)) \log_{\diffeo(\Vector^\sumIndA)}\diffeo(\Vector^*), \quad \text{for } \sumIndA=1, \ldots, \dataPointNum.
    \label{eq:thm-bary-stab-curvature-op}
\end{equation}
and define $\betaLogmPoint: \Real \to \Real$ as 
\begin{equation}
    \betaLogmPoint(\kappa) := \left\{\begin{matrix}
\frac{\sqrt{-\kappa}}{\sinh(\sqrt{-\kappa})}, & \kappa <0, \\
1, & \kappa = 0, \\
\frac{\sqrt{\kappa}}{\sin(\sqrt{\kappa})}, & \kappa >0.
\end{matrix}\right.
\label{eq:beta-log-point}
\end{equation}
Finally, consider for $\sumIndA=1, \ldots, \dataPointNum$ open neighbourhoods $\mathcal{U}(\Vector^\sumIndA)\subset \Real^\dimInd$ on which $\diffeo$ has local Lipschitz constants $\operatorname{Lip}_{\Vector^\sumIndA}(\diffeo)$ and consider an open neighbourhood and $\mathcal{V}(\diffeo(\Vector^*)) \subset \manifold$ on which $\diffeo^{-1}$ has local Lipschitz constant $\operatorname{Lip}_{\diffeo(\Vector^*)}(\diffeo^{-1})$.

Then, as $\{\VectorB^{\sumIndA}\}_{\sumIndA=1}^\dataPointNum \to \{\Vector^\sumIndA\}_{\sumIndA=1}^\dataPointNum$ the approximate Riemannian barycentre $\tilde{\VectorB}_{\Vector^*}\in \Real^\dimInd$ defined as in \cref{eq:approximate-barycentre}
behaves as
\begin{equation}
    \|\tilde{\VectorB}_{\Vector^*} - \Vector^*\|_2 \leq \frac{\operatorname{Lip}_{\diffeo(\Vector^*)}(\diffeo^{-1})}{\dataPointNum} \sum_{\sumIndA=1}^{\dataPointNum} \betaLogmPoint(\kappa^\sumIndA_{\max}) \operatorname{Lip}_{\Vector^\sumIndA}(\diffeo)\|\VectorB^\sumIndA - \Vector^\sumIndA\|_2 + o(\|\VectorB^\sumIndA - \Vector^\sumIndA\|_2).
\end{equation}

\end{theorem}

\begin{proof}
    Let  $\{\Theta^{\sumIndB}_{\diffeo(\Vector^\sumIndA)}\}_{\sumIndB=1}^{\dimInd}\subset \tangent_{\diffeo(\Vector^\sumIndA)} \manifold$ for $\sumIndA=1, \ldots, \dataPointNum$ be orthonormal frames that diagonalize the operators in \cref{eq:thm-bary-stab-curvature-op} corresponding to the eigenvalues $\{\kappa^\sumIndA_{\sumIndB}\}_{\sumIndB=1}^{\dimInd}$. Next, as $\{\VectorB^{\sumIndA}\}_{\sumIndA=1}^\dataPointNum \to \{\Vector^\sumIndA\}_{\sumIndA=1}^\dataPointNum$ we may assume that $\VectorB^{\sumIndA} \in \mathcal{U}(\Vector^\sumIndA)$ and $\tilde{\VectorB}_{\Vector^*} \in \mathcal{V}(\diffeo(\Vector^*))$. Finally, before moving on to proving the statement, consider the Taylor approximation
    \begin{multline}
        \log_{\diffeo(\Vector^*)} \diffeo(\VectorB^\sumIndA) = \log_{\diffeo(\Vector^*)} \diffeo(\Vector^\sumIndA) + D_{\diffeo(\Vector^\sumIndA)} \log_{\diffeo(\Vector^*)} (\cdot) [\log_{\diffeo(\Vector^\sumIndA)} \diffeo(\VectorB^\sumIndA)] + \mathcal{O}(\distance_\manifold(\diffeo(\VectorB^\sumIndA),\diffeo( \Vector^\sumIndA))^2)\\
        \overset{\text{\cite[Lemma~1]{bergmann2019recent}}}{=} \log_{\diffeo(\Vector^*)} \diffeo(\Vector^\sumIndA) + \sum_{\sumIndB=1}^\dimInd \betaLogmPoint(\kappa^\sumIndA_{\sumIndB})\bigl(\log_{\diffeo(\Vector^\sumIndA)} \diffeo(\VectorB^\sumIndA), \Theta^{\sumIndB}_{\diffeo(\Vector^\sumIndA)}\bigr)_{\diffeo(\Vector^\sumIndA)} \mathcal{P}_{\diffeo(\Vector^*) \leftarrow \diffeo(\Vector^\sumIndA)} \Theta^{\sumIndB}_{\diffeo(\Vector^\sumIndA)} \\
        + \mathcal{O}(\distance_\manifold(\diffeo(\VectorB^\sumIndA),\diffeo( \Vector^\sumIndA))^2)
        \label{eq:thm-bary-stability-taylor-log}
    \end{multline}
    and note that
    \begin{multline}
         \frac{1}{\dataPointNum} \sum_{\sumIndA=1}^\dataPointNum \log_{\diffeo(\Vector^*)} \diffeo(\VectorB^\sumIndA) = \frac{1}{\dataPointNum} \sum_{\sumIndA=1}^\dataPointNum \sum_{\sumIndB=1}^\dimInd \betaLogmPoint(\kappa^\sumIndA_{\sumIndB})\bigl(\log_{\diffeo(\Vector^\sumIndA)} \diffeo(\VectorB^\sumIndA), \Theta^{\sumIndB}_{\diffeo(\Vector^\sumIndA)}\bigr)_{\diffeo(\Vector^\sumIndA)} \mathcal{P}_{\diffeo(\Vector^*) \leftarrow \diffeo(\Vector^\sumIndA)} \Theta^{\sumIndB}_{\diffeo(\Vector^\sumIndA)} \\
         + \mathcal{O}(\distance_\manifold(\diffeo(\VectorB^\sumIndA),\diffeo( \Vector^\sumIndA))^2),
         \label{eq:thm-bary-stability-taylor-log-sum}
    \end{multline}
    because $\frac{1}{\dataPointNum} \sum_{\sumIndA=1}^\dataPointNum \log_{\diffeo(\Vector^*)} \diffeo(\Vector^\sumIndA) = 0_{\diffeo(\Vector^*)}$ due to first-order optimality conditions.
    
    Subsequently,
    \begin{subequations}
    \allowdisplaybreaks
    \begin{multline}
        \|\tilde{\VectorB}_{\Vector^*} - \Vector^*\|_2 = \Bigl\|\diffeo^{-1}\Bigl(\exp_{\diffeo(\Vector^*)} \Bigl(\frac{1}{\dataPointNum} \sum_{\sumIndA=1}^\dataPointNum \log_{\diffeo(\Vector^*)} \diffeo(\VectorB^\sumIndA) \Bigr) \Bigr) - \diffeo^{-1}(\diffeo(\Vector^*)) \Bigr\|_2 \\
        \leq \operatorname{Lip}_{\diffeo(\Vector^*)}(\diffeo^{-1}) \distance_\manifold \Bigl(\exp_{\diffeo(\Vector^*)} \Bigl(\frac{1}{\dataPointNum} \sum_{\sumIndA=1}^\dataPointNum \log_{\diffeo(\Vector^*)} \diffeo(\VectorB^\sumIndA) \Bigr), \diffeo(\Vector^*) \Bigr)\\
        = \operatorname{Lip}_{\diffeo(\Vector^*)}(\diffeo^{-1}) \Bigl\| \frac{1}{\dataPointNum} \sum_{\sumIndA=1}^\dataPointNum \log_{\diffeo(\Vector^*)} \diffeo(\VectorB^\sumIndA) \Bigr\|_{\diffeo(\Vector^*)}\\
        \overset{\text{\cref{eq:thm-bary-stability-taylor-log-sum}}}{=} 
        \frac{\operatorname{Lip}_{\diffeo(\Vector^*)}(\diffeo^{-1})}{\dataPointNum} \Bigl\| \sum_{\sumIndA,\sumIndB=1}^{\dataPointNum,\dimInd} \betaLogmPoint(\kappa^\sumIndA_{\sumIndB})\bigl(\log_{\diffeo(\Vector^\sumIndA)} \diffeo(\VectorB^\sumIndA), \Theta^{\sumIndB}_{\diffeo(\Vector^\sumIndA)}\bigr)_{\diffeo(\Vector^\sumIndA)} \mathcal{P}_{\diffeo(\Vector^*) \leftarrow \diffeo(\Vector^\sumIndA)} \Theta^{\sumIndB}_{\diffeo(\Vector^\sumIndA)} + \mathcal{O}(\distance_\manifold(\diffeo(\VectorB^\sumIndA),\diffeo( \Vector^\sumIndA))^2) \Bigr\|_{\diffeo(\Vector^*)}\\
        \leq \frac{\operatorname{Lip}_{\diffeo(\Vector^*)}(\diffeo^{-1})}{\dataPointNum} \sum_{\sumIndA=1}^\dataPointNum \Bigl\|  \sum_{\sumIndB=1}^\dimInd \betaLogmPoint(\kappa^\sumIndA_{\sumIndB})\bigl(\log_{\diffeo(\Vector^\sumIndA)} \diffeo(\VectorB^\sumIndA), \Theta^{\sumIndB}_{\diffeo(\Vector^\sumIndA)}\bigr)_{\diffeo(\Vector^\sumIndA)} \mathcal{P}_{\diffeo(\Vector^*) \leftarrow \diffeo(\Vector^\sumIndA)} \Theta^{\sumIndB}_{\diffeo(\Vector^\sumIndA)} + \mathcal{O}(\distance_\manifold(\diffeo(\VectorB^\sumIndA),\diffeo( \Vector^\sumIndA))^2) \Bigr\|_{\diffeo(\Vector^*)} \\
        \leq \frac{\operatorname{Lip}_{\diffeo(\Vector^*)}(\diffeo^{-1})}{\dataPointNum} \sum_{\sumIndA=1}^\dataPointNum \Bigl\|  \sum_{\sumIndB=1}^\dimInd \betaLogmPoint(\kappa^\sumIndA_{\sumIndB})\bigl(\log_{\diffeo(\Vector^\sumIndA)} \diffeo(\VectorB^\sumIndA), \Theta^{\sumIndB}_{\diffeo(\Vector^\sumIndA)}\bigr)_{\diffeo(\Vector^\sumIndA)} \mathcal{P}_{\diffeo(\Vector^*) \leftarrow \diffeo(\Vector^\sumIndA)} \Theta^{\sumIndB}_{\diffeo(\Vector^\sumIndA)} \Bigr\|_{\diffeo(\Vector^*)} + o(\distance_\manifold(\diffeo(\VectorB^\sumIndA),\diffeo( \Vector^\sumIndA)))\\
        = \frac{\operatorname{Lip}_{\diffeo(\Vector^*)}(\diffeo^{-1})}{\dataPointNum} \sum_{\sumIndA=1}^\dataPointNum  \sqrt{\sum_{\sumIndB=1}^\dimInd \betaLogmPoint(\kappa^\sumIndA_{\sumIndB})^2 \bigl(\log_{\diffeo(\Vector^\sumIndA)} \diffeo(\VectorB^\sumIndA), \Theta^{\sumIndB}_{\diffeo(\Vector^\sumIndA)}\bigr)_{\diffeo(\Vector^\sumIndA)}^2}  + o(\distance_\manifold(\diffeo(\VectorB^\sumIndA),\diffeo( \Vector^\sumIndA))) \\
        \leq \frac{\operatorname{Lip}_{\diffeo(\Vector^*)}(\diffeo^{-1})}{\dataPointNum} \sum_{\sumIndA=1}^\dataPointNum  \sqrt{ \betaLogmPoint(\kappa^\sumIndA_{\max})^2\sum_{\sumIndB=1}^\dimInd \bigl(\log_{\diffeo(\Vector^\sumIndA)} \diffeo(\VectorB^\sumIndA), \Theta^{\sumIndB}_{\diffeo(\Vector^\sumIndA)}\bigr)_{\diffeo(\Vector^\sumIndA)}^2}  + o(\distance_\manifold(\diffeo(\VectorB^\sumIndA),\diffeo( \Vector^\sumIndA)))\\
        = \frac{\operatorname{Lip}_{\diffeo(\Vector^*)}(\diffeo^{-1})}{\dataPointNum} \sum_{\sumIndA=1}^\dataPointNum \betaLogmPoint(\kappa^\sumIndA_{\max})\distance_\manifold(\diffeo(\VectorB^\sumIndA),\diffeo( \Vector^\sumIndA))  + o(\distance_\manifold(\diffeo(\VectorB^\sumIndA),\diffeo( \Vector^\sumIndA)))\\
        \leq \frac{\operatorname{Lip}_{\diffeo(\Vector^*)}(\diffeo^{-1})}{\dataPointNum} \sum_{\sumIndA=1}^{\dataPointNum} \betaLogmPoint(\kappa^\sumIndA_{\max}) \operatorname{Lip}_{\Vector^\sumIndA}(\diffeo)\|\VectorB^\sumIndA - \Vector^\sumIndA\|_2 + o(\|\VectorB^\sumIndA - \Vector^\sumIndA\|_2).
    \end{multline}
    \end{subequations}
\end{proof}

\section{Non-linear compression on $(\Real^\dimInd, (\cdot,\cdot)^\diffeo)$ and applications}
\label{sec:non-linear-compression-pullback}

Now that we know from \cref{thm:pull-back-mappings,thm:well-posedness-barycentre} that we need to map the data manifold into a geodesic subspace and from \cref{thm:stability-geodesics-symmetric,thm:stability-barycentre} that we need to preserve local isometry around the data manifold while doing so, we shift our focus in this section towards how to modify existing algorithms for symmetric Riemannian manifolds to $(\Real^\dimInd, (\cdot,\cdot)^\diffeo)$. In particular, we will consider how we can piggyback off of existing theory for data compression on symmetric Riemannian manifolds and how a recent efficient algorithm for low rank approximation enables construction of the \emph{Riemannian autoencoder} (RAE) and \emph{curvature corrected Riemannian autoencoder} (CC-RAE) mapping, which are non-linear compression mappings that have additional nice mathematical properties such as an interpretable latent space that traditional neural network-based autoencoders do not have. Although the focus is on algorithm design, we will also observe several times that the above-mentioned best practices for diffeomorphisms $\diffeo:\Real^\dimInd \to \manifold$ are necessary for useful and efficient algorithms.



\subsection{Diffeomorphisms for efficient compression}

In the following we will consider the general problem of efficiently compressing a data set $\{\Vector^\sumIndA\}_{\sumIndA=1}^\dataPointNum \subset \Real^\dimInd$ through solving
\begin{equation}
    \inf_{(\tangentVector_{\VectorC}^{1}, \ldots, \tangentVector_{\VectorC}^{\dataPointNum})\in \mathcal{S}_{\VectorC}((\Real^{\dimInd})^\dataPointNum)} \Bigl\{\sum_{\sumIndA=1}^\dataPointNum \| \exp_{\VectorC}^\diffeo(\tangentVector_{\VectorC}^{\sumIndA}) - \Vector^\sumIndA\|_2^2 \Bigr\},
    \label{eq:global-approximation-problem-l2}
\end{equation}
where $\mathcal{S}_{\VectorC}\bigl((\Real^{\dimInd})^\dataPointNum\bigr) \subset \tangent_{\VectorC}(\Real^{\dimInd})^\dataPointNum$ is a suitable class of approximators defined by a tangent space subset at $\VectorC \in \Real^\dimInd$. Solving \cref{eq:global-approximation-problem-l2} directly -- given a $\mathcal{S}_{\VectorC}\bigl((\Real^{\dimInd})^\dataPointNum\bigr)$ -- will in general be hard. Instead, we will relax \cref{eq:global-approximation-problem-l2} several times into a more tractable problem using the requirements on the diffeomorphism from the previous section and using recent ideas from approximation of manifold-valued data on symmetric Riemannian manifolds. Then, in the second part of this section we will consider a concrete example for $\mathcal{S}_{\VectorC}\bigl((\Real^{\dimInd})^\dataPointNum\bigr)$ and solve that using a recently proposed algorithm.

For the first relaxation, we argue that instead of solving \cref{eq:global-approximation-problem-l2} we can solve
\begin{equation}
    \inf_{(\tangentVector_{\VectorC}^{1}, \ldots, \tangentVector_{\VectorC}^{\dataPointNum})\in \mathcal{S}_{\VectorC}((\Real^{\dimInd})^\dataPointNum)} \Bigl\{\sum_{\sumIndA=1}^N \distance_{\Real^\dimInd}^\diffeo (\exp_{\VectorC}^\diffeo(\tangentVector_{\VectorC}^{\sumIndA}), \Vector^\sumIndA)^2 \Bigr\}.
    \label{eq:global-approximation-problem}
\end{equation}
The following results give a hint as to why this might work. That is, each term in the summand of \cref{eq:global-approximation-problem-l2} can be upper and lower bounded by multiples of the corresponding term in \cref{eq:global-approximation-problem}.

\begin{proposition}
    Let $(\manifold, (\cdot,\cdot))$ be a $\dimInd$-dimensional symmetric Riemannian manifold and let $\diffeo:\Real^\dimInd \to \manifold$ be a smooth diffeomorphism such that $\diffeo(\Real^\dimInd) \subset \manifold$ is a geodesically convex set. Furthermore, let $\Vector \in \Real^\dimInd$ be any point and consider open neighbourhoods $\mathcal{U}(\Vector)\subset \Real^\dimInd$ on which $\diffeo$ has local Lipschitz constants $\operatorname{Lip}_{\Vector}(\diffeo)$ and $\mathcal{V}(\diffeo(\Vector)) \subset \manifold$ on which $\diffeo^{-1}$ has local Lipschitz constant $\operatorname{Lip}_{\diffeo(\Vector)}(\diffeo^{-1})$. Finally, let $\tangentVector_\VectorC \in \tangent_\VectorC \Real^\dimInd$ be a tangent vector at any point $\VectorC \in \Real^\dimInd$ such that $\exp_{\VectorC}^\diffeo(\tangentVector_\VectorC) \in \mathcal{U}(\Vector)$ and $\diffeo(\exp_{\VectorC}^\diffeo(\tangentVector_\VectorC)) \in \mathcal{V}(\diffeo(\Vector))$.

    Then,
    \begin{equation}
        \frac{1}{\operatorname{Lip}_{\Vector} (\diffeo)} \distance_{\Real^\dimInd}^\diffeo (\exp_{\VectorC}^\diffeo(\tangentVector_\VectorC),  \Vector) \leq \| \exp_{\VectorC}^\diffeo(\tangentVector_\VectorC) - \Vector\|_2 \leq \operatorname{Lip}_{\diffeo(\Vector)} (\diffeo^{-1}) \distance_{\Real^\dimInd}^\diffeo (\exp_{\VectorC}^\diffeo(\tangentVector_\VectorC), \Vector).
        \label{eq:lem-sandwich-approximation-error}
    \end{equation}
\end{proposition}

\begin{proof}
    The bounds in \cref{eq:lem-sandwich-approximation-error} follow directly from
    \begin{multline}
        \frac{1}{\operatorname{Lip}_{\Vector} (\diffeo)} \distance_{\Real^\dimInd}^\diffeo (\exp_{\VectorC}^\diffeo(\tangentVector_\VectorC),  \Vector) \overset{\text{\cref{thm:pull-back-mappings} (iv)}}{=} \frac{1}{\operatorname{Lip}_{\Vector} (\diffeo)} \distance_{\manifold}(\diffeo (\exp_{\VectorC}^\diffeo(\tangentVector_\VectorC)), \diffeo (\Vector)) \leq \frac{\operatorname{Lip}_{\Vector} (\diffeo)}{\operatorname{Lip}_{\Vector} (\diffeo)} \|\exp_{\VectorC}^\diffeo(\tangentVector_\VectorC)) - \Vector \|_2 \\
        = \|\exp_{\VectorC}^\diffeo(\tangentVector_\VectorC)) - \Vector \|_2 = \|\diffeo^{-1}(\diffeo(\exp_{\VectorC}^\diffeo(\tangentVector_\VectorC))) ) - \diffeo^{-1}(\diffeo(\Vector)) \|_2 \leq \operatorname{Lip}_{\diffeo(\Vector)} (\diffeo^{-1}) \distance_{\manifold}(\diffeo (\exp_{\VectorC}^\diffeo(\tangentVector_\VectorC)), \diffeo (\Vector) \\
        \overset{\text{\cref{thm:pull-back-mappings} (iv)}}{=}  \operatorname{Lip}_{\diffeo(\Vector)} (\diffeo^{-1}) \distance_{\Real^\dimInd}^\diffeo (\exp_{\VectorC}^\diffeo(\tangentVector_\VectorC),  \Vector).
    \end{multline}
\end{proof}


\begin{corollary}
\label{cor:global-approximation-l2-sandwich}
Let $(\manifold, (\cdot,\cdot))$ be a $\dimInd$-dimensional symmetric Riemannian manifold and let $\diffeo:\Real^\dimInd \to \manifold$ be a smooth diffeomorphism such that $\diffeo(\Real^\dimInd) \subset \manifold$ is a geodesically convex set. Furthermore, let $\{\Vector^\sumIndA\}_{\sumIndA=1}^\dataPointNum \subset \Real^\dimInd$ be a data set and consider for $\sumIndA=1, \ldots, \dataPointNum$ open neighbourhoods $\mathcal{U}(\Vector^\sumIndA)\subset \Real^\dimInd$
on which $\diffeo$ has local Lipschitz constants $\operatorname{Lip}_{\Vector^\sumIndA}(\diffeo)$ and open neighbourhoods $\mathcal{V}(\diffeo(\Vector^\sumIndA)) \subset \manifold$ on which $\diffeo^{-1}$ has local Lipschitz constant $\operatorname{Lip}_{\diffeo(\Vector^\sumIndA)}(\diffeo^{-1})$. Finally, let $ \{\tangentVector_{\VectorC}^{\sumIndA}\}_{\sumIndA=1}^\dataPointNum \subset \tangent_{\VectorC}\Real^{\dimInd}$ be a set of tangent vectors at any point $\VectorC \in \Real^\dimInd$ such that $\exp_{\VectorC}^\diffeo(\tangentVector_\VectorC^\sumIndA) \in \mathcal{U}(\Vector^\sumIndA)$ and $\diffeo(\exp_{\VectorC}^\diffeo(\tangentVector_\VectorC^\sumIndA)) \in \mathcal{V}(\diffeo(\Vector^\sumIndA))$ for $\sumIndA=1, \ldots, \dataPointNum$.

Then,
\begin{equation}
    \frac{1}{\operatorname{Lip}_{\{\Vector^\sumIndA\}_{\sumIndA=1}^\dataPointNum} (\diffeo)^2 } \sum_{\sumIndA=1}^N \distance_{\Real^\dimInd}^\diffeo (\exp_{\VectorC}^\diffeo(\tangentVector_{\VectorC}^{\sumIndA}), \Vector^\sumIndA)^2 \leq \sum_{\sumIndA=1}^\dataPointNum \| \exp_{\VectorC}^\diffeo(\tangentVector_{\VectorC}^{\sumIndA}) - \Vector^\sumIndA\|_2^2 \leq \operatorname{Lip}_{\{\diffeo(\Vector^\sumIndA)\}_{\sumIndA=1}^\dataPointNum} (\diffeo^{-1})^2 \sum_{\sumIndA=1}^N \distance_{\Real^\dimInd}^\diffeo (\exp_{\VectorC}^\diffeo(\tangentVector_{\VectorC}^{\sumIndA}), \Vector^\sumIndA)^2,
\end{equation}
where 
\begin{equation}
    \operatorname{Lip}_{\{\Vector^\sumIndA\}_{\sumIndA=1}^\dataPointNum} (\diffeo) :=\max_{\sumIndA=1, \ldots, \dataPointNum} \{\operatorname{Lip}_{\Vector^\sumIndA} (\diffeo) \} \quad \text{and} \quad \operatorname{Lip}_{\{\diffeo(\Vector^\sumIndA)\}_{\sumIndA=1}^\dataPointNum} (\diffeo^{-1}) := \max_{\sumIndA=1, \ldots, \dataPointNum} \{\operatorname{Lip}_{\diffeo(\Vector^\sumIndA)} (\diffeo^{-1}) \}.
\end{equation}
    
\end{corollary}

From \cref{cor:global-approximation-l2-sandwich} we then see that if $\diffeo$ is a local isometry on each of the data points (so around the data manifold), the upper and the lower bounds are equal. In other words, \cref{eq:global-approximation-problem-l2} and \cref{eq:global-approximation-problem} are equivalent under the isometry requirement, which thus allows us to solve for \cref{eq:global-approximation-problem} instead. 

The reason why \cref{eq:global-approximation-problem} is more manageable than \cref{eq:global-approximation-problem-l2} comes down to the following. In recent work, it has been shown that solving a problem of the form \cref{eq:global-approximation-problem}, can be relaxed (again) without losing on global-geometry awareness, generality on the type of approximation task and computational feasibility. For that we need the space $(\Real^\dimInd, (\cdot,\cdot)^\diffeo)$ to be symmetric, which we know to be the case by \cref{thm:local-symmetry}. In particular, we have the following identity for the loss function in \cref{eq:global-approximation-problem}

\begin{equation}
    \sum_{\sumIndA=1}^N \distance_{\Real^\dimInd}^\diffeo (\exp_{\VectorC}^\diffeo(\tangentVector_{\VectorC}^{\sumIndA}), \Vector^\sumIndA)^2 \overset{\text{\cite[Thm. 3.4]{diepeveen2023curvature}}}{=} \sum_{\sumIndA=1}^N \sum_{\sumIndB=1}^{\dimInd} \betaExptVector (\kappa^\sumIndA_{\sumIndB})^2 \Bigl(\Bigl(  \tangentVector_{\VectorC}^{\sumIndA} - \log^\diffeo_{\VectorC} (\Vector^\sumIndA), \Psi^{\sumIndA,\sumIndB}_{\VectorC} \Bigr)^\diffeo_{\VectorC}\Bigr)^2 + \mathcal{O}\bigl((\|\tangentVector_{\VectorC}^{\sumIndA} - \log^\diffeo_{\VectorC} (\Vector^\sumIndA)\|^\diffeo_{\VectorC})^3\bigr),
\end{equation}
where $\{\kappa^\sumIndA_{\sumIndB}\}_{\sumIndB=1}^{\dimInd} \subset \Real$ and $\{\Psi^{\sumIndA,\sumIndB}_\VectorC\}_{\sumIndB=1}^{\dimInd} \subset \tangent_\VectorC \Real^\dimInd$ are the eigenvalues and corresponding eigenvectors of the operators
    \begin{equation}
    \Psi_{\VectorC} \mapsto \curvature^\diffeo_\VectorC (\Psi_{\VectorC}, \log^\diffeo_{\VectorC} (\Vector^\sumIndA)) \log^\diffeo_{\VectorC} (\Vector^\sumIndA), \quad \text{for } \sumIndA=1, \ldots, \dataPointNum,
    \label{eq:thm-approx-curvature-op}
\end{equation}
and the function $\betaExptVector: \Real \to \Real$ is defined as 
\begin{equation}
    \betaExptVector(\kappa) := \left\{\begin{matrix}
\frac{\sinh(\sqrt{-\kappa})}{\sqrt{-\kappa}}, & \kappa <0, \\
1, & \kappa = 0, \\
\frac{\sin(\sqrt{\kappa})}{\sqrt{\kappa}}, & \kappa >0.
\end{matrix}\right.
\label{eq:beta-log-point}
\end{equation}

In \cite{diepeveen2023curvature} the authors show that solving
\begin{equation}
     \inf_{(\tangentVector_{\VectorC}^{1}, \ldots, \tangentVector_{\VectorC}^{\dataPointNum})\in \mathcal{S}_{\VectorC}((\Real^{\dimInd})^\dataPointNum)} \Bigl\{ \sum_{\sumIndA=1}^N \sum_{\sumIndB=1}^{\dimInd} \betaExptVector (\kappa^\sumIndA_{\sumIndB})^2 \Bigl(\Bigl(  \tangentVector_{\VectorC}^{\sumIndA} - \log^\diffeo_{\VectorC} (\Vector^\sumIndA), \Psi^{\sumIndA,\sumIndB}_{\VectorC} \Bigr)^\diffeo_{\VectorC}\Bigr)^2 \Bigr\}
    \label{eq:curvature-corrected-approximation-problem}
\end{equation}
will yield a good approximation of a minimiser of \cref{eq:global-approximation-problem}, and in turn by \cref{cor:global-approximation-l2-sandwich} a good approximation of a minimiser of our original problem \cref{eq:global-approximation-problem-l2}.

To solve \cref{eq:curvature-corrected-approximation-problem}, we need to have the eigenvalues and the corresponding eigenvectors of the operators \cref{eq:thm-approx-curvature-op}. The next lemma tells us that we can just pull them back from the embedding space.
\begin{proposition}
Let $(\manifold, (\cdot,\cdot))$ be a $\dimInd$-dimensional symmetric Riemannian manifold and let $\diffeo:\Real^\dimInd \to \manifold$ be a smooth diffeomorphism such that $\diffeo(\Real^\dimInd) \subset \manifold$ is a geodesically convex set. Furthermore, let $\Vector,\VectorB \in \Real^\dimInd$ be distinct points and let $\{\Theta_{\diffeo(\Vector)}^\sumIndB\}_{\sumIndB=1}^\dimInd \subset \tangent_{\diffeo(\Vector)}\manifold$ be the orthonormal basis that diagonalizes the operator
\begin{equation}
    \Theta_{\diffeo(\Vector)} \mapsto \curvature_{\diffeo(\Vector)}(\Theta_{\diffeo(\Vector)}, \log_{\diffeo(\Vector)}\diffeo(\VectorB)) \log_{\diffeo(\Vector)}\diffeo(\VectorB)
\end{equation}
with corresponding eigenvalues $\{\kappa_{\sumIndB}\}_{\sumIndB=1}^\dimInd \subset \Real$.

Then, $\{\diffeo^{-1}_*[\Theta_{\diffeo(\Vector)}^\sumIndB]\}_{\sumIndB=1}^\dimInd \subset \tangent_{\Vector}\Real^\dimInd$ is an orthonormal basis that diagonalizes 
\begin{equation}
    \Psi_{\Vector} \mapsto \curvature^\diffeo_\Vector (\Psi_{\Vector}, \log^\diffeo_{\Vector} \VectorB) \log^\diffeo_{\Vector} \VectorB
\end{equation}
with corresponding eigenvalues $\{\kappa_\sumIndB\}_{\sumIndB=1}^\dimInd\subset \Real$.
\end{proposition}

\begin{proof}
    We will check that $\{\diffeo^{-1}_*[\Theta_{\diffeo(\Vector)}^\sumIndB]\}_{\sumIndB=1}^\dimInd \subset \tangent_{\Vector}\Real^\dimInd$ is an orthonormal basis and that it diagonalizes the curvature operator. Indeed,
\begin{equation}
    (\diffeo^{-1}_*[\Theta_{\diffeo(\Vector)}^\sumIndA],\diffeo^{-1}_*[\Theta_{\diffeo(\Vector)}^\sumIndB])^\diffeo_{\Vector} 
    =  (\diffeo_*[\diffeo^{-1}_*[\Theta_{\diffeo(\Vector)}^\sumIndA]], \diffeo_*[\diffeo^{-1}_*[\Theta_{\diffeo(\Vector)}^\sumIndB]])_{\diffeo(\Vector)}
    = (\Theta_{\diffeo(\Vector)}^\sumIndA,\Theta_{\diffeo(\Vector)}^\sumIndB)_{\diffeo(\Vector)} 
    = \delta_{\sumIndA, \sumIndB}
\end{equation}
and
\begin{multline}
     \curvature^\diffeo(\diffeo^{-1}_*[\Theta_{\diffeo(\Vector)}^\sumIndB], \log^\diffeo_\Vector\VectorB) \log^\diffeo_\Vector\VectorB  \overset{\text{\cref{thm:pull-back-mappings} (ii)}}{=} \curvature^\diffeo(\diffeo^{-1}_*[\Theta_{\diffeo(\Vector)}^\sumIndB], \diffeo^{-1}_*[\log_{\diffeo(\Vector)} \diffeo(\VectorB)]) \diffeo^{-1}_*[\log_{\diffeo(\Vector)} \diffeo(\VectorB)] \\
     \overset{\cref{eq:pull-back-curvature}}{=} \diffeo^{-1}_* [R(\diffeo_* [\diffeo^{-1}_*[\Theta_{\diffeo(\Vector)}^\sumIndB]], \diffeo_* [\diffeo^{-1}_*[\log_{\diffeo(\Vector)} \diffeo(\VectorB)]])\diffeo_* [\diffeo^{-1}_*[ \log_{\diffeo(\Vector)} \diffeo(\VectorB)]] ] \\
     = \diffeo^{-1}_* [R(\Theta_{\diffeo(\Vector)}^\sumIndB, \log_{\diffeo(\Vector)} \diffeo(\VectorB)) \log_{\diffeo(\Vector)} \diffeo(\VectorB) ]
     = \kappa_\sumIndB \diffeo^{-1}_* [\Theta_{\diffeo(\Vector)}^\sumIndB].
\end{multline}
\end{proof}


\subsection{The curvature corrected Riemannian autoencoder}
\label{sec:cc-rae}

In this part we consider a particular type of compression, i.e., through low rank approximation, and see that the solution can be used to construct a non-linear projection mapping that has a non-linear encoder-decoder structure that reminisces of neural network-parameterized autoencoder architectures \cite{demers1992non,kingma2013auto}.

We will first consider low rank approximation itself. That is, we want to find an approximation that lives in the set
\begin{equation}
    \mathcal{S}_{\VectorC}^r ((\Real^{\dimInd})^\dataPointNum):= \{(\tangentVector_{\VectorC}^{1}, \ldots, \tangentVector_{\VectorC}^{\dataPointNum}) \in \tangent_{\VectorC} (\Real^{\dimInd})^\dataPointNum \; \mid \; \tangentVector_{\VectorC}^{\sumIndA} = \sum_{\sumIndC=1}^r \mathbf{U}_{\sumIndA,\sumIndC} \CCnetworkWeightVector^\sumIndC_{\VectorC}, \text{ where } \mathbf{U} \in \Real^{\dataPointNum \times r}, (\CCnetworkWeightVector^1_{\VectorC}, \ldots,\CCnetworkWeightVector^r_{\VectorC}) \in \tangent_{\VectorC} (\Real^{\dimInd})^r \}.
\end{equation}

In the following we will use the curvature-corrected singular value decomposition developed in \cite[Alg.~2]{diepeveen2023curvature} to get a cheap approximate minimiser of \cref{eq:curvature-corrected-approximation-problem}. That is, we first compute a minimiser to
\begin{equation}
     \inf_{(\tangentVector_{\VectorC}^{1}, \ldots, \tangentVector_{\VectorC}^{\dataPointNum})\in \mathcal{S}_{\VectorC}((\Real^{\dimInd})^\dataPointNum)} \Bigl\{ \sum_{\sumIndA=1}^N  (\|  \tangentVector_{\VectorC}^{\sumIndA} - \log^\diffeo_{\VectorC} (\Vector^\sumIndA)\|^\diffeo_{\VectorC})^2 \Bigr\},
    \label{eq:tangent-space-approximation-problem}
\end{equation}
and subsequently use it for approximating a minimiser of the curvature corrected error \cref{eq:curvature-corrected-approximation-problem}.


Truncated singular value decomposition in the tangent space solves \cref{eq:tangent-space-approximation-problem}. To compute this, consider any basis $\{\Phi_{\VectorC}^{\sumIndD}\}_{\sumIndD=1}^{\dimInd} \subset \tangent_\VectorC \Real^{\dimInd}$ that is orthonormal with respect to $(\cdot, \cdot)^\diffeo$ and define $\Matrix\in \Real^{\dataPointNum \times \dimInd}$ as
\begin{equation}
    \Matrix_{\sumIndA, \sumIndD} = (\log_\VectorC^\diffeo (\Vector^{\sumIndA}), \Phi_{\VectorC}^{\sumIndD})_\VectorC^\diffeo, \quad \text{for } \sumIndA = 1, \ldots, \dataPointNum, \; \sumIndD = 1, \ldots , \dimInd.
\end{equation}
Next, consider its the singular value decomposition (SVD)
\begin{equation}
    \Matrix = \mathbf{U} \Sigma \mathbf{W}^\top,
\end{equation}
where $\mathbf{U} \in \Real^{\dataPointNum \times R}$, $\Sigma = \operatorname{diag}(\sigma_1, \ldots, \sigma_{R})\in \Real^{R \times R}$ with $\sigma_1\geq \ldots \geq \sigma_{R}$, $\mathbf{W} \in \Real^{\dimInd \times R}$ and where $R := \operatorname{rank}(\Matrix)$. 
The low-rank approximation that is optimal in the tangent space error \cref{eq:tangent-space-approximation-problem}, is the one constructed from the first $r$ singular vectors. The key observation that is important to highlight is the following: if the data manifold is mapped to a low-dimensional geodesic subspace, we find by \cref{thm:pull-back-mappings} (i) that our data set has low rank. So we can really choose a low value for $r$ under such a diffeomorphism.

However, minimising \cref{eq:tangent-space-approximation-problem} will generally not minimise \cref{eq:curvature-corrected-approximation-problem}. Following the approach in \cite{diepeveen2023curvature}, we use the first $r$ rows of $\mathbf{U}$ from the tangent space SVD and correct the tangent vectors to get a good approximate minimiser. That is, our new approximation is
\begin{equation}
    \log^\diffeo_\VectorC (\Vector^\sumIndA) \approx \sum_{\sumIndC=1}^{r} \mathbf{U}_{\sumIndA,\sumIndC} \hat{\CCnetworkWeightVector}^{\sumIndC}_{\VectorC} := \sum_{\sumIndC=1}^{r} \mathbf{U}_{\sumIndA,\sumIndC} \sum_{\sumIndD=1}^{\dimInd}\hat{\mathbf{V}}_{\sumIndC,\sumIndD} \Phi^\sumIndD_{\VectorC},
    \label{eq:logxi-approx-low-rank}
\end{equation}
where
\begin{equation}
    \hat{\mathbf{V}} := \argmin_{\mathbf{V}\in \Real^{r \times \dimInd}} \Bigl\{ \sum_{\sumIndA=1}^N \sum_{\sumIndB=1}^{\dimInd} \betaExptVector (\kappa^\sumIndA_{\sumIndB})^2 \Bigl(\Bigl(  \sum_{\sumIndC,\sumIndD=1}^{r,\dimInd} \mathbf{U}_{\sumIndA,\sumIndC} \mathbf{V}_{\sumIndC,\sumIndD}\Phi^\sumIndD_{\VectorC} - \log^\diffeo_{\VectorC} (\Vector^\sumIndA), \Psi^{\sumIndA,\sumIndB}_{\VectorC} \Bigr)^\diffeo_{\VectorC}\Bigr)^2 \Bigr\},
    \label{eq:cc-coefs-diffeo-low-rank}
\end{equation}
which is well-posed if $\kappa^\sumIndA_{\sumIndB} < \pi^2$ for all $\sumIndA = 1, \ldots, \dataPointNum$ and $\sumIndB = 1, \ldots, \dimInd$ by \cite[Prop.~5.2]{diepeveen2023curvature} and can be solved in closed-form.


Subsequently, we can use the curvature corrected low rank approximation for a data-driven non-linear embedding. In particular, consider the $(\cdot, \cdot)^\diffeo$-orthonormal vectors $\{\CCnetworkWeightVector^{\sumIndC}_{\VectorC} \}_{\sumIndC=1}^r$ such that $\operatorname{span}(\{\CCnetworkWeightVector^{\sumIndC}_{\VectorC} \}_{\sumIndC=1}^r) = \operatorname{span}(\{\hat{\CCnetworkWeightVector}^{\sumIndC}_{\VectorC} \}_{\sumIndC=1}^r)$ -- which can be constructed through Gram-Schmidt orthogonalisation --, and consider the mappings $E:\Real^\dimInd \to \Real^r$ defined coordinate-wise as
\begin{equation}
    E(\Vector)_\sumIndC := (\log_{\VectorC}^{\diffeo} (\Vector), \CCnetworkWeightVector^\sumIndC_{\VectorC}
    )_\VectorC^\diffeo, \quad \sumIndC = 1, \ldots, r,
\end{equation}
and $D:\Real^r \to \Real^\dimInd$ defined as
\begin{equation}
    D(\latentVector):= \exp_{\VectorC}^{\diffeo} \Bigl( \sum_{\sumIndC=1}^r \latentVector_\sumIndC \CCnetworkWeightVector^\sumIndC_{\VectorC} \Bigr).
\end{equation}



Then, noting that we can rewrite
     \begin{equation}
        \sum_{\sumIndC=1}^{r} \mathbf{U}_{\sumIndA,\sumIndC} \hat{\CCnetworkWeightVector}^{\sumIndC}_{\VectorC} = \sum_{\sumIndC=1}^{r} \tilde{\mathbf{U}}_{\sumIndA,\sumIndC} \CCnetworkWeightVector^\sumIndC_{\VectorC}, \quad \text{for } \sumIndA = 1, \ldots, \dataPointNum,
        \label{eq:rem-change-of-coords-ccvectors}
    \end{equation}
     where $\tilde{\mathbf{U}} \in \Real^{\dataPointNum\times r}$ is defined as 
     \begin{equation}
         \tilde{\mathbf{U}}_{\sumIndA,\sumIndC}:= \sum_{\sumIndC'=1}^r \mathbf{U}_{\sumIndA,\sumIndC'} (\hat{\CCnetworkWeightVector}^{\sumIndC'}_{\VectorC}, \CCnetworkWeightVector^{\sumIndC}_{\VectorC})_\VectorC^\diffeo \quad \text{for }  \sumIndA = 1, \ldots, \dataPointNum, \; \sumIndC = 1, \ldots , r,
     \end{equation}
     we now know that the value $\| D(E(\Vector^\sumIndA)) - \Vector^\sumIndA\|_2^2$ is small for every $\sumIndA=1, \ldots, \dataPointNum$, since
\begin{equation}
     \| D(E(\Vector^\sumIndA)) - \Vector^\sumIndA\|_2^2 = \Bigl\| \exp_{\VectorC}^\diffeo \Bigl(  \sum_{\sumIndC=1}^r (\log_{\VectorC}^{\diffeo} (\Vector^\sumIndA), \CCnetworkWeightVector^\sumIndC_{\VectorC}
    )_\VectorC^\diffeo  \CCnetworkWeightVector^\sumIndC_{\VectorC}\Bigr) - \Vector^\sumIndA\Bigr\|_2^2 \overset{\text{\cref{eq:logxi-approx-low-rank,eq:rem-change-of-coords-ccvectors}}}{\approx} \Bigl\| \exp_{\VectorC}^\diffeo \Bigl(  \sum_{\sumIndC=1}^r \tilde{\mathbf{U}}_{\sumIndA,\sumIndC}  \CCnetworkWeightVector^\sumIndC_{\VectorC}\Bigr) - \Vector^\sumIndA\Bigr\|_2^2,
    \label{eq:RAE-error-on-data}
\end{equation}
and the right-hand side is small because the $\tilde{\mathbf{U}} \in \Real^{\dataPointNum \times r}$ and $\{\CCnetworkWeightVector^{\sumIndC}_{\VectorC} \}_{\sumIndC=1}^r \subset \tangent_{\VectorC} (\Real^{\dimInd})$ are constructed to approximately minimise \cref{eq:global-approximation-problem-l2}. 

In other words, a suitable diffeomorphism combined with a low rank approximation scheme on manifolds gives rise to a natural non-linear encoder-decoder, which we call the \emph{curvature corrected Riemannian autoencoder} (CC-RAE) for the case where curvature corrected singular value decomposition \cite[Alg.~2]{diepeveen2023curvature} is used for low rank approximation. If the curvature correction step is skipped, i.e., we take $\CCnetworkWeightVector^\sumIndC_{\VectorC} := \sum_{\sumIndD=1}^\dimInd \mathbf{W}_{\sumIndD\sumIndC} \Phi_{\VectorC}^{\sumIndD}$, we call the non-linear encoder-decoder a \emph{Riemannian autoencoder} (RAE). Choosing an RAE over a CC-RAE can be admissible if the curvature effects are small, because then the minimisers of \cref{eq:tangent-space-approximation-problem} and \cref{eq:curvature-corrected-approximation-problem} are close to each other \cite[Thm. 4.6]{diepeveen2023curvature}.


\begin{remark}
\label{rem:rae-interpretability}
    Our analysis has already suggested that $\diffeo$ should be a local isomorphism on the data set. The additional upshot of local isometry (combined with mapping the data manifold into a geodesic subspace) for the CC-RAE is that the latent space has meaning. That is, local and global geometry and distances in the latent space $\Real^r$ correspond to local and global geometry and distances in the actual data space (up to curvature). So if we have a suitable diffeomorphism, the CC-RAE does not suffer from the interpretability issues that classical neural network-based autoencoder architectures \cite{demers1992non,kingma2013auto} have.
\end{remark}
\section{Learning suitable diffeomorphisms}
\label{sec:learning-diffeo}

From \cref{sec:basic-processing-pullback,sec:non-linear-compression-pullback} we have learned  that diffeomorphisms $\diffeo:\Real^\dimInd \to \manifold$ need to map the data manifold into a geodesic subspace of $(\manifold, (\cdot, \cdot))$ while preserving local isometry around the data manifold for proper, stable and efficient data analysis on $(\Real^\dimInd, (\cdot,\cdot)^\diffeo)$. What we have not addressed so far is whether such mappings exist in the first place and if they exist, how to construct them. Although a rigorous answer to the question of existence is open and beyond the scope of this article, we can use insights from empirical successes for a heuristic construction strategy. 

We will use a variational approach to constructing diffeomorphisms from a data set  $\{\Vector^\sumIndA\}_{\sumIndA=1}^\dataPointNum \subset \Real^\dimInd$. That is, we choose a class of suitable diffeomorphisms and subsequently try to find a best one in some energy. We propose to consider diffeomorphisms $\diffeo: \Real^\dimInd \to \manifold_{\dimInd'} \times \Real^{\dimInd - \dimInd'}$ of the form
\begin{equation}
    \diffeo := (\diffeoC^{-1}, \mathbf{I}_{\dimInd-\dimInd'}) \circ \diffeoB \circ \mathbf{O} \circ T_{\VectorC},
    \label{eq:diffeo-class}
\end{equation}
where $\diffeoC: \mathcal{U} \to \Real^{\dimInd'}$ is a chart on a (geodesically convex) subset $\mathcal{U} \subset \manifold_{\dimInd'}$ of a $\dimInd'$-dimensional Riemannian manifold $(\manifold_{\dimInd'}, (\cdot, \cdot)')$, $\diffeoB: \Real^{\dimInd}\to \Real^{\dimInd}$ is a real-valued diffeomorphism, $\mathbf{O} \in \mathbb{O}(\dimInd)$ is a orthogonal matrix, and $T_{\VectorC}: \Real^{\dimInd}\to \Real^{\dimInd}$ is given by $T_{\VectorC} (\Vector) = \Vector - \VectorC$. 

The motivation for the above class of diffeomorphisms is easiest understood by considering what parts can be modeled and what should typically be learned. We first note that the low-dimensional manifold $\manifold_{\dimInd'}$ is the space we want to map the data manifold into. This space can be chosen by finding the non-linear embedding manifold among several candidate manifolds (see \cref{sec:intro-pull-back-geometry}) that has minimal metric distortion with respect to the approximate geodesic distances $\{\distance_{\sumIndA,\sumIndA'} \}_{\sumIndA,\sumIndA'=1}^\dataPointNum$ between all pairs of points $\Vector^\sumIndA$ and $\Vector^{\sumIndA'}$ -- constructed from completion of local Euclidean distances to neighbouring data points \cite{tenenbaum2000global} or something more sophisticated, e.g., that still preserves isometry \cite{budninskiy2019parallel} or gives more noise-robust geodesics \cite{little2022balancing}. Next, assuming that the data manifold and $\manifold_{\dimInd'}$ are actually homeomorphic\footnote{This is not always true if we can find a non-linear embedding.}, the second coordinates need to be mapped to zero, i.e., we want 
\begin{equation}
    \pi_2(\diffeo(\Vector^\sumIndA)) = \mathbf{0}_{\dimInd - \dimInd'} \in \Real^{\dimInd - \dimInd'}, \quad \text{for all }\sumIndA=1, \ldots, \dataPointNum 
    \label{eq:diffeo-mapping-goal}
\end{equation}
To accomplish \cref{eq:diffeo-mapping-goal} locally, assume for the moment that $\diffeoB = \mathbf{I}_{\dimInd}$ and consider the case that the vector $\VectorC$ in \cref{eq:diffeo-class} is a data point, i.e., $\VectorC = \Vector^{\sumIndA'}$ for some $\sumIndA' \in [\dataPointNum]$. 
Then, for any $\mathbf{O}$ we have
\begin{equation}
    \pi_2(\diffeo(\Vector^{\sumIndA'})) = \pi_2(((\diffeoC^{-1}, \mathbf{I}_{\dimInd-\dimInd'}) \circ \mathbf{I}_{\dimInd} \circ \mathbf{O} \circ T_{\VectorC})(\Vector^{\sumIndA'})) = \pi_2(((\diffeoC^{-1}, \mathbf{I}_{\dimInd-\dimInd'}) \circ \mathbf{I}_{\dimInd} \circ \mathbf{O})(\mathbf{0}_{\dimInd}))  = \mathbf{0}_{\dimInd- \dimInd'} \in \Real^{\dimInd - \dimInd'},
\end{equation}
but we do not necessarily get zero for data points close to $\Vector^{\sumIndA'}$. We can accomplish the latter -- at least approximately -- by choosing a particular matrix $\mathbf{O}$ through (local) principal component analysis. That is, for some $r>0$ we compute
\begin{equation}
    \sum_{\Vector \in \{\Vector^\sumIndA\}_{\sumIndA=1}^\dataPointNum\cap \ball_{r} (\Vector^{\sumIndA'})} (\Vector - \Vector^{\sumIndA'}) \otimes (\Vector - \Vector^{\sumIndA'}) = \mathbf{U} \Lambda \mathbf{U}^\top, 
\end{equation}
where $\mathbf{U}  \in \mathbb{O}(\dimInd)$ and $\Lambda = \operatorname{diag} (\lambda_1, \ldots, \lambda_\dimInd) \in \Real^{\dimInd\times \dimInd}$ with $\lambda_1 \geq \ldots \geq \lambda_\dimInd$, and set $\mathbf{O} := \mathbf{U}^\top$. Finally, to accomplish \cref{eq:diffeo-mapping-goal} globally, we drop the assumption that $\diffeoB = \mathbf{I}_{\dimInd}$ and realize that we need a $\diffeoB$ that is approximately identity close to $\Vector^{\sumIndA'}$, but deviates from identity at the points farther away from $\Vector^{\sumIndA'}$. For that we can choose $\diffeoB:= \diffeoB_\networkParams$ to be a neural network such as \cite{behrmann2019invertible,gruffaz2021learning,salman2018deep,younes2018diffeomorphic} with parameters $\networkParams \in \networkParamsSet$ in some parameter set $\networkParamsSet$.






Then, for a suitable mapping $\diffeoC: \mathcal{U} \to \Real^{\dimInd'}$, fixed $\VectorC \in \{\Vector^\sumIndA\}_{\sumIndA=1}^\dataPointNum $ and corresponding orthogonal matrix $\mathbf{O}$ obtained from a (local) PCA around $\VectorC$, we propose to find an optimal diffeomorphism $\diffeo_\networkParams$ of the form \cref{eq:diffeo-class} with learnable $\diffeoB:= \diffeoB_\networkParams$ through solving the minimisation problem
\begin{equation}
    \inf_{\networkParams \in \networkParamsSet} \Bigl\{ \frac{1}{\dataPointNum(\dataPointNum-1)}\sum_{\substack{\sumIndA, \sumIndB = 1, \\ \sumIndA \neq \sumIndA'}}^{\dataPointNum(\dataPointNum-1)} \Bigl(  \distance_{\Real^\dimInd}^{\diffeo_\networkParams}(\Vector^{\sumIndA}, \Vector^{\sumIndA'}) - 
    \distance_{\sumIndA,\sumIndA'}\Bigr)^2  + \alpha_{\mathrm{sub}} \frac{1}{\dataPointNum}\sum_{\sumIndA=1}^\dataPointNum \|\pi_2 (\diffeo_\networkParams(\Vector^{\sumIndA}))\|_1 + \alpha_{\mathrm{iso}} \frac{1}{\dataPointNum}\sum_{\sumIndA=1}^\dataPointNum \| \bigl( (\mathbf{e}^\sumIndB, \mathbf{e}^{\sumIndB'})_{\Vector^\sumIndA}^{\diffeo_\networkParams} \bigr)_{\sumIndB,\sumIndB'=1}^\dimInd- \mathbf{I}_\dimInd \|_F^2\Bigr\},
    \label{eq:pwd-learning-problem}
\end{equation}
with $\diffeoB_\networkParams$ initialized close to identity, where $\distance_{\Real^\dimInd}^{\diffeo_\networkParams}: \Real^\dimInd \times \Real^\dimInd \to \Real$ is the distance on $(\Real^\dimInd, (\cdot, \cdot)^{\diffeo_\networkParams})$ (see \cref{thm:pull-back-mappings} (iv)) and where $\alpha_{\mathrm{sub}},\alpha_{\mathrm{iso}} >0$. The rationale behind each of the three terms comes down to the following. \emph{The first term takes global geometry into account} and is motivated by the success of Siamese networks \cite{bromley1993signature,pai2019dimal} for learning isometric embeddings and by the theoretical observation that global distances are important in learning embeddings \cite{fefferman2020reconstruction}. \emph{The second term enforces that the data manifold is mapped to $\manifold_{\dimInd'}$}, which is motivated by our analysis. Finally, \emph{the third term enforces local isometry}, which is also motivated by our analysis.





\section{Numerics}
\label{sec:numerics-pullback-manifolds}

In this section we test the suitability of pullback Riemannian geometry for data analysis tasks from \cref{sec:basic-processing-pullback,sec:non-linear-compression-pullback}. Remember from \cref{sec:intro-pull-back-geometry} that ideal data analysis tools should be able to \emph{interpolate and extrapolate over non-linear paths through the data, compute non-linear means on such non-linear paths, and perform low-rank approximation over curved subspaces spanned by such non-linear paths} as argued in \cite{diepeveen2023riemannian}. So naturally, we will test how well the Riemannian interpretation of these tools on $(\Real^\dimInd, (\cdot,\cdot)^\diffeo)$ -- geodesic interpolation, the barycentre as a non-linear mean, and tangent space low-rank approximation from any data point using the logarithmic mapping -- capture the data geometry for a properly chosen diffeomorphism $\diffeo:\Real^\dimInd \to \manifold$ and and Riemannian geometry of $(\manifold, (\cdot, \cdot))$. In particular, we want to get insight into (i) the role of the Riemannian geometry of $(\manifold, (\cdot, \cdot))$, and see (ii) to what extent the construction in \cref{sec:learning-diffeo} is suitable for constructing a diffeomorphism $\diffeo:\Real^\dimInd \to \manifold$.

\begin{figure}[h!]
    \centering
    \begin{subfigure}{0.31\linewidth}
    \centering
        \includegraphics[width=\linewidth]{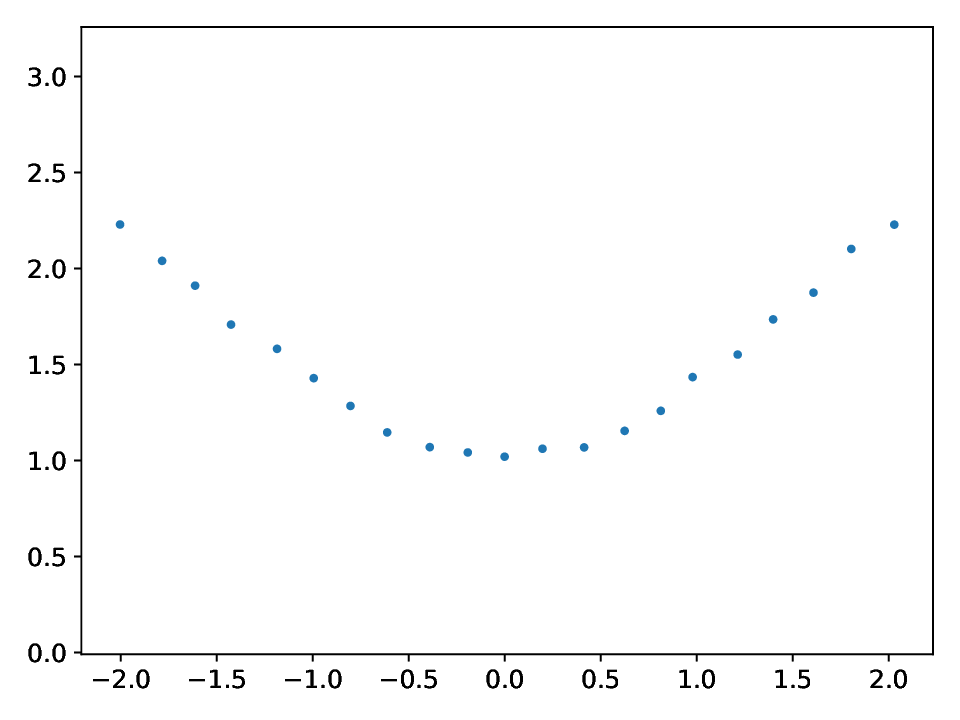}
    \caption{}
    \label{fig:hyperbolic-example-data}
    \end{subfigure}
    \hfill
    \begin{subfigure}{0.31\linewidth}
    \centering
        \includegraphics[width=\linewidth]{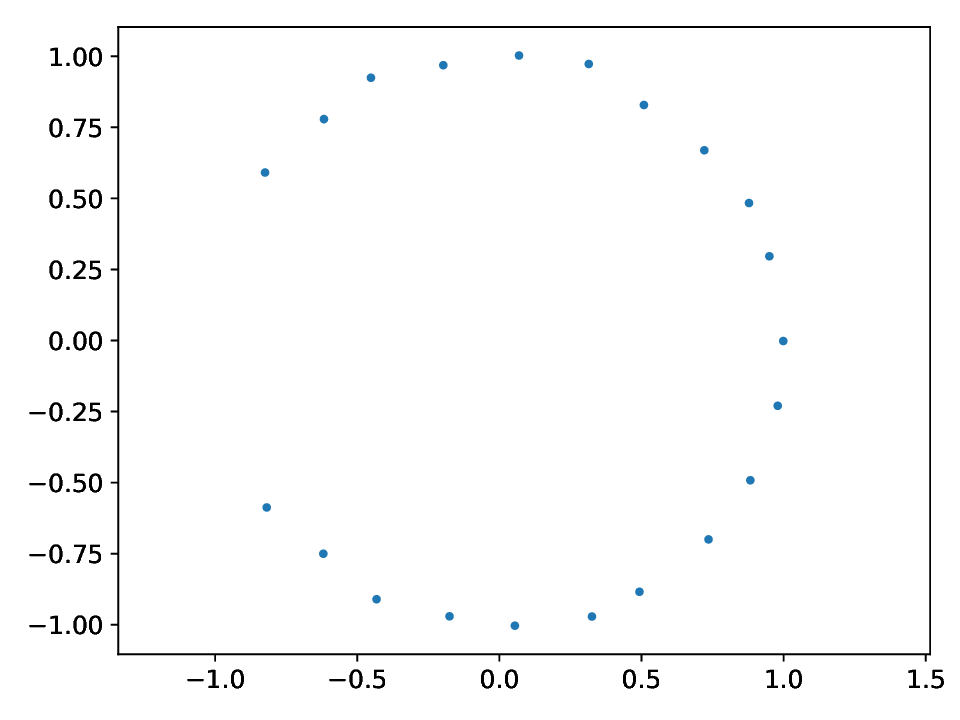}
    \caption{}
    \label{fig:cirlc-example-data}
    \end{subfigure}
    \hfill
    \begin{subfigure}{0.31\linewidth}
    \centering
        \includegraphics[width=\linewidth]{experiments/diffeo-effects/s_data.eps}
    \caption{}
    \label{fig:spiral-example-data}
    \end{subfigure}
    \caption{Three toy data sets.}
    \label{fig:data-for-pullback-geometry}
\end{figure}

\paragraph{Data.}
For numerical evaluation of (i) and (ii), data analysis under pullback geometry can be interpreted best if we can visualize geodesics, barycentres, low-rank approximations and their errors. For that reason we consider the $\Real^2$-valued toy data sets (a), (b) and (c) in \cref{fig:data-for-pullback-geometry} that have underlying dimension 1.


\paragraph{Outline of experiments.}
In the following we will showcase that data sets (a) and (b) can be analysed under a modelled hyperbolic and spherical pullback geometry, and that the data set (c) can be analysed under a learned Euclidean pullback geometry. Then, considering the differences between (a) and (b) in \cref{sec:numerics-pullback-manifolds-curvature-effects} will give us insight into (i) and considering (c) for different ways of learning the diffeomorphism -- with out without the best practices predicted by theory -- in \cref{sec:numerics-pullback-manifolds-diffeo-effects} will give us insight into (ii). It should be noted that for data sets (a) and (b) we can also learn a Euclidean pullback geometry -- since both data sets can be embedded into 1-dimensional Euclidean space --, but for now that is not the goal of these experiments.




\paragraph{General experimental settings.}
Throughout \cref{sec:numerics-pullback-manifolds-curvature-effects,sec:numerics-pullback-manifolds-diffeo-effects}, the manifolds mappings in \cref{thm:pull-back-mappings} and combinations thereof are used for analysing data sets (a), (b) and (c) under their respective pullback geometries. In addition, Riemannian barycentres as defined in \cref{eq:diffeo-barycentre-formal} are computed in the embedding manifold and mapped back as suggested in \cref{eq:diffeo-barycentre} by \cref{thm:well-posedness-barycentre}. For the barycentre problem in the embedding manifold, Riemannian gradient descent with unit step size is used to solve it and is terminated when the relative Riemannian gradient at iteration $\ell$ satisfies
\begin{equation}
    \frac{\|\Grad \bigl(\frac{1}{2\dataPointNum} \sum_{\sumIndA = 1}^\dataPointNum \distance_{\manifold}(\cdot, \diffeo(\Vector^\sumIndA))^2\bigr)\|_{\mPoint^{\ell}} }{\|\Grad\bigl( \frac{1}{2\dataPointNum} \sum_{\sumIndA = 1}^\dataPointNum \distance_{\manifold}(\cdot, \diffeo(\Vector^\sumIndA))^2 \bigr) \|_{\mPoint^{0}}} = \frac{\|\frac{1}{\dataPointNum} \sum_{\sumIndA = 1}^\dataPointNum \log_{\mPoint^{\ell}} (\diffeo(\Vector^\sumIndA))\|_{\mPoint^{\ell}}}{ \|\frac{1}{\dataPointNum} \sum_{\sumIndA = 1}^\dataPointNum \log_{\mPoint^{0}} (\diffeo(\Vector^\sumIndA)) \|_{\mPoint^{0}} } < 10^{-3}.
\end{equation}
The exception is the case of Euclidean pullback geometry, where we have a closed-form solution \cref{eq:euclidean-diffeo-mean} due to \cref{cor:euclidean-diffeo-mean}.



Finally, all of the experiments are implemented using \texttt{PyTorch} in Python 3.8 and run on a 2 GHz Quad-Core Intel Core i5 with 16GB RAM. 

\subsection{Curvature effects}
\label{sec:numerics-pullback-manifolds-curvature-effects}

For the first experiment we choose specific pullback geometries under which we perform data analysis tasks from \cref{sec:basic-processing-pullback,sec:non-linear-compression-pullback} on data sets (a) and (b) to get insight into (i), i.e., the role of the pulled back geometry. In particular, under pullback geometry with positive and negative curvature we will consider stability of geodesics and barycentres and consider low rank approximation whose quality is tested through constructing Riemannian autoencoders\footnote{For both data sets there is no deterioration due to curvature in the low rank approximation. So we can skip the curvature correction step outlined in \cref{sec:cc-rae} and get a good minimiser from a basic truncated singular value decomposition in the tangent space.}.
 


\subsubsection{Hyperbolic pullback}
\label{sec:numerics-pullback-manifolds-curvature-effects-hyperbolic}
For analysing data set (a), we consider a fully modelled diffeomorphism of the form \cref{eq:diffeo-class} into the unit hyperboloid and pull back the standard hyperbolic Riemannian structure. That is, we consider $\diffeo^{\mathrm{a}}: \Real^2 \to \Hyperboloid^{2}$ with $\VectorC^{\mathrm{a}} := (0, 1) \in \Real^2$, $\mathbf{O}^{\mathrm{a}}:= \mathbf{I}_2 \in \Real^{2\times 2}$, $\diffeoB^{\mathrm{a}} :\Real^{2} \to \Real^{2}$ being identity, and $\diffeoC^{\mathrm{a}}: \Hyperboloid^{2}\to \Real^2$ being the diffeomorphism \cref{eq:hyperboloid-chart} in \cref{app:pullback-manifolds-numerics-hyperboloid}, and pull back the Minkowski inner product on $\Real^{\dimInd+1}$. Then, our data analysis is carried out on $(\Real^2, (\cdot,\cdot)^{\diffeo^{\mathrm{a}}})$.

\paragraph{Interpolation.}

\begin{figure}[h!]
    \centering
    \begin{subfigure}{0.48\linewidth}
    \centering
        \includegraphics[width=\linewidth]{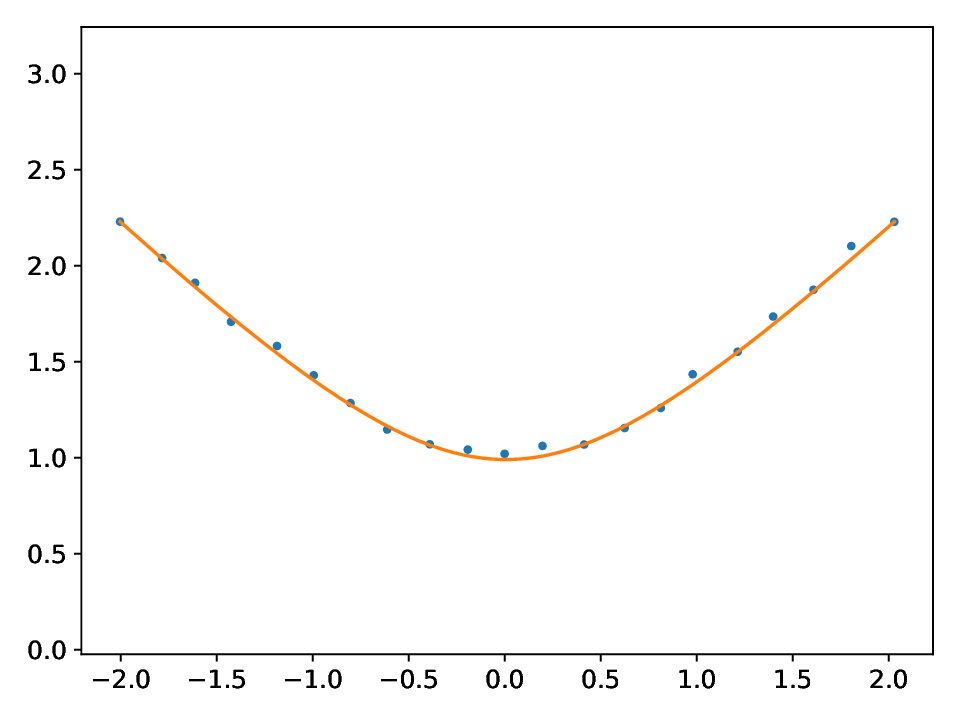}
    \caption{Geodesic interpolation.}
    \label{fig:hyperbolic-interpolation-example-geodesic}
    \end{subfigure}
    \hfill
    \begin{subfigure}{0.48\linewidth}
    \centering
        \includegraphics[width=\linewidth]{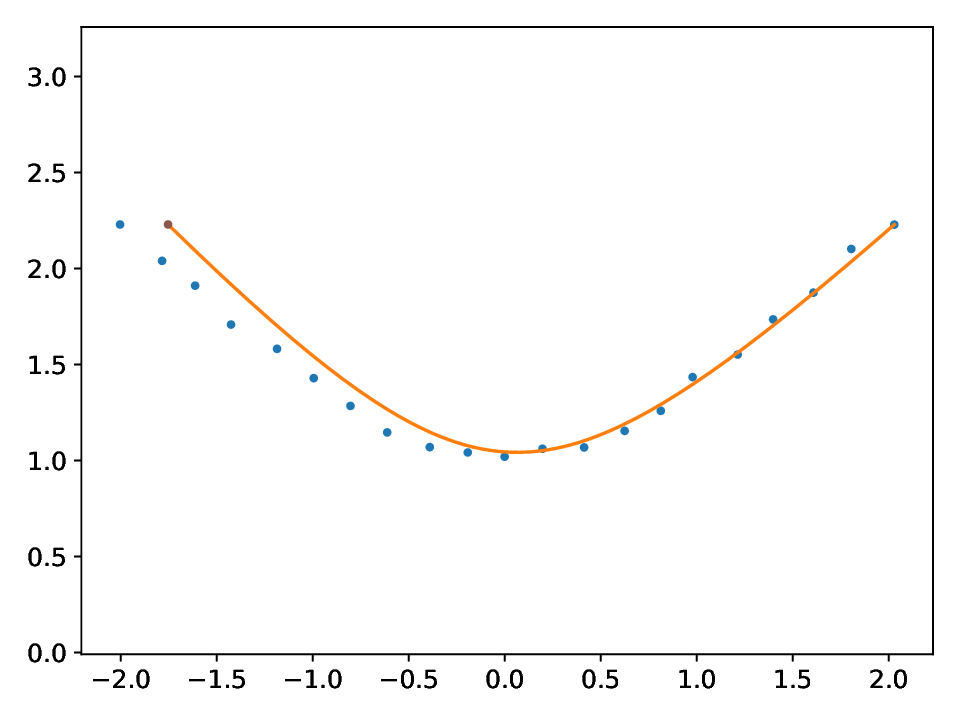}
    \caption{Perturbed geodesic interpolation.}
    \label{fig:hyperbolic-interpolation-example-variation}
    \end{subfigure}
    \caption{Geodesic interpolation and perturbed geodesic interpolation of the end points of data set (a) on $(\Real^2, (\cdot,\cdot)^{\diffeo^{\mathrm{a}}})$ indicates that the chosen pullback geometry is suitable for data analysis on this data set and is stable with respect to small perturbations.}
    \label{fig:hyperbolic-interpolation-example}
\end{figure}

First, we consider interpolation and its stability through interpolating between the end points of the data set (a) and vary one of the end points with an out of distribution point, see \cref{fig:hyperbolic-interpolation-example}. \Cref{fig:hyperbolic-interpolation-example-geodesic} showcases that the pulled back geometry does exactly what it should do. That is, geodesics go straight through the data set. \Cref{thm:stability-geodesics-symmetric} suggests that for negative curvature we expect stability with respect to changing end points, which is also exactly what is observed in \cref{fig:hyperbolic-interpolation-example-variation}.

\begin{figure}[h!]
    \centering
    \begin{subfigure}{0.48\linewidth}
    \centering
        \includegraphics[width=\linewidth]{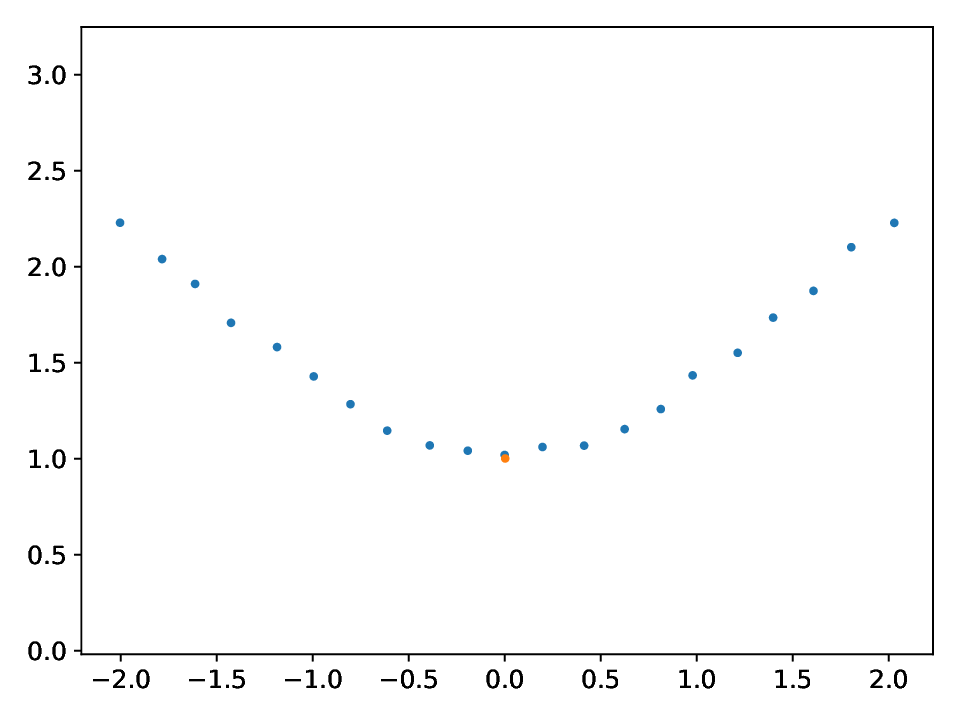}
    \caption{The data barycentre.}
    \label{fig:hyperbolic-barycentre-example-standard}
    \end{subfigure}
    \hfill
    \begin{subfigure}{0.48\linewidth}
    \centering
        \includegraphics[width=\linewidth]{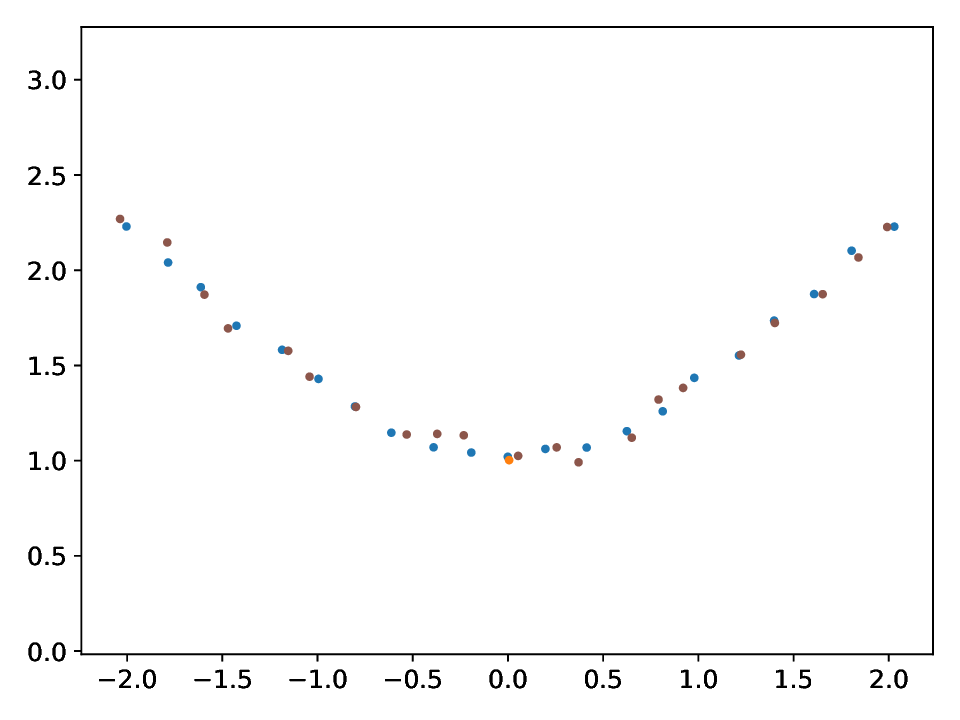}
    \caption{The perturbed data barycentre.}
    \label{fig:hyperbolic-barycentre-example-variation}
    \end{subfigure}
    \caption{The data barycentre and perturbed data barycentre of data set (a) on $(\Real^2, (\cdot,\cdot)^{\diffeo^{\mathrm{a}}})$ indicates that the chosen pullback geometry is suitable for data analysis on this data set and is stable with respect to small perturbations.}
    \label{fig:hyperbolic-barycentre-example}
\end{figure}

\paragraph{Barycentre.} Next, we consider the Riemannian barycentre and test its stability through small perturbations of the data, see \cref{fig:hyperbolic-barycentre-example}. Similarly as before, \cref{fig:hyperbolic-barycentre-example-standard} shows that the Riemannian barycentre is right were it is expected to be, i.e., within the data set and close to $\VectorC^{\mathrm{a}}$ due to the symmetry of the data around that point. Small variations of the data do not give rise to instabilities either as shown in \cref{fig:hyperbolic-barycentre-example-variation}. This is again in line with expectations due to \cref{thm:stability-barycentre}.

\begin{figure}[h!]
    \centering
    \begin{subfigure}{0.48\linewidth}
    \centering
        \includegraphics[width=\linewidth]{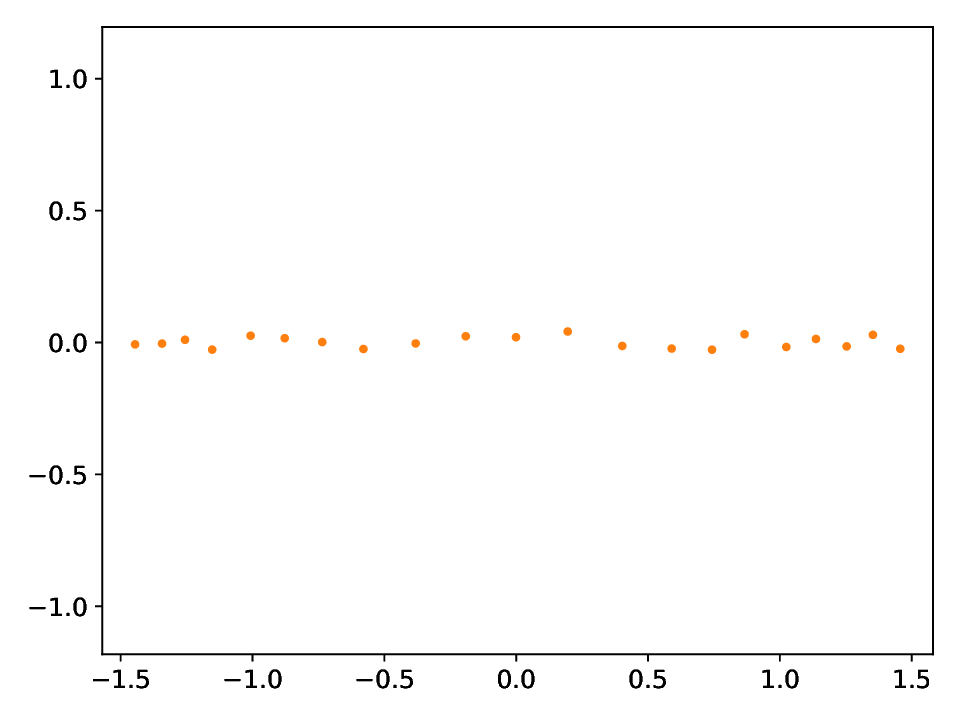}
    \caption{Data set (a) seen from the tangent space $\tangent_{\VectorC^{\mathrm{a}}} \Real^2$.}
    \label{fig:hyperbolic-rae-example-logs}
    \end{subfigure}
    \hfill
    \begin{subfigure}{0.48\linewidth}
    \centering
        \includegraphics[width=\linewidth]{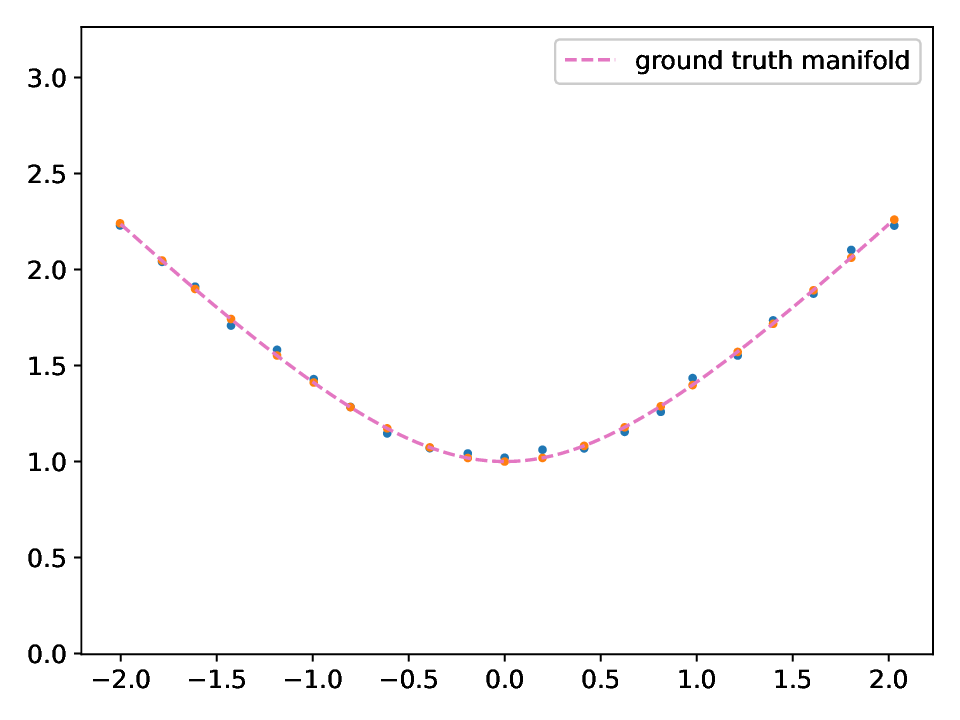}
    \caption{Data set (a) approximated by the RAE at $\VectorC^{\mathrm{a}}$.}
    \label{fig:hyperbolic-rae-example-rae}
    \end{subfigure}
    \caption{Low rank approximation of data set (a) on $(\Real^2, (\cdot,\cdot)^{\diffeo^{\mathrm{a}}})$ and the Riemannian autoencoder constructed from it indicate that the chosen pullback geometry is suitable for data analysis on this data set.}
    \label{fig:hyperbolic-rae-example}
\end{figure}

\paragraph{Riemannian autoencoder.} Finally, we compute the logarithmic mappings from $\VectorC^{\mathrm{a}}$ to all of the data, do a low rank approximation and construct a RAE, which is used to project the original data onto the learned manifold. Both results are shown in \cref{fig:hyperbolic-rae-example}. Once again, the results are by and large in line with expectations. That is, the data look 1-dimensional and linear on the tangent space (\cref{fig:hyperbolic-rae-example-logs}) and the RAE finds the correct data manifold (\cref{fig:hyperbolic-rae-example-rae}). Surprisingly, we retrieve the correct manifold without curvature correction, which can make a big difference \cite{diepeveen2023curvature}.

\paragraph{Conclusions.} Negative curvature evidently seems to give rise to stability for interpolation and barycentres and its known challenges for low rank approximation -- although not being an issue here -- can in general can be accounted for when passing to a CC-RAE. Overall, the above numerical experiments are in line with the theory in the sense that they both suggest that pulling back negatively curved spaces can be beneficial for downstream data analysis.

\subsubsection{Spherical pullback}
\label{sec:numerics-pullback-manifolds-curvature-effects-circle}
For analysing data set (b), we consider a fully modelled diffeomorphism of the form \cref{eq:diffeo-class} into the unit sphere and pull back the standard spherical Riemannian structure. That is, we consider $\diffeo^{\mathrm{b}}: \Real^2 \to \Sphere^{2}$ with $\VectorC^{\mathrm{b}} := (1, 0) \in \Real^2$, $\mathbf{O}^{\mathrm{b}}:= \mathbf{I}_2 \in \Real^{2\times 2}$, $\diffeoB^{\mathrm{b}} :\Real^{2} \to \Real^{2}$ being identity, and $\diffeoC^{\mathrm{b}}: \Sphere^{2}\setminus\{(0, 0, 1)\}  \to \Real^2$ being the diffeomorphism \cref{eq:sphere-chart} in \cref{app:pullback-manifolds-numerics-sphere}, and pull back the Euclidean inner product on $\Real^{\dimInd+1}$. Then, our data analysis is carried out on $(\Real^2, (\cdot,\cdot)^{\diffeo^{\mathrm{b}}})$.

\begin{figure}[h!]
    \centering
    \begin{subfigure}{0.48\linewidth}
    \centering
        \includegraphics[width=\linewidth]{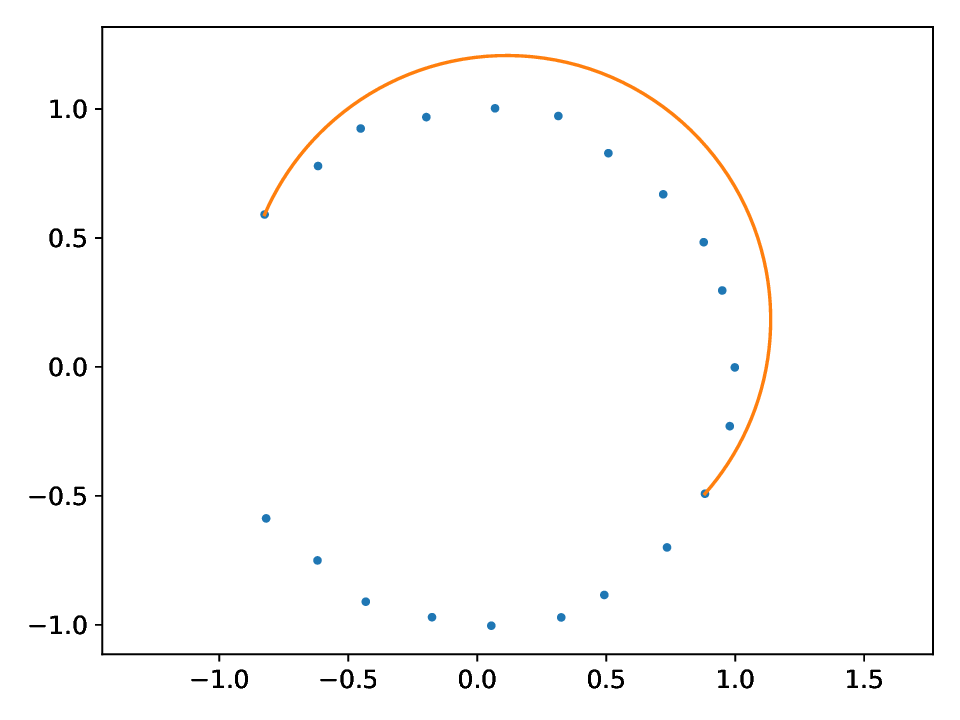}
    \caption{Geodesic interpolation.}
    \label{fig:circle-interpolation-example-geodesic}
    \end{subfigure}
    \hfill
    \begin{subfigure}{0.48\linewidth}
    \centering
        \includegraphics[width=\linewidth]{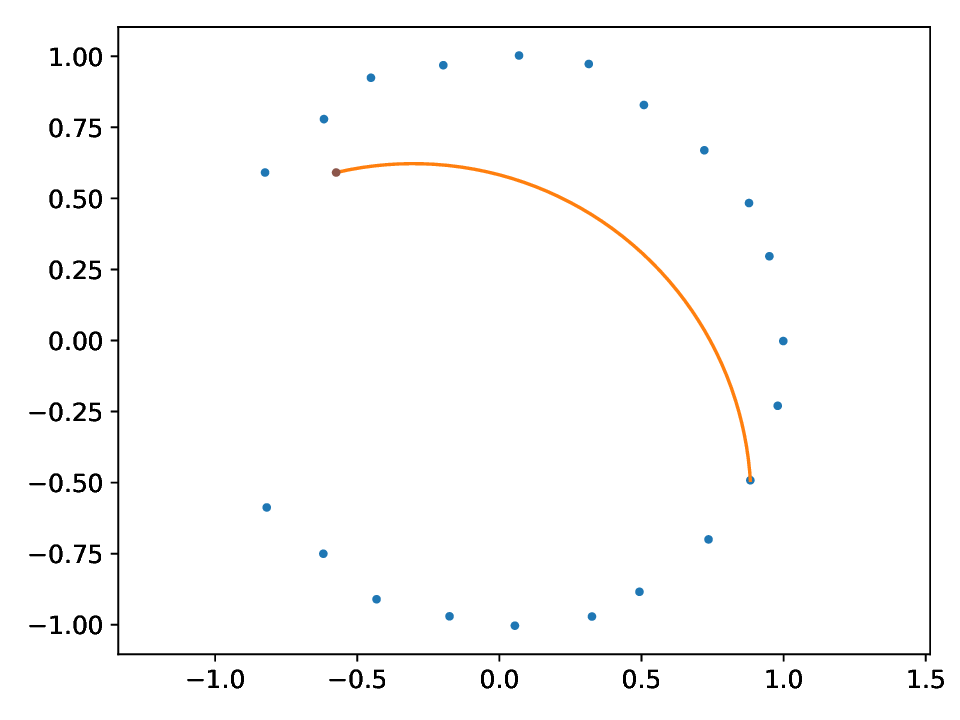}
    \caption{Perturbed geodesic interpolation.}
    \label{fig:circle-interpolation-example-variation}
    \end{subfigure}
    \caption{Geodesic interpolation and perturbed geodesic interpolation of the end points of data set (b) on $(\Real^2, (\cdot,\cdot)^{\diffeo^{\mathrm{b}}})$ indicates that the chosen pullback geometry is suitable for data analysis on this data set, but is unstable with respect to small perturbations.}
    \label{fig:circle-interpolation-example}
\end{figure}

\paragraph{Interpolation.}
Once again, we consider interpolation and its stability through interpolating between two points of the data set (b) and vary one of the end points with an out of distribution point, see \cref{fig:circle-interpolation-example}. Now, \cref{fig:circle-interpolation-example-geodesic} showcases that the pulled back geometry does almost what it should do. That is, the geodesic passes relatively close along the data set, but does not go through it like in the hyperbolic case in \cref{sec:numerics-pullback-manifolds-curvature-effects-hyperbolic}. The reason that we do not go exactly through the data set can be contributed to the positive curvature of the pulled back geometry that amplifies the small noise in the data (\cref{thm:stability-geodesics-symmetric}). \Cref{fig:circle-interpolation-example-variation} is in line with this suspicion, because for a change of end points we get entirely different geodesics.

\begin{figure}[h!]
    \centering
    \begin{subfigure}{0.48\linewidth}
    \centering
        \includegraphics[width=\linewidth]{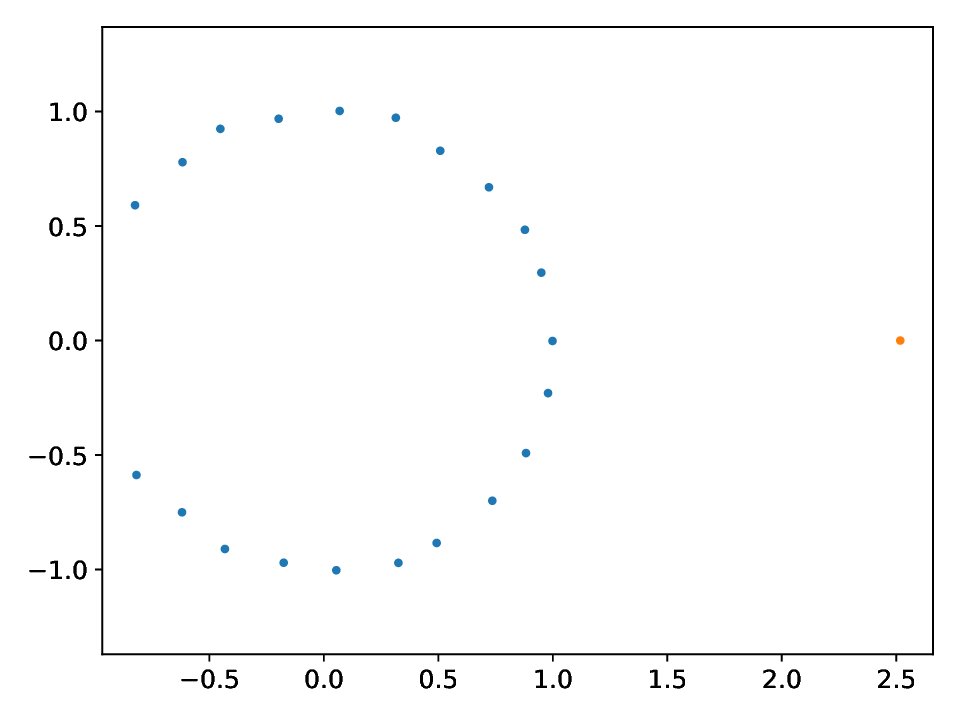}
    \caption{The data barycentre.}
    \label{fig:circle-barycentre-example-standard}
    \end{subfigure}
    \hfill
    \begin{subfigure}{0.48\linewidth}
    \centering
        \includegraphics[width=\linewidth]{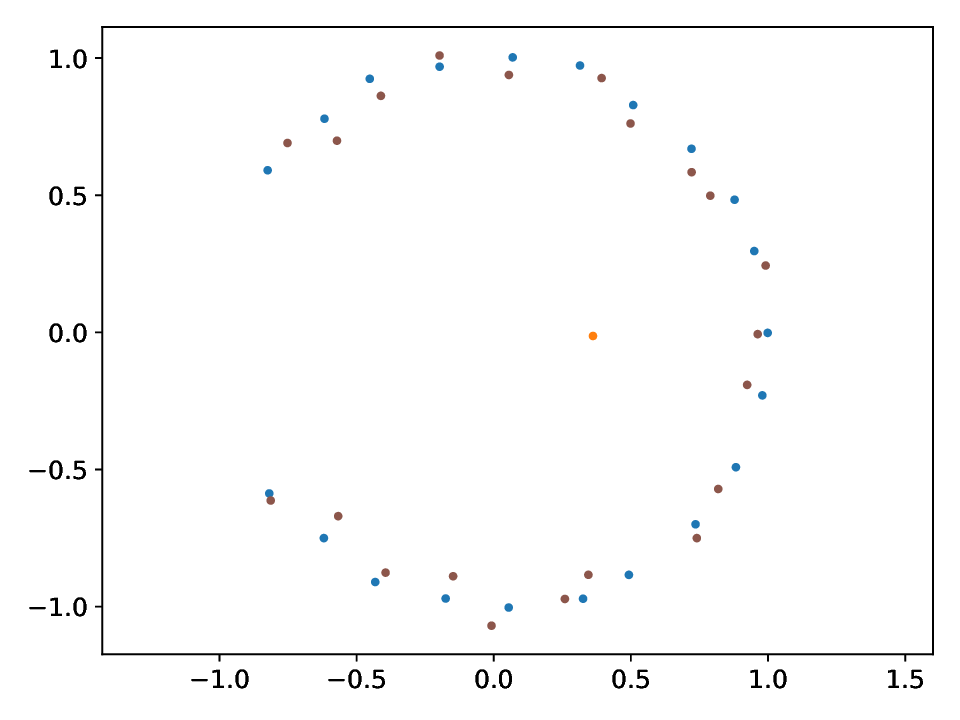}
    \caption{The perturbed data barycentre.}
    \label{fig:circle-barycentre-example-variation}
    \end{subfigure}
    \caption{The data barycentre and perturbed data barycentre of data set (b) on $(\Real^2, (\cdot,\cdot)^{\diffeo^{\mathrm{b}}})$ indicates that the chosen pullback geometry is unstable with respect to small perturbations.}
    \label{fig:circle-barycentre-example}
\end{figure}

\paragraph{Barycentre.}
Next, we also consider the Riemannian barycentre and test its stability through small perturbations of the data, see \cref{fig:circle-barycentre-example}. Contrary to the hyperbolic case in \cref{sec:numerics-pullback-manifolds-curvature-effects-hyperbolic}, \cref{fig:circle-barycentre-example-standard} shows that the Riemannian barycentre is far away from were it is expected to be, i.e., within the data set and close to $\VectorC^{\mathrm{b}}$ due to the symmetry of the data around that point. Similarly to interpolation, we expect that this is caused by the positive curvature \cref{thm:stability-barycentre}. This suspicion is once again backed up when considering small variations of the data, which also give rise to instabilities as shown in \cref{fig:circle-barycentre-example-variation}, where the perturbed data barycentre is very far away from the original one\footnote{The Riemannian gradient descent scheme for computing barycentre of the perturbed data is initialized from the original data barycentre to make sure that differences cannot be contributed to different local minima due to the non-convexity of the barycentre problem on the sphere.}.

\begin{figure}[h!]
    \centering
    \begin{subfigure}{0.48\linewidth}
    \centering
        \includegraphics[width=\linewidth]{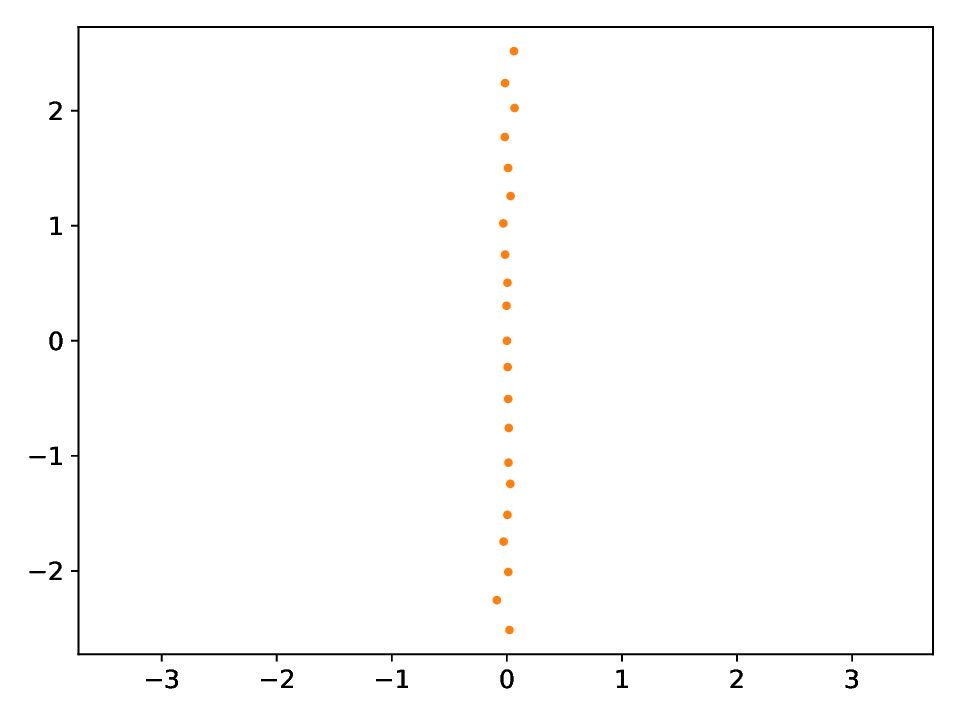}
    \caption{Data set (b) seen from the tangent space $\tangent_{\VectorC^{\mathrm{b}}} \Real^2$.}
    \label{fig:circle-rae-example-logs}
    \end{subfigure}
    \hfill
    \begin{subfigure}{0.48\linewidth}
    \centering
        \includegraphics[width=\linewidth]{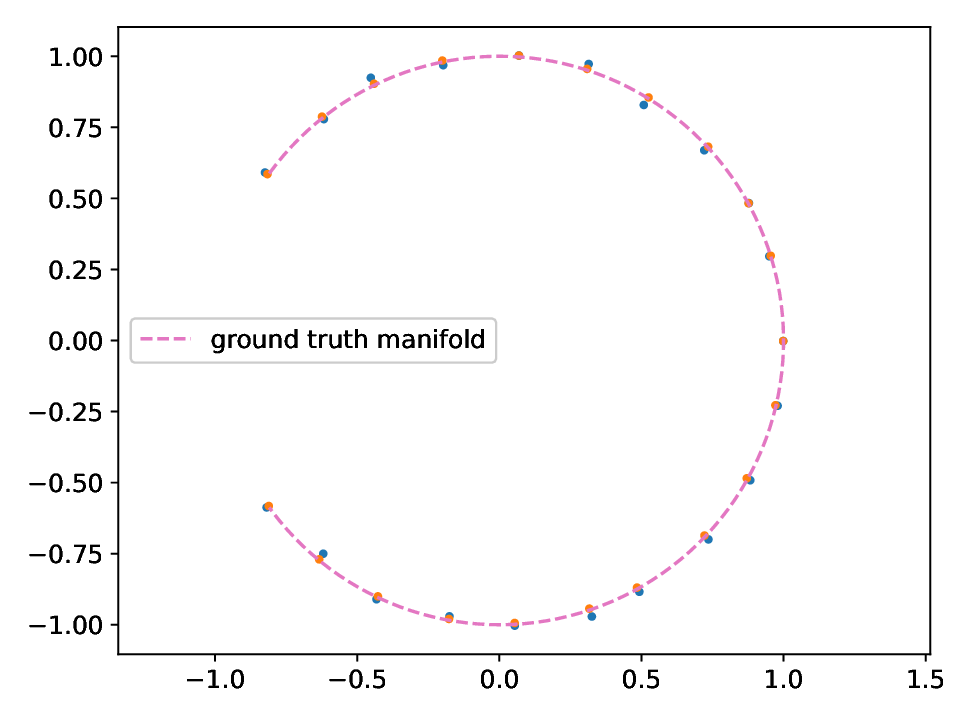}
    \caption{Data set (b) approximated by the RAE at $\VectorC^{\mathrm{b}}$.}
    \label{fig:circle-rae-example-rae}
    \end{subfigure}
    \caption{Low rank approximation of data set (b) on $(\Real^2, (\cdot,\cdot)^{\diffeo^{\mathrm{b}}})$ and the Riemannian autoencoder constructed from it indicate that the chosen pullback geometry is suitable for data analysis on this data set.}
    \label{fig:circle-rae-example}
\end{figure}

\paragraph{Riemannian autoencoder.}Finally, we compute the logarithmic mappings from $\VectorC^{\mathrm{b}}$ to all of the data, do a low rank approximation and construct a RAE, which is used to project the original data onto the learned manifold. Both results are shown in \cref{fig:circle-rae-example}. This time, the results are in line with expectations. That is, the data look 1-dimensional and linear on the tangent space (\cref{fig:circle-rae-example-logs}) and the RAE finds the correct data manifold (\cref{fig:circle-rae-example-rae}). Unsurprisingly, we retrieve the correct manifold without curvature correction, which normally is not expected to make a big difference anyways \cite{diepeveen2023curvature}.

\paragraph{Conclusions.} Positive curvature evidently seems to give rise to several instabilities for interpolation and barycentres that are hard to correct for. In contrast, it performs remarkably well for low rank approximation, where in general can be accounted for curvature. Overall, pulling back positively curved spaces can be problematic for downstream data analysis.

\subsection{Diffeomorphism effects}
\label{sec:numerics-pullback-manifolds-diffeo-effects}
Next, we consider learning pullback geometries under which we perform data analysis tasks from \cref{sec:basic-processing-pullback,sec:non-linear-compression-pullback} on data set (c) to get insight into (ii), i.e., the role of the diffeomorphism. In particular, under pullback geometry with learned diffeomorphisms -- obtained from training on \cref{eq:pwd-learning-problem} with and without subspace regularisation and local isometry regularisation -- we will consider stability of geodesics and barycentres and consider low rank approximation whose quality is tested through constructing Riemannian autoencoders\footnote{Since the pulled back Euclidean space has zero curvature we can skip the curvature correction step outlined in \cref{sec:cc-rae} and get a good minimiser from a basic truncated singular value decomposition in the tangent space.}.

For analysing data set (c), we will learn several diffeomorphism $\diffeo^{\mathrm{c}}_{\networkParams^1}, \diffeo^{\mathrm{c}}_{\networkParams^2}, \diffeo^{\mathrm{c}}_{\networkParams^3}, \diffeo^{\mathrm{c}}_{\networkParams^4}: \Real^2 \to \Real^{2}$ of the form \cref{eq:diffeo-class} into Euclidean space and pull back the standard Euclidean Riemannian structure. That is, for $\sumIndD = 1, 2, 3, 4$ and 
\begin{equation*}
    (\alpha_{\mathrm{sub}}^1,\alpha_{\mathrm{iso}}^1) = (10, 0.01), \quad (\alpha_{\mathrm{sub}}^2,\alpha_{\mathrm{iso}}^2) = (10, 0), \quad (\alpha_{\mathrm{sub}}^3,\alpha_{\mathrm{iso}}^3) = (0, 0.01), \quad (\alpha_{\mathrm{sub}}^4,\alpha_{\mathrm{iso}}^4) =  (0,0),
\end{equation*}
we consider $\diffeo^{\mathrm{c}}_{\networkParams^\sumIndD}: \Real^2 \to \Real^{2}$ with $\VectorC^{\mathrm{c}} := (0, 0) \in \Real^2$, $\mathbf{O}^{\mathrm{c}} \in \Real^{2\times 2}$ obtained from PCA at $\VectorC^{\mathrm{c}}$ as outlined in \cref{sec:learning-diffeo}, $\diffeoB_{\networkParams^\sumIndD}^{\mathrm{c}} :\Real^{2} \to \Real^{2}$ being an invertible residual network \cite{behrmann2019invertible} -- that is trained to minimise \cref{eq:pwd-learning-problem} for the above $\sumIndD$-dependent choice of $\alpha_{\mathrm{sub}}^\sumIndD,\alpha_{\mathrm{iso}}^\sumIndD$, but fixed approximate geodesic distances $\{\distance_{\sumIndA,\sumIndA'} \}_{\sumIndA,\sumIndA'=1}^\dataPointNum$ between all pairs of points $\Vector^\sumIndA$ and $\Vector^{\sumIndA'}$ that are constructed from completion of local Euclidean distances to neighbouring data points \cite{tenenbaum2000global} --, and $\diffeoC^{\mathrm{c}}: \Real\to \Real$ being identity, and pull back the $\ell^2$ inner product on $\Real$. Then, our data analysis is carried out on $(\Real^2, (\cdot,\cdot)^{\diffeo_{\networkParams^\sumIndD}^{\mathrm{c}}})$ for $\sumIndD = 1, 2, 3, 4$. 

\begin{figure}[h!]
    \centering
    \begin{subfigure}{0.48\linewidth}
    \centering
        \includegraphics[width=\linewidth]{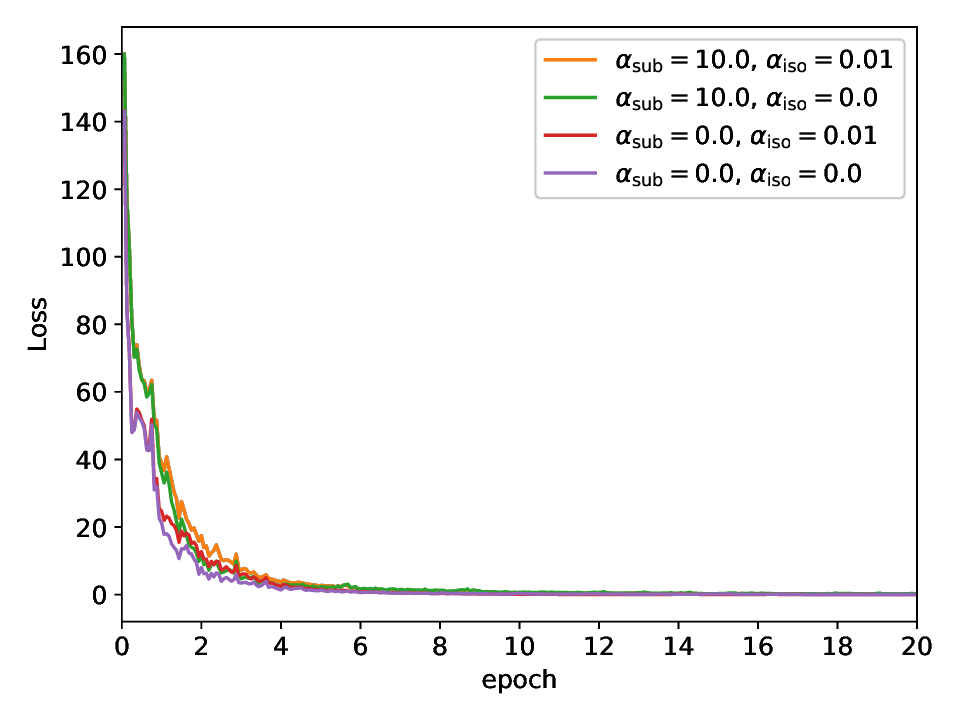}
    \caption{The progressions of the training losses \cref{eq:pwd-learning-problem}.}
    \label{fig:convergence-of-iresnets}
    \end{subfigure}
    \hfill
    \begin{subfigure}{0.48\linewidth}
    \centering
        \includegraphics[width=\linewidth]{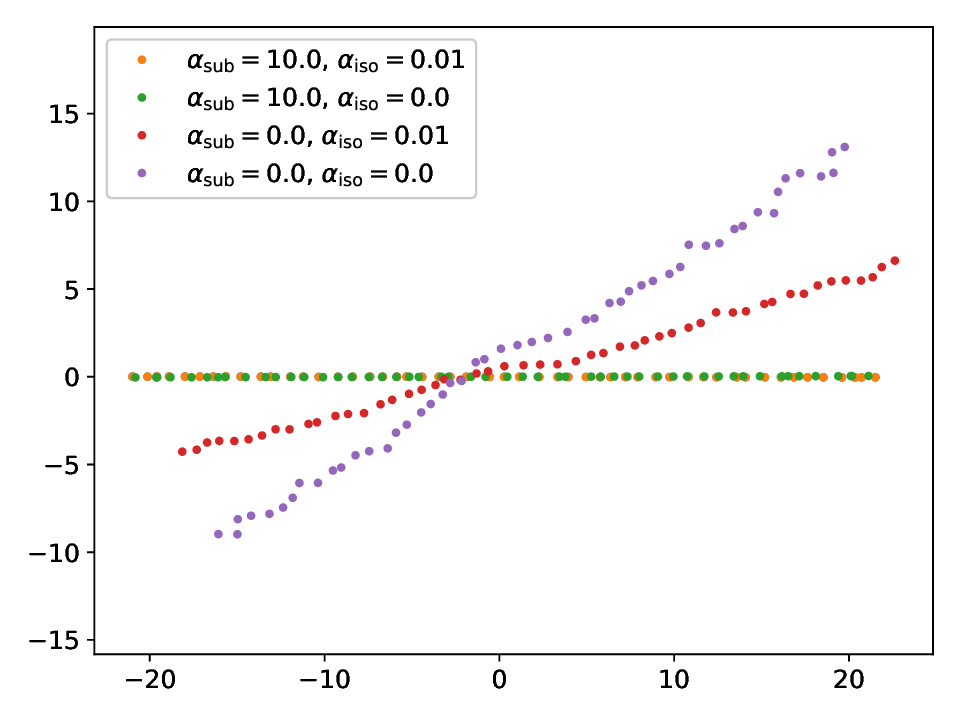}
    \caption{The learned diffeomorphisms mapping data set (c).}
    \label{fig:spiral-diffeo-range-on-data}
    \end{subfigure}
    \caption{The progressions of the losses of the learned several diffeomorphisms $\diffeo^{\mathrm{c}}_{\networkParams^1}, \diffeo^{\mathrm{c}}_{\networkParams^2}, \diffeo^{\mathrm{c}}_{\networkParams^3}, \diffeo^{\mathrm{c}}_{\networkParams^4}$ along with where the data sets are mapped to suggest that the networks have been properly trained after only 20 epochs.}
\end{figure}

For each $\sumIndD = 1, 2, 3, 4$, the invertible residual network $\diffeo^{\mathrm{c}}_{\networkParams^\sumIndD}$ consists of $100$ residual blocks, each consisting of $2$ hidden layers with width $10$ and ELU non-linearities \cite{clevert2015fast}. To enforce invertibility we do not restrict the bias terms, but normalize the weight matrices by their approximate spectral norm, which is computed through $10$ power iterations, and multiply them by a constant $c = 0.8$. The networks are trained for $20$ epochs with batch size $64$ using the ADAM optimiser \cite{kingma2014adam} with learning rate $\tau = 10^{-3}$ and exponential decay rates for the moment estimates $\beta_1 = 0.9$ and $\beta_2=0.99$. The progressions of the respective losses \cref{eq:pwd-learning-problem} in \cref{fig:convergence-of-iresnets} suggest that the networks already have converged after only 20 epochs. As a sanity check consider \cref{fig:spiral-diffeo-range-on-data} showing that the data set (c) has been (approximately) mapped into geodesic subspaces -- affine linear subspaces -- of $\Real^2$, which is a basic requirement for our diffeomorphisms we have encountered throughout \cref{sec:basic-processing-pullback,sec:non-linear-compression-pullback}. 

\begin{figure}[h!]
    \centering
    \begin{subfigure}{0.48\linewidth}
    \centering
        \includegraphics[width=\linewidth]{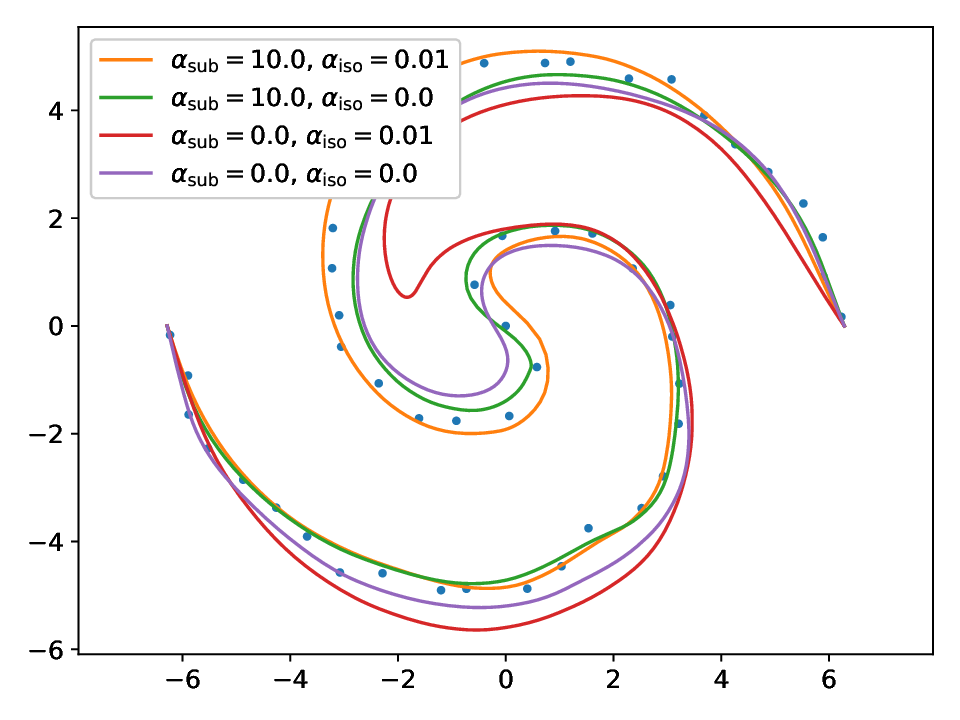}
    \caption{Geodesic interpolation.}
    \label{fig:spiral-interpolation-example-geodesic}
    \end{subfigure}
    \hfill
    \begin{subfigure}{0.48\linewidth}
    \centering
        \includegraphics[width=\linewidth]{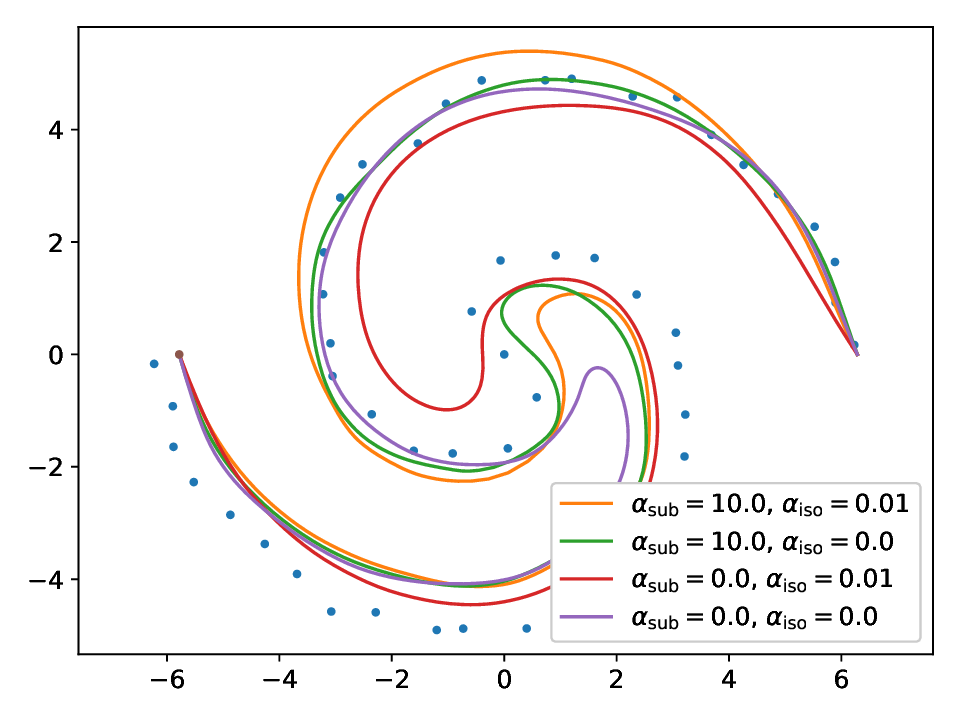}
    \caption{Perturbed geodesic interpolation.}
    \label{fig:spiral-interpolation-example-variation}
    \end{subfigure}
    \caption{Geodesic interpolations and perturbed geodesic interpolations of the end points of data set (c) on $(\Real^2, (\cdot,\cdot)^{\diffeo^{\mathrm{c}}_{\networkParams^\sumIndD}})$ for $\sumIndD=1,2,3,4$ indicate that pullback geometry trained with both subspace and isometry regularisation (orange) is suitable for data analysis on this data set and is stable with respect to small perturbations.}
    \label{fig:spiral-interpolation-example}
\end{figure}

\paragraph{Interpolation.}
Here too, we consider interpolation and its stability through interpolating between the end points of the data set (c) and vary one of the end points with an out of distribution point, see \cref{fig:spiral-interpolation-example}. \Cref{fig:spiral-interpolation-example-geodesic} shows that the learned pullback geometries each do what they should do to a different extent. In particular, visually the $\diffeo_{\networkParams^1}^{\mathrm{c}}$-geodesic (orange) -- trained with both subspace and isometry regularisation -- goes almost straight through the data set, the  $\diffeo_{\networkParams^2}^{\mathrm{c}}$-geodesic (green) and the $\diffeo_{\networkParams^4}^{\mathrm{c}}$-geodesic (purple) -- both trained without isometry regularisation -- are close runner ups, but the $\diffeo_{\networkParams^3}^{\mathrm{c}}$-geodesic (red) -- trained with only the isometry loss -- does not result in good interpolation. \Cref{tab:spiral-interpolation-example-errros} confirms our visual observations. That is, considering the \emph{geodesic error} (see \cref{eq:app-geodesic-error} in \cref{app:pullback-manifold-error-metrics-numerics}), we see that the $\diffeo_{\networkParams^1}^{\mathrm{c}}$-geodesic really is better than the other geodesics. Arguably, the lack of isometry regularisation in the training of $\diffeo_{\networkParams^2}^{\mathrm{c}}$ and $\diffeo_{\networkParams^4}^{\mathrm{c}}$ is a root cause, because now the  $\diffeo_{\networkParams^2}^{\mathrm{c}}$-geodesic and the $\diffeo_{\networkParams^4}^{\mathrm{c}}$-geodesic do not have constant speed in the $\ell^2$ sense, i.e., there are parts of the curve where the geodesic goes faster than in other parts, which causes errors. This can be backed up by considering the geodesic variation. Even though \cref{fig:spiral-interpolation-example-variation} suggests that especially $\diffeo_{\networkParams^2}^{\mathrm{c}}$ does not suffer from instabilities, \cref{tab:spiral-interpolation-example-errros} shows that it is actually performing worst in terms of \emph{variation error} (see \cref{eq:app-geodesic-variation-error} in \cref{app:pullback-manifold-error-metrics-numerics}) with $\diffeo_{\networkParams^4}^{\mathrm{c}}$ as a runner up, which can only be explained through errors between corresponding times along the geodesic and not the shape of the overall geodesic. This is unsurprising as \cref{thm:stability-geodesics-symmetric} suggests that for the networks trained without isometry regularisation we expect instability with respect to changing end points.

\begin{table}[h!]
\centering
\caption{The geodesic error \cref{eq:app-geodesic-error} and variation error \cref{eq:app-geodesic-variation-error} of the geodesics in \cref{fig:spiral-interpolation-example} with respect to data set (c) indicate that pullback geometry trained with both subspace and isometry regularisation is suitable for data analysis on this data set and is stable with respect to small perturbations.}
    \label{tab:spiral-interpolation-example-errros}
\begin{tabular}{c|cc|cc}
\toprule
 & $\alpha_{\mathrm{sub}}$ & $\alpha_{\mathrm{iso}}$ & geodesic error  & variation error \\
 \midrule
 $\diffeo_{\networkParams^1}^{\mathrm{c}}$ & $10$ & $0.01$  & $0.25 \pm 0.14$ & $0.62 \pm 0.45$  \\
 $\diffeo_{\networkParams^2}^{\mathrm{c}}$ & $10$ & $0$  & $0.38 \pm 0.24$ & $1.03 \pm 0.84$ \\
 $\diffeo_{\networkParams^3}^{\mathrm{c}}$ & $0$ & $0.01$ & $0.73 \pm 0.68$ & $0.68 \pm 0.49$ \\
 $\diffeo_{\networkParams^4}^{\mathrm{c}}$ & $0$ & $0$ & $0.37 \pm 0.21$ & $0.91 \pm 0.61$ \\
 \bottomrule
\end{tabular}
\end{table}

\begin{figure}[h!]
    \centering
    \begin{subfigure}{0.48\linewidth}
    \centering
        \includegraphics[width=\linewidth]{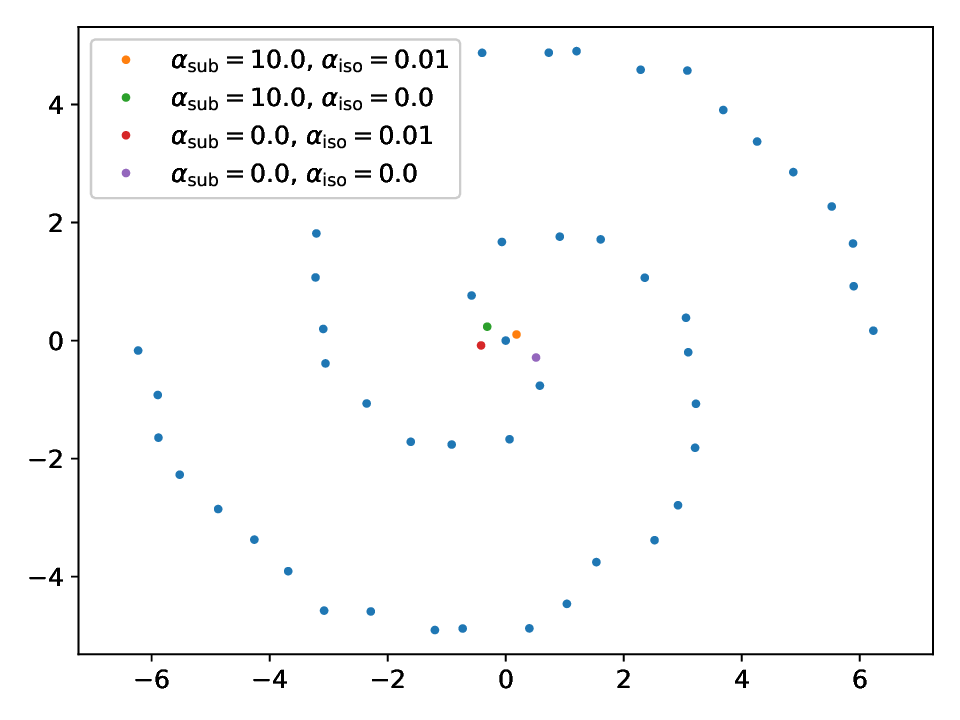}
    \caption{The data barycentre.}
    \label{fig:spiral-barycentre-example-standard}
    \end{subfigure}
    \hfill
    \begin{subfigure}{0.48\linewidth}
    \centering
        \includegraphics[width=\linewidth]{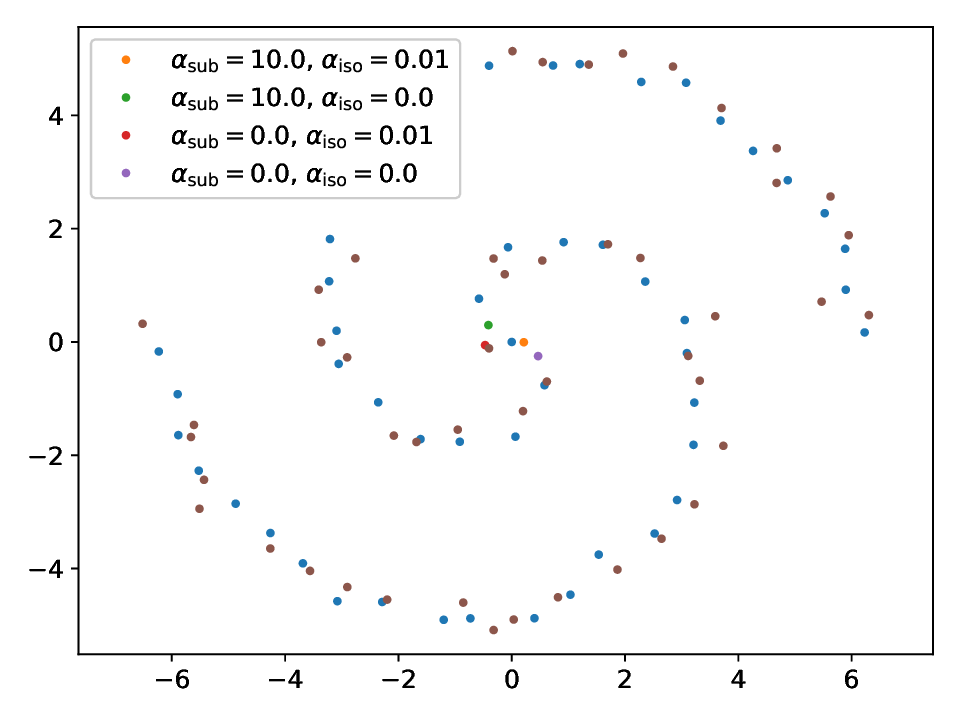}
    \caption{The perturbed data barycentre.}
    \label{fig:spiral-barycentre-example-variation}
    \end{subfigure}
    \caption{The data barycentres and perturbed data barycentres of data set (c) on $(\Real^2, (\cdot,\cdot)^{\diffeo^{\mathrm{c}}_{\networkParams^\sumIndD}})$ for $\sumIndD=1,2,3,4$ indicate that all chosen pullback geometries are suitable for data analysis on this data set and are stable with respect to small perturbations.}
    \label{fig:spiral-barycentre-example}
\end{figure}

\paragraph{Barycentre.}
Next, we consider the Riemannian barycentre and test its stability through small perturbations of the data, see \cref{fig:spiral-barycentre-example}. Here the results are somewhat more straightforward compared to evaluating geodesic interpolation. That is, similarly to the hyperbolic case in \cref{sec:numerics-pullback-manifolds-curvature-effects-hyperbolic}, \cref{fig:spiral-barycentre-example-standard} shows that the Riemannian barycentre is were it is expected to be, i.e., within the data set and close to $\VectorC^{\mathrm{c}}$ due to the symmetry of the data around that point. In addition, small variations of the data do not give rise to instabilities \cref{fig:spiral-barycentre-example-variation}, which is not surprising for the $\diffeo_{\networkParams^1}^{\mathrm{c}}$-barycentre and the $\diffeo_{\networkParams^3}^{\mathrm{c}}$-barycentre, but is surprising for the $\diffeo_{\networkParams^2}^{\mathrm{c}}$-barycentre and the $\diffeo_{\networkParams^4}^{\mathrm{c}}$-barycentre due to \cref{thm:stability-barycentre}.



\begin{figure}[h!]
    \centering
    \begin{subfigure}{0.48\linewidth}
    \centering
        \includegraphics[width=\linewidth]{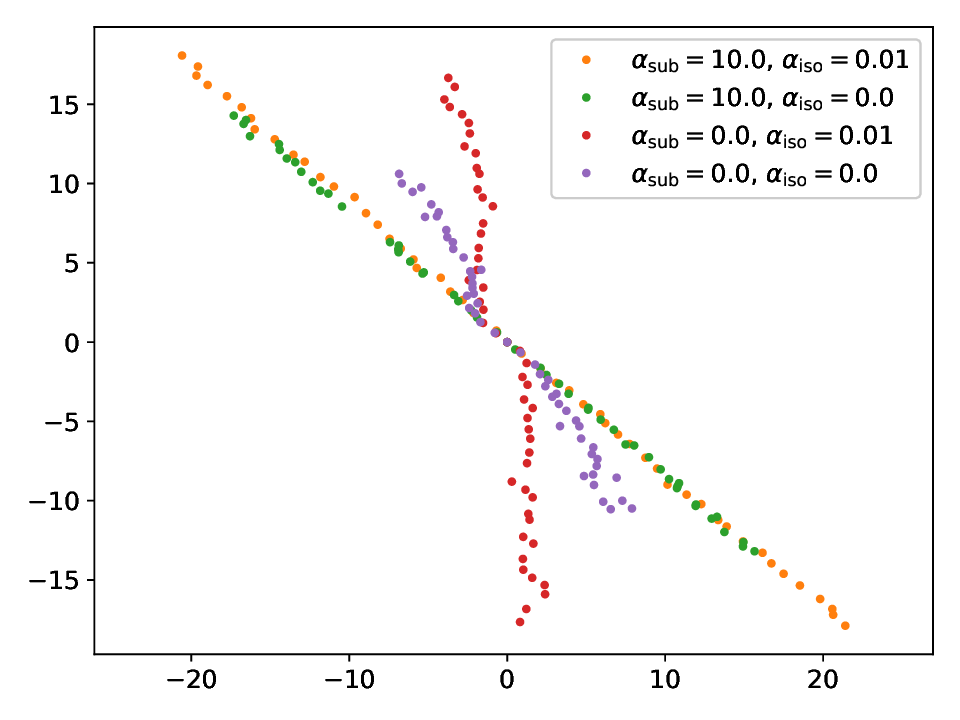}
    \caption{Data set (c) seen from the tangent space $\tangent_{\VectorC^{\mathrm{c}}} \Real^2$.}
    \label{fig:spiral-rae-example-logs}
    \end{subfigure}
    \hfill
    \begin{subfigure}{0.48\linewidth}
    \centering
        \includegraphics[width=\linewidth]{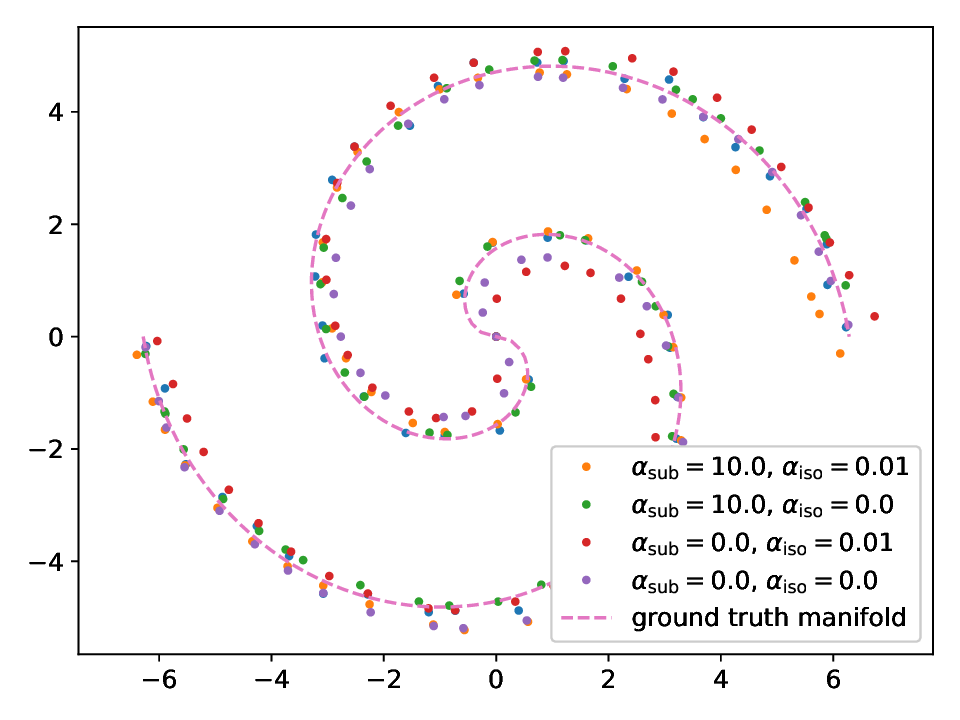}
    \caption{Data set (c) approximated by the RAEs at $\VectorC^{\mathrm{c}}$.}
    \label{fig:spiral-rae-example-rae}
    \end{subfigure}
    \caption{Low rank approximations of data set (c) on $(\Real^2, (\cdot,\cdot)^{\diffeo^{\mathrm{c}}_{\networkParams^\sumIndD}})$ for $\sumIndD=1,2,3,4$ and the Riemannian autoencoders constructed from them indicate that pullback geometry trained with both subspace and isometry regularisation (orange) is suitable for data analysis on this data set and is stable with respect to small perturbations.}
    \label{fig:spiral-rae-example}
\end{figure}

\paragraph{Riemannian autoencoder.}
Finally, we compute the logarithmic mappings from $\VectorC^{\mathrm{c}}$ to all of the data, do a low rank approximation and construct an RAE, which is used to project the original data onto the learned manifolds. Both results are shown in \cref{fig:spiral-rae-example}. Once again, the results are by and large in line with expectations. That is, for the $\diffeo_{\networkParams^1}^{\mathrm{c}}$-logarithmic mappings (orange) and the $\diffeo_{\networkParams^2}^{\mathrm{c}}$-logarithmic mappings (green) the data look 1-dimensional and linear on the tangent space, whereas under the $\diffeo_{\networkParams^3}^{\mathrm{c}}$-logarithmic mappings (red) and the $\diffeo_{\networkParams^4}^{\mathrm{c}}$-logarithmic mappings (purple) the data still looks somewhat non-linear (\cref{fig:spiral-rae-example-logs}). This discrepancy can once again be attributed to the network training, as we only really expect that we see 1-dimensionality for sure if we train with subspace loss. Upon closer inspection of the $\diffeo_{\networkParams^1}^{\mathrm{c}}$-and $\diffeo_{\networkParams^2}^{\mathrm{c}}$-logarithmic mappings the effect of isometry regularisation is also visible. That is, for the $\diffeo_{\networkParams^2}^{\mathrm{c}}$ case -- trained without isometry regularisation -- there are clear differences in the adjacent distances to subsequent tangent vectors (there are even several large gaps), whereas for the $\diffeo_{\networkParams^1}^{\mathrm{c}}$ case -- trained with full regularisation -- we do not observe this behaviour\footnote{Note that this is essentially the same as saying that $\diffeo_{\networkParams^2}^{\mathrm{c}}$-geodesics do not have constant speed through the data set in an $\ell^2$ sense, which we already observed before.}. Interestingly, all RAEs find the correct data manifold similarly well (\cref{fig:spiral-rae-example-rae}). Having said that, we can only expect that the latent space is isometric to Euclidean space in the case of $\diffeo_{\networkParams^1}^{\mathrm{c}}$, which can be important for interpretability (see \cref{rem:rae-interpretability}).
 



\paragraph{Conclusions}
Lack of isometry evidently seems to give rise to sub-optimal interpolation results and instabilities thereof. In addition, if the data manifold not being mapped into a geodesic subspace, we do not get that data looks as low-dimensional as it is. So overall, the above numerical experiments are in line with the theory that pulling back a learned diffeomorphism that is trained with both subspace and isometry regularisation can be beneficial for downstream data analysis.


\section{Conclusions}
\label{sec:conclusions-pullback-manifolds}

With this investigation, we hope to contribute to the development of a data-driven Riemannian geometry for interpretable, stable and efficient data analysis. In particular, within the broader scope of data processing, we believe that the proposed mathematical framework has important implications on how to construct data-driven Riemannian geometry and has important implications for handling data in general. 

In this work we aimed to address two main questions: \emph{how should we go about constructing diffeomorphisms into a symmetric Riemannian manifold?} and \emph{how should we use or modify successful algorithms on symmetric Riemannian manifolds for data analysis in the pullback geometry framework?}

\paragraph{Characterisation of diffeomorphisms for proper and stable data analysis.}
Regarding the first question, we have seen that diffeomorphisms should map the data manifold into low-dimensional geodesic subspaces of the pulled back Riemannian manifold for proper interpolation (\cref{thm:pull-back-mappings}) through and barycentres (\cref{thm:well-posedness-barycentre}) within the data set, but should also do it in such a way that it is a local isometry on the data in order to get $\ell^2$-stability (\cref{thm:stability-geodesics-symmetric,thm:stability-barycentre}). In addition, in \cref{thm:stability-geodesics-symmetric,thm:stability-barycentre} we have shown and quantified that the only other instabilities are due to curvature of the pulled back Riemannian manifold. This interplay between curvature effects and diffeomorphism effects on geodesics and barycentres has been observed in numerical experiments, in which we highlighted the potential detrimental effects of pulling back positive curvature Riemannian manifolds.

\paragraph{Characterisation of diffeomorphisms for efficient data analysis.}

Subsequently, regarding the second question we have seen that we can piggyback off of existing theory for data compression on symmetric Riemannian manifolds and use a recent efficient algorithm for low rank approximation to construct non-linear compression mappings, i.e., the \emph{Riemannian autoencoder} (RAE) and the \emph{curvature corrected Riemannian autoencoder} (CC-RAE). The developed ideas have shown promising results that were in line with expectations when tested in numerical experiments, i.e., the success of the methods are somewhat independent of the underlying curvature of the pulled back geometry -- and can in the worst case be corrected for --, but for interpretability of the latent space we want diffeomorphisms to map the data manifold isometrically into low-dimensional geodesic subspaces.

\paragraph{Construction of diffeomorphisms for proper, stable and efficient data analysis.}

Finally, having learned that diffeomorphisms need to map the data manifold into a geodesic subspace of the pulled back Riemannian manifold while preserving local isometry around the data manifold for proper, stable and efficient data analysis under pullback geometry, we have addressed how to construct such diffeomorphisms in a general setting using insights from several empirical successes. In particular, we have proposed a learning problem \cref{eq:pwd-learning-problem} to train invertible neural networks, which was then tested as a proof of concept through the above data analysis tasks. The numerical results were once again in line with the above expectations in the sense that low-dimensional embedding and local isometry are key for data analysis on pullback manifolds.

\section*{Acknowledgments}

We thank Andrea Bertozzi for helpful discussions and acknowledge support from the European Union Horizon 2020 research and innovation programme under the Marie Skodowska-Curie grant agreement No. 777826 NoMADS. 


\bibliographystyle{plain}
\bibliography{references} 

\clearpage
\appendix
\section{Proofs for the results from \cref{sec:prelim-pullback-geometry-theoretical-results}}
\label{app:proofs-prelim-pullback-geometry-theoretical-results}

\subsection{Proof of \cref{thm:pull-back-mappings}}
\paragraph{Auxiliary lemma.}

\begin{lemma}
\label{lem:connection}
Let $(\manifold, (\cdot,\cdot))$ be a $\dimInd$-dimensional Riemannian manifold with Levi-Civita connection denoted by $\nabla_{(\cdot)} (\cdot):\mathcal{X}(\manifold) \times \mathcal{X}(\manifold) \to \mathcal{X}(\manifold)$ and let $\diffeo:\Real^\dimInd \to \manifold$ be a smooth diffeomorphism. 
The Levi-Civita connection $\nabla_{(\cdot)}^\diffeo (\cdot): \mathcal{X}(\Real^\dimInd) \times \mathcal{X}(\Real^\dimInd) \to \mathcal{X}(\Real^\dimInd)$ on $(\Real^\dimInd, (\cdot,\cdot)^\diffeo)$ is given by
\begin{equation}
    \nabla_{\tangentVector}^\diffeo \tangentVectorB = \diffeo^{-1}_*[\nabla_{\diffeo_* [\tangentVector]} \diffeo_*[\tangentVectorB]] \quad \text{where $\tangentVector, \tangentVectorB\in \mathcal{X}(\Real^\dimInd)$}.
    \label{eq:levi-civita-pullback}
\end{equation}
\end{lemma}

\begin{proof}
First notice that the mapping $\diffeo^{-1}_*[\nabla_{\diffeo_* [\cdot]} \diffeo_*[\cdot]]: \mathcal{X}(\Real^\dimInd) \times \mathcal{X}(\Real^\dimInd) \to \mathcal{X}(\Real^\dimInd)$ defines a connection. Indeed, it is $C^{\infty}(\manifold)$-linear in the first variable 
\begin{multline}
    \diffeo^{-1}_*[\nabla_{\diffeo_* [f\tangentVector]} \diffeo_*[\tangentVectorB]] = \diffeo^{-1}_*[\nabla_{(f \circ \diffeo^{-1}) \diffeo_* [\tangentVector]} \diffeo_*[\tangentVectorB]] = \diffeo^{-1}_*[(f \circ \diffeo^{-1}) \nabla_{\diffeo_* [\tangentVector]} \diffeo_*[\tangentVectorB]] \\
    = f\diffeo^{-1}_*[\nabla_{\diffeo_* [\tangentVector]} \diffeo_*[\tangentVectorB]], \quad f\in C^{\infty}(\manifold)
\end{multline}
and satisfies the Leibniz rule in the second variable
\begin{multline}
    \diffeo^{-1}_*[\nabla_{\diffeo_* [\tangentVector]} \diffeo_*[f\tangentVectorB]] = \diffeo^{-1}_*[\nabla_{\diffeo_* [\tangentVector]} (f \circ \diffeo^{-1}) \diffeo_*[\tangentVectorB]] = \diffeo^{-1}_*[(f \circ \diffeo^{-1})\nabla_{\diffeo_* [\tangentVector]} \diffeo_*[\tangentVectorB] + (\diffeo_* [\tangentVector] (f \circ \diffeo^{-1})) \diffeo_*[\tangentVectorB]] \\
    =   f\diffeo^{-1}_*[\nabla_{\diffeo_* [\tangentVector]} \diffeo_*[\tangentVectorB]] + (\tangentVector f) \tangentVectorB, \quad f\in C^{\infty}(\manifold).
\end{multline}

Next, to proof the claim note that the connection in \cref{eq:levi-civita-pullback} is the Levi-Civita connection if and only if it satisfies the Koszul formula
\begin{equation}
    2 (\nabla_{\tangentVector}^\diffeo \tangentVectorB, \tangentVectorC)^\diffeo = \tangentVector(\tangentVectorB, \tangentVectorC)^\diffeo + \tangentVectorB(\tangentVector, \tangentVectorC)^\diffeo - \tangentVectorC(\tangentVector, \tangentVectorB)^\diffeo - ([\tangentVectorB, \tangentVector], \tangentVectorC)^\diffeo - ([\tangentVector, \tangentVectorC], \tangentVectorB)^\diffeo - ([\tangentVectorB, \tangentVectorC], \tangentVector)^\diffeo
    \label{eq:koszul-pullback}
\end{equation}
for all $\tangentVector, \tangentVectorB, \tangentVectorC\in \mathcal{X}(\Real^\dimInd)$.

It remains to show that \cref{eq:levi-civita-pullback} satisfies \cref{eq:koszul-pullback}. We will consider the first and the fourth term on the right-hand-side of \cref{eq:koszul-pullback}, i.e., $\tangentVector(\tangentVectorB, \tangentVectorC)^\diffeo$ and $([\tangentVectorB, \tangentVector], \tangentVectorC)^\diffeo$. The other terms will follow analogously.

\begin{multline}
    \tangentVector(\tangentVectorB, \tangentVectorC)^\diffeo = \tangentVector(\tangentVectorB, \tangentVectorC)^\diffeo_{(\cdot)} = \tangentVector(D_{(\cdot)} \diffeo [\tangentVectorB], D_{(\cdot)} \diffeo[\tangentVectorC])_{\diffeo(\cdot)}  = \tangentVector(D_{\diffeo^{-1}(\diffeo(\cdot))} \diffeo [\tangentVectorB], D_{\diffeo^{-1}(\diffeo(\cdot))} \diffeo[\tangentVectorC])_{\diffeo(\cdot)} \\
    = D_{(\cdot)} \diffeo [\tangentVector](D_{\diffeo^{-1}(\cdot)} \diffeo [\tangentVectorB], D_{\diffeo^{-1}(\cdot)} \diffeo[\tangentVectorC])_{(\cdot)} = \diffeo_*[\tangentVector] (\diffeo_* [\tangentVectorB], \diffeo_* [\tangentVectorC])  
    \label{eq:koszul-pullback-1}
\end{multline}

\begin{multline}
    ([\tangentVectorB, \tangentVector], \tangentVectorC)^\diffeo = ([\tangentVectorB, \tangentVector], \tangentVectorC)^\diffeo_{(\cdot)} = (D_{(\cdot)} \diffeo [[\tangentVectorB, \tangentVector]],D_{(\cdot)} \diffeo [ \tangentVectorC])_{\diffeo(\cdot)} = ([D_{(\cdot)} \diffeo [\tangentVectorB], D_{(\cdot)} \diffeo [\tangentVector]],D_{(\cdot)} \diffeo [ \tangentVectorC])_{\diffeo(\cdot)}\\
    = ([D_{\diffeo^{-1}(\diffeo(\cdot))} \diffeo [\tangentVectorB], D_{\diffeo^{-1}(\diffeo(\cdot))} \diffeo [\tangentVector]],D_{\diffeo^{-1}(\diffeo(\cdot))} \diffeo [ \tangentVectorC])_{\diffeo(\cdot)} = ([D_{\diffeo^{-1}(\cdot)} \diffeo [\tangentVectorB], D_{\diffeo^{-1}(\cdot)} \diffeo [\tangentVector]],D_{\diffeo^{-1}(\cdot)} \diffeo [ \tangentVectorC])_{(\cdot)}\\
    = ([\diffeo_* [\tangentVectorB], \diffeo_*[\tangentVector]],\diffeo_*[ \tangentVectorC])
    \label{eq:koszul-pullback-2}
\end{multline}
where we used that $\diffeo_*[\tangentVectorB, \tangentVector] = [\diffeo_*[\tangentVectorB], \diffeo_*[\tangentVector]]$.

Then, substituting \cref{eq:koszul-pullback-1,eq:koszul-pullback-2} into \cref{eq:koszul-pullback} gives the desired result
\begin{multline}
   2 (\nabla_{\tangentVector}^\diffeo \tangentVectorB, \tangentVectorC)^\diffeo = \diffeo_*[\tangentVector] (\diffeo_* [\tangentVectorB], \diffeo_* [\tangentVectorC])   + \diffeo_*[\tangentVectorB] (\diffeo_* [\tangentVector], \diffeo_* [\tangentVectorC])   - \diffeo_*[\tangentVectorC] (\diffeo_* [\tangentVector], \diffeo_* [\tangentVectorB])   \\
   - ([\diffeo_* [\tangentVectorB], \diffeo_*[\tangentVector]],\diffeo_*[ \tangentVectorC]) - ([\diffeo_* [\tangentVector], \diffeo_*[\tangentVectorC]],\diffeo_*[ \tangentVectorB]) - ([\diffeo_* [\tangentVectorB], \diffeo_*[\tangentVectorC]],\diffeo_*[ \tangentVector])\\
   = 2(\nabla_{\diffeo_*[ \tangentVector]} \diffeo_* [\tangentVectorB], \diffeo_*[\tangentVectorC]) = 2(\diffeo_* [\diffeo^{-1}_* [\nabla_{\diffeo_*[ \tangentVector]} \diffeo_* [\tangentVectorB]]], \diffeo_*[\tangentVectorC]) = 2 (\diffeo^{-1}_* [\nabla_{\diffeo_*[ \tangentVector]} \diffeo_* [\tangentVectorB]], \tangentVectorC)^\diffeo.
\end{multline}

\end{proof}

\paragraph{Proof of the proposition.}

\begin{proof}[Proof of \cref{thm:pull-back-mappings}]
    We will show the identities in the same order as in the statement.
    \begin{enumerate}[label=(\roman*)]
        \item The set of piece-wise smooth curves mapping into a smooth manifold $\manifoldB$ is denoted by $PC^\infty([0,1], \manifoldB)$. Next, define the subset $PC_{\mPoint,\mPointB}^\infty([0,1], \manifoldB)$ as 
        \begin{equation}
            PC_{\mPoint,\mPointB}^\infty([0,1], \manifoldB) := \{c \in PC^\infty([0,1], \manifoldB)\; \mid \; c(0) =  \mPoint, c(1) =  \mPointB \}. 
        \end{equation}
        Then, for proving \cref{eq:thm-geodesic-remetrized}, note that geodesics minimise by definition the variational problem 
        \begin{multline}
            \inf_{c^\diffeo_{\Vector, \VectorB} \in PC_{\Vector,\VectorB}^\infty([0,1], \Real^\dimInd)} \Bigl\{ \int_{0}^1 \|\dot{c}^\diffeo_{\Vector, \VectorB}\|^\diffeo_{c^\diffeo_{\Vector, \VectorB}(t)} \; \mathrm{d}t \Bigr\} 
            = \inf_{c^\diffeo_{\Vector, \VectorB} \in PC_{\Vector,\VectorB}^\infty([0,1], \Real^\dimInd)} \Bigl\{  \int_{0}^1 \|\diffeo_*[\dot{c}^\diffeo_{\Vector, \VectorB}]\|_{\diffeo(c^\diffeo_{\Vector, \VectorB}(t))} \; \mathrm{d}t  \Bigr\} \\
            = \inf_{\substack{c_{\diffeo(\Vector), \diffeo(\VectorB)} \in PC_{\diffeo(\Vector), \diffeo(\VectorB)}^\infty([0,1], \manifold)\\ s.t. \; c_{\diffeo(\Vector), \diffeo(\VectorB)} =  \diffeo \circ c^\diffeo_{\Vector, \VectorB}, \\ c^\diffeo_{\Vector, \VectorB} \in PC_{\Vector,\VectorB}^\infty([0,1], \Real^\dimInd)}} \Bigl\{  \int_{0}^1 \|\dot{c}_{\diffeo(\Vector), \diffeo(\VectorB)}\|_{c_{\diffeo(\Vector), \diffeo(\VectorB)}(t)} \; \mathrm{d}t  \Bigr\} \\
            \geq \inf_{c_{\diffeo(\Vector), \diffeo(\VectorB)} \in PC_{\diffeo(\Vector), \diffeo(\VectorB)}^\infty([0,1], \manifold)} \int_{0}^1 \|\dot{c}_{\diffeo(\Vector), \diffeo(\VectorB)}\|_{c_{\diffeo(\Vector), \diffeo(\VectorB)}(t)} \; \mathrm{d}t
            \label{eq:geodesic-lower-bound}
        \end{multline}
        The lower bound is minimised by the geodesic $\geodesic_{\diffeo(\Vector), \diffeo(\VectorB)}$. Then, through geodesic convexity of $\diffeo(\Real^\dimInd)$ the curve $\geodesic^\diffeo_{\Vector, \VectorB} := \diffeo^{-1}\circ \geodesic_{\diffeo(\Vector), \diffeo(\VectorB)}$ is well-defined. Moreover, it is easy to check that $\geodesic^\diffeo_{\Vector, \VectorB}$ attains the lower bound in \cref{eq:geodesic-lower-bound}. In other words, $\geodesic^\diffeo_{\Vector, \VectorB}$ is a length-minimising geodesic on the Riemannian manifold $(\Real^\dimInd, (\cdot,\cdot)^\diffeo)$, which proofs the claim.
        \item For proving \cref{eq:thm-log-remetrized}, consider
        \begin{equation}
            \log^\diffeo_{\Vector} \VectorB = \dot{\geodesic}^\diffeo_{\Vector, \VectorB}(0) \overset{\cref{eq:thm-geodesic-remetrized}}{=}  \diffeo^{-1}_{*}[\dot{\geodesic}_{\diffeo(\Vector),\diffeo(\VectorB)}(0)] = \diffeo^{-1}_{*}[\log_{\diffeo(\Vector)} \diffeo(\VectorB)],
        \end{equation}
        which proves the claim.
        \item For proving \cref{eq:thm-exp-remetrized}, it is sufficient to show that the proposed mapping is the inverse of the logarithmic mapping, i.e., 
        \begin{equation}
            \diffeo^{-1}(\exp_{\diffeo(\Vector)}(\diffeo_*[\log^\diffeo_{\Vector} \VectorB])) = \VectorB \quad \text{and} \quad \log^\diffeo_{\Vector} (\diffeo^{-1}(\exp_{\diffeo(\Vector)}(\diffeo_*[\tangentVector_\Vector])) = \tangentVector_\Vector.
        \end{equation}
        First note that for any $\Vector, \VectorB\in \Real^\dimInd$ and any $\tangentVector_\Vector \in \mathcal{G}_{\Vector}$ both mappings $\diffeo^{-1}(\exp_{\diffeo(\Vector)}(\diffeo_*[\log^\diffeo_{\Vector} \VectorB]))$ and $\log^\diffeo_{\Vector} (\diffeo^{-1}(\exp_{\diffeo(\Vector)}(\diffeo_*[\tangentVector_\Vector]))$ are well-defined. Next, direct evaluation gives
        \begin{multline}
            \diffeo^{-1}(\exp_{\diffeo(\Vector)}(\diffeo_*[\log^\diffeo_{\Vector} \VectorB])) \overset{\cref{eq:thm-log-remetrized}}{=} \diffeo^{-1}(\exp_{\diffeo(\Vector)}(\diffeo_*[\diffeo^{-1}_{*}[\log_{\diffeo(\Vector)} \diffeo(\VectorB)]])) \\
            = \diffeo^{-1}(\exp_{\diffeo(\Vector)}(\log_{\diffeo(\Vector)} \diffeo(\VectorB)])) = \diffeo^{-1}(\diffeo(\VectorB)) = \VectorB,
        \end{multline}
        and 
        \begin{multline}
            \log^\diffeo_{\Vector} (\diffeo^{-1}(\exp_{\diffeo(\Vector)}(\diffeo_*[\tangentVector_\Vector]))) \overset{\cref{eq:thm-log-remetrized}}{=} \diffeo^{-1}_{*}[\log_{\diffeo(\Vector)} (\diffeo(\diffeo^{-1}(\exp_{\diffeo(\Vector)}(\diffeo_*[\tangentVector_\Vector]))))] \\
            = \diffeo^{-1}_{*}[\log_{\diffeo(\Vector)} (\exp_{\diffeo(\Vector)}(\diffeo_*[\tangentVector_\Vector]))] = = \diffeo^{-1}_{*}[\diffeo_*[\tangentVector_\Vector]] = \tangentVector_\Vector,
        \end{multline}
        which proves the claim.
        \item For proving \cref{eq:thm-distance-remetrized}, consider
        \begin{multline}
            \distance^\diffeo_{\Real^\dimInd}(\Vector, \VectorB) = \inf_{c^\diffeo_{\Vector, \VectorB} \in PC_{\Vector,\VectorB}^\infty([0,1], \Real^\dimInd)} \Bigl\{ \int_{0}^1 \|\dot{c}^\diffeo_{\Vector, \VectorB}\|^\diffeo_{c^\diffeo_{\Vector, \VectorB}(t)} \; \mathrm{d}t \Bigr\} \\
            \overset{(i)}{=} \inf_{c_{\diffeo(\Vector), \diffeo(\VectorB)} \in PC_{\diffeo(\Vector), \diffeo(\VectorB)}^\infty([0,1], \manifold)} \int_{0}^1 \|\dot{c}_{\diffeo(\Vector), \diffeo(\VectorB)}\|_{c_{\diffeo(\Vector), \diffeo(\VectorB)}(t)} \; \mathrm{d}t  = \distance_{\manifold}(\diffeo(\Vector), \diffeo(\VectorB)),
        \end{multline}
        which proves the claim.
        \item For proving \cref{eq:thm-parallel-transport-remetrized}, we must show that
        \begin{equation}
            \nabla_{\dot{\geodesic}^\diffeo_{\Vector, \VectorB}(t)}^\diffeo \diffeo^{-1}_* [\mathcal{P}_{\diffeo(\geodesic^\diffeo_{\Vector, \VectorB}(t)) \leftarrow \diffeo(\Vector)}  \diffeo_*[\tangentVector_\Vector]] = 0, \quad \text{for $t\in [0,1]$.}
        \end{equation}
        Direct evaluation through \cref{lem:connection} gives
        \begin{multline}
            \nabla_{\dot{\geodesic}^\diffeo_{\Vector, \VectorB}(t)}^\diffeo \diffeo^{-1}_* [\mathcal{P}_{\diffeo(\geodesic^\diffeo_{\Vector, \VectorB}(t)) \leftarrow \diffeo(\Vector)}  \diffeo_*[\tangentVector_\Vector]] \overset{\cref{eq:levi-civita-pullback}}{=} \diffeo^{-1}_*[\nabla_{\diffeo_* [\dot{\geodesic}^\diffeo_{\Vector, \VectorB}(t)]} \diffeo_*[\diffeo^{-1}_* [\mathcal{P}_{\diffeo(\geodesic^\diffeo_{\Vector, \VectorB}(t)) \leftarrow \diffeo(\Vector)}  \diffeo_*[\tangentVector_\Vector]] ]]\\
            = \diffeo^{-1}_*[\nabla_{\diffeo_* [\dot{\geodesic}^\diffeo_{\Vector, \VectorB}(t)]} \mathcal{P}_{\diffeo(\geodesic^\diffeo_{\Vector, \VectorB}(t)) \leftarrow \diffeo(\Vector)}  \diffeo_*[\tangentVector_\Vector]] \overset{\cref{eq:thm-geodesic-remetrized}}{=} \diffeo^{-1}_*[\nabla_{\dot{\geodesic}_{\diffeo(\Vector), \diffeo(\VectorB)}(t)} \mathcal{P}_{\geodesic_{\diffeo(\Vector), \diffeo(\VectorB)}(t) \leftarrow \diffeo(\Vector)}  \diffeo_*[\tangentVector_\Vector]]\\
            = \diffeo^{-1}_*[ 0 ] = 0,
        \end{multline}
        which proves the claim.
    \end{enumerate}
\end{proof}

\subsection{Proof of \cref{thm:completeness}}

\begin{proof}[Proof of \cref{thm:completeness}]
    Choose any point $\Vector\in \Real^\dimInd$. By the Hopf-Rinow Theorem it is sufficient to show that the exponential mapping $\exp^\diffeo_\Vector$ is defined on all of $\tangent_\Vector \Real^\dimInd$. 
    
    Using identity (iii) from \cref{thm:pull-back-mappings} we obtain the following equivalency:

    \begin{equation}
    \begin{gathered}
        \exp^\diffeo_{\Vector} \text{ is defined on all of $\tangentVector_\Vector \in \tangent_\Vector \Real^\dimInd$}, \\
        \quad \Leftrightarrow \quad \\
        \exp_{\diffeo(\Vector)} \text{ is defined on all of $\tangent_{\diffeo(\Vector)} \manifold$ and } \exp_{\diffeo(\Vector)} (\tangent_{\diffeo(\Vector)} \manifold) \subset \diffeo(\Real^\dimInd) .
    \end{gathered}
    \end{equation}
    The lower statement holds by assumption. Indeed, using the completeness of $\manifold$ and the Hopf-Rinow Theorem once more we have that $\exp_{\diffeo(\Vector)}$ is defined on all of $\tangent_{\diffeo(\Vector)} \manifold$ and that $\exp_{\diffeo(\Vector)} (\tangent_{\diffeo(\Vector)} \manifold) = \manifold$, which is equal to $\diffeo(\Real^\dimInd)$ due to the fact that $\diffeo$ is a global diffeomorphism.
    
\end{proof}

\subsection{Proof of \cref{thm:local-symmetry}}

\begin{proof}[Proof of \cref{thm:local-symmetry}]
    For proving the claim, we have to construct a geodesic reflection and open neighbourhoods in $\Real^\dimInd$, and show that this geodesic reflection is a local isometry.

    Choose $\Vector\in \Real^\dimInd$, let $\reflection_{\diffeo(\Vector)}: \mathcal{U}(\diffeo(\Vector)) \rightarrow \manifold$ be a geodesic reflection on a neighbourhood $\mathcal{U}(\diffeo(\Vector))$ and consider an open subset $\mathcal{V}(\diffeo(\Vector)) \subset\{\mPoint \in \manifold \mid \mPoint \in \mathcal{U}(\diffeo(\Vector)) \cap \diffeo (\Real^\dimInd), \;  \reflection_{\diffeo(\Vector)}(\mPoint) \in  \diffeo (\Real^\dimInd)\}$. Next, define the mapping $\reflection_{\Vector}^\diffeo: \diffeo^{-1}(\mathcal{V}(\diffeo(\Vector))) \rightarrow \Real^\dimInd$ as 
    \begin{equation}
        \reflection_{\Vector}^\diffeo := \diffeo^{-1} \circ \reflection_{\diffeo(\Vector)} \circ \; \diffeo.
        \label{eq:reflection-pullback}
    \end{equation}
    We will show that $\reflection_{\Vector}^\diffeo$ is a geodesic reflection on $\diffeo^{-1}(\mathcal{V}(\diffeo(\Vector)))$, which is a neighbourhood of $\Vector$. For that, we check the properties in \cref{eq:geodesic-reflection-def-properties}, both of which follow from direct evaluation:
    \begin{equation}
        \reflection_{\Vector}^\diffeo (\Vector) = \diffeo^{-1} (\reflection_{\diffeo(\Vector)}(\diffeo(\Vector))) \overset{\text{reflection}}{=} \diffeo^{-1} (\diffeo(\Vector)) = \Vector
    \end{equation}
    and
    \begin{multline}
        D_{\Vector}\reflection_{\Vector}^\diffeo  \overset{\text{chain rule}}{=} D_{\diffeo(\Vector)}\diffeo^{-1}  \circ D_{\diffeo(\Vector)} \reflection_{\diffeo(\Vector)}  \circ \; D_{\Vector} \diffeo 
        \overset{\text{reflection}}{=} D_{\diffeo(\Vector)}\diffeo^{-1} \circ (-\operatorname{id}_{\diffeo(\Vector)}) \circ \; D_{\Vector} \diffeo  \\
        = - D_{\diffeo(\Vector)}\diffeo^{-1}   \circ \; D_{\Vector} \diffeo   = - \operatorname{id}_{\Vector},
    \end{multline}
    which proves the claim that $\reflection_{\Vector}^\diffeo$ is a geodesic reflection on $\diffeo^{-1}(\mathcal{V}(\diffeo(\Vector)))$.

    It remains to check isometry of $\reflection_{\Vector}^\diffeo$ on $\diffeo^{-1}(\mathcal{V}(\diffeo(\Vector)))$. Choose $\VectorB, \VectorC \in \diffeo^{-1}(\mathcal{V}(\diffeo(\Vector)))$, then 
    \begin{multline}
        \distance^\diffeo_{\Real^\dimInd}(\reflection_{\Vector}^\diffeo(\VectorB), \reflection_{\Vector}^\diffeo(\VectorC)) \overset{\substack{\text{\cref{thm:pull-back-mappings} (iv)}}}{=} \distance_{\manifold}(\diffeo (\reflection_{\Vector}^\diffeo(\VectorB)), \diffeo(\reflection_{\Vector}^\diffeo(\VectorC))) \\
        \overset{\cref{eq:reflection-pullback}}{=} \distance_{\manifold}(\diffeo (\diffeo^{-1} (\reflection_{\diffeo(\Vector)}( \diffeo(\VectorB)))), \diffeo(\diffeo^{-1} (\reflection_{\diffeo(\Vector)}( \diffeo(\VectorC))))) = \distance_{\manifold}(\reflection_{\diffeo(\Vector)}( \diffeo(\VectorB)), \reflection_{\diffeo(\Vector)}( \diffeo(\VectorC))) \\
        \overset{\text{symmetry \cref{eq:symmetry-properties}}}{=} \distance_{\manifold}(\diffeo(\VectorB), \diffeo(\VectorC)) \overset{\substack{\text{\cref{thm:pull-back-mappings} (iv)}}}{=} \distance^\diffeo_{\Real^\dimInd}(\VectorB, \VectorC)
    \end{multline}
    
\end{proof}

\subsection{Proof of \cref{thm:hadamard}}

\paragraph{Auxiliary lemma.}

\begin{lemma}
\label{lem:curvature-operator}
Let $(\manifold, (\cdot,\cdot))$ be a $\dimInd$-dimensional Riemannian manifold with curvature tensor denoted by $\curvature:\mathcal{X}(\manifold) \times \mathcal{X}(\manifold) \times \mathcal{X}(\manifold) \to \mathcal{X}(\manifold)$ and let $\diffeo:\Real^\dimInd \to \manifold$ be a smooth diffeomorphism. The curvature tensor $\curvature^\diffeo: \mathcal{X}(\Real^\dimInd) \times \mathcal{X}(\Real^\dimInd) \times \mathcal{X}(\Real^\dimInd) \to \mathcal{X}(\Real^\dimInd)$ on $(\Real^\dimInd, (\cdot,\cdot)^\diffeo)$ is given by 
    \begin{equation}
    \curvature^\diffeo(\tangentVector, \tangentVectorB) \tangentVectorC = \diffeo^{-1}_*[\curvature(\diffeo_* [\tangentVector], \diffeo_* [\tangentVectorB])\diffeo_* [\tangentVectorC]].
    \label{eq:pull-back-curvature}
\end{equation}
\end{lemma}
\begin{proof}
    Direct evaluation through \cref{lem:connection} gives
    \begin{multline}
        \curvature^\diffeo(\tangentVector, \tangentVectorB) \tangentVectorC = \nabla_{\tangentVector}^\diffeo \nabla_{\tangentVectorB}^\diffeo \tangentVectorC - \nabla_{\tangentVectorB}^\diffeo \nabla_{\tangentVector}^\diffeo \tangentVectorC - \nabla_{[\tangentVector, \tangentVectorB]}^\diffeo \tangentVectorC \\
        \overset{\cref{eq:levi-civita-pullback}}{=} \diffeo^{-1}_*[\nabla_{\diffeo_* [\tangentVector]} \diffeo_*[ \diffeo^{-1}_*[\nabla_{\diffeo_* [\tangentVectorB]} \diffeo_*[ \tangentVectorC]]]] - \diffeo^{-1}_*[\nabla_{\diffeo_* [\tangentVectorB]} \diffeo_*[ \diffeo^{-1}_*[\nabla_{\diffeo_* [\tangentVector]} \diffeo_*[ \tangentVectorC]]]] - \diffeo^{-1}_*[\nabla_{\diffeo_*[[\tangentVector, \tangentVectorB]]} \diffeo_* [\tangentVectorC]]\\
        = \diffeo^{-1}_*[\nabla_{\diffeo_* [\tangentVector]} \nabla_{\diffeo_* [\tangentVectorB]} \diffeo_*[ \tangentVectorC]] - \diffeo^{-1}_*[\nabla_{\diffeo_* [\tangentVectorB]} \nabla_{\diffeo_* [\tangentVector]} \diffeo_*[ \tangentVectorC]] - \diffeo^{-1}_*[\nabla_{[\diffeo_*[\tangentVector], \diffeo_*[\tangentVectorB]]} \diffeo_* [\tangentVectorC]]\\
        = \diffeo^{-1}_*[\curvature(\diffeo_* [\tangentVector], \diffeo_* [\tangentVectorB])\diffeo_* [\tangentVectorC]]
    \end{multline}
\end{proof}

\paragraph{Proof of the proposition.}

\begin{proof}[Proof of \cref{thm:hadamard}]
    To proof the claim we need to show that $(\Real^\dimInd,  (\cdot,\cdot)^\diffeo)$ is complete, simply connected and has non-positive sectional curvature. As $\Real^\dimInd$ is simply connected it remains to show (i) completeness and (ii) non-positive sectional curvature.

    (i) Note that Hadamard manifolds are complete by definition and have unique geodesics, i.e., these spaces are geodesically convex. So completeness follows directly from \cref{thm:completeness}.

    (ii) Let $\Vector\in \Real^\dimInd$ be any point and choose a two-dimensional subspace $\sigma_{\Vector} \subset \tangent_\Vector \Real^\dimInd$ spanned by a basis $\{\tangentVector_\mPoint, \tangentVectorB_\mPoint\}$. Direct evaluation and \cref{lem:curvature-operator} gives
    \begin{multline}
        \sectionalcurvature_\Vector^\diffeo(\sigma_\Vector) =
        \frac{( \curvature_\Vector^\diffeo(\tangentVector_\Vector, \tangentVectorB_\Vector) \tangentVectorB_\Vector, \tangentVector_\Vector)_\Vector^\diffeo}{(\|\tangentVector_\Vector\|_\Vector^\diffeo)^2 (\|\tangentVectorB_\Vector\|_\Vector^\diffeo)^2 - (( \tangentVector_\Vector, \tangentVectorB_\Vector)_\Vector^\diffeo)^2} 
        = \frac{( \diffeo_*[\curvature_\Vector^\diffeo(\tangentVector_\Vector, \tangentVectorB_\Vector) \tangentVectorB_\Vector], \diffeo_*[\tangentVector_\Vector])_{\diffeo(\Vector)}}{\|\diffeo_*[\tangentVector_\Vector]\|_{\diffeo(\Vector)}^2 \|\diffeo_*[\tangentVectorB_\Vector]\|_{\diffeo(\Vector)}^2 - ( \diffeo_*[\tangentVector_\Vector], \diffeo_*[\tangentVectorB_\Vector])_{\diffeo(\Vector)}^2} \\
        \overset{\cref{eq:pull-back-curvature}}{=} \frac{( \diffeo_*[\diffeo^{-1}_*[\curvature_{\diffeo(_\Vector)}(\diffeo_* [\tangentVector_\Vector], \diffeo_* [\tangentVectorB_\Vector])\diffeo_* [\tangentVectorB_\Vector]]], \diffeo_*[\tangentVector_\Vector])_{\diffeo(\Vector)}}{\|\diffeo_*[\tangentVector_\Vector]\|_{\diffeo(\Vector)}^2 \|\diffeo_*[\tangentVectorB_\Vector]\|_{\diffeo(\Vector)}^2 - ( \diffeo_*[\tangentVector_\Vector], \diffeo_*[\tangentVectorB_\Vector])_{\diffeo(\Vector)}^2} \\
        = \frac{( \curvature_{\diffeo(_\Vector)}(\diffeo_* [\tangentVector_\Vector], \diffeo_* [\tangentVectorB_\Vector])\diffeo_* [\tangentVectorB_\Vector], \diffeo_*[\tangentVector_\Vector])_{\diffeo(\Vector)}}{\|\diffeo_*[\tangentVector_\Vector]\|_{\diffeo(\Vector)}^2 \|\diffeo_*[\tangentVectorB_\Vector]\|_{\diffeo(\Vector)}^2 - ( \diffeo_*[\tangentVector_\Vector], \diffeo_*[\tangentVectorB_\Vector])_{\diffeo(\Vector)}^2} = \sectionalcurvature_{\diffeo(_\Vector)}(\diffeo_* [\sigma_{\Vector}]) \leq 0.
    \end{multline}
\end{proof}
\section{Supplementary material to \cref{sec:numerics-pullback-manifolds}}

\subsection{Pullback manifolds from \cref{sec:numerics-pullback-manifolds-curvature-effects}}
\label{app:pullback-manifolds-numerics}

\subsubsection{The unit hyperboloid}
\label{app:pullback-manifolds-numerics-hyperboloid}
\paragraph{Riemannian geometry.}
The $\dimInd$-dimensional unit hyperboloid $\Hyperboloid^{\dimInd}$ is given by
\begin{equation}
    \Hyperboloid^{\dimInd} := \{\mPoint\in \Real^{\dimInd+1} \, \mid \, \|\mPoint\|_M^2 = \mPoint_1^2 + \ldots + \mPoint_{\dimInd}^2 - \mPoint_{\dimInd+1}^2 = -1, \mPoint_{\dimInd+1}>0\}
\end{equation}
and its tangent space $\tangent_\mPoint \Hyperboloid^{\dimInd}$ at $\mPoint\in \Hyperboloid^\dimInd$ is defined as
\begin{equation}
    \tangent_\mPoint \Hyperboloid^{\dimInd} := \{\tangentVector_\mPoint \in \Real^{\dimInd+1} \, \mid \, (\tangentVector_\mPoint, \mPoint)_M = 0) \}.
\end{equation}
We use the Minkowski inner product on $\Real^{\dimInd+1}$ to construct the Riemannian manifold $(\Hyperboloid^{\dimInd}, (\cdot, \cdot)_M)$.

\paragraph{Canonical chart.}
The mapping $\diffeoC:\Hyperboloid^{\dimInd} \to \Real^\dimInd$ given by 
\begin{equation}
    \diffeoC(\mPoint):= (\mPoint_1, 
    \ldots, \mPoint_\dimInd ) \in \Real^{\dimInd}
    \label{eq:hyperboloid-chart}
\end{equation}
provides a standard chart that covers all of the manifold. Its inverse is given by
\begin{equation}
    \diffeoC^{-1}(\Vector) := (\Vector_1, \ldots, \Vector_\dimInd, \sqrt{\|\Vector\|_2^2 + 1} ) \in \Hyperboloid^{\dimInd}.
\end{equation}

\subsubsection{The unit sphere}
\label{app:pullback-manifolds-numerics-sphere}
\paragraph{Riemannian geometry.}
The $\dimInd$-dimensional unit sphere $\Sphere^{\dimInd}$ is given by
\begin{equation}
    \Sphere^{\dimInd} := \{\mPoint\in \Real^{\dimInd+1} \, \mid \, \|\mPoint\|_2^2 = \mPoint_1^2 + \ldots + \mPoint_{\dimInd+1}^2 = 1 \}
\end{equation}
and its tangent space $\tangent_\mPoint \Sphere^{\dimInd}$ at $\mPoint\in \Sphere^\dimInd$ is defined as
\begin{equation}
    \tangent_\mPoint \Sphere^{\dimInd} := \{\tangentVector_\mPoint \in \Real^{\dimInd+1} \, \mid \, (\tangentVector_\mPoint, \mPoint)_2 = 0) \}.
\end{equation}
We use the standard Euclidean inner product on $\Real^{\dimInd+1}$ to construct the Riemannian manifold $(\Sphere^{\dimInd}, (\cdot, \cdot)_2)$.

\paragraph{Chart.}
The mapping $\diffeoC: \Sphere^{\dimInd}\setminus\{(0, \ldots, 0, 1)\} \to \Real^\dimInd$ given by 
\begin{equation}
    \diffeoC(\mPoint):= \biggl(\frac{\mPoint_1}{1 - \mPoint_{\dimInd+1}}, 
    \ldots, \frac{\mPoint_\dimInd}{1 - \mPoint_{\dimInd+1}} \biggr) \in \Real^{\dimInd}
    \label{eq:sphere-chart}
\end{equation}
provides a standard chart that covers almost all of the manifold. Its inverse is given by
\begin{equation}
    \diffeoC^{-1}(\Vector) := \biggl(\frac{2\Vector_1}{1 + \|\Vector\|_2^2 }, \ldots, \frac{2\Vector_\dimInd}{1 + \|\Vector\|_2^2}, \frac{\|\Vector\|_2^2 -1}{1 +\|\Vector\|_2^2} \biggr) \in \Sphere^{\dimInd}.
\end{equation}

\subsection{Error metrics for the evaluation of pulled back geometries}
\label{app:pullback-manifold-error-metrics-numerics}
For comparison of the different pullback geometries we can consider several error metrics as a sanity check. 

\paragraph{Geodesic errors.}
In particular, given a discrete geodesic on $\Real^\dimInd$, i.e., an ordered 1-dimensional data set $\Vector^1, \ldots, \Vector^\dataPointNum\in \Real^{\dimInd}$, we can test whether such a discrete geodesic is an approximate geodesic on a pullback manifold $(\Real^\dimInd, (\cdot,\cdot)^\diffeo)$ through considering the \emph{geodesic error}
\begin{equation}
    \frac{1}{\dataPointNum}\sum_{\sumIndC=1}^{\dataPointNum} \|\geodesic^\diffeo_{\Vector^1, \Vector^{\dataPointNum}} (t_{\sumIndC}) - \Vector^{\sumIndC} \|_2, \quad \text{where } t_\sumIndC := \left\{\begin{matrix}
0 & \text{if } \sumIndC=1, \\
\frac{\sum_{\sumIndA=1}^{\sumIndC - 1} \|\Vector^{\sumIndA} - \Vector^{\sumIndA + 1}\|_2}{\sum_{\sumIndA=1}^{\dataPointNum - 1} \|\Vector^{\sumIndA} - \Vector^{\sumIndA + 1}\|_2} & \text{if }  \sumIndC=2, \ldots, \dataPointNum, \\
\end{matrix}\right.
\label{eq:app-geodesic-error}
\end{equation}
and test stability of the pullback geodesics through considering the \emph{(geodesic) variation error} with respect to $\Vector^1$ to a new point $\VectorC\in \Real^\dimInd$ close to $\Vector^1$
\begin{equation}
    \frac{1}{\dataPointNum}\sum_{\sumIndC=1}^{\dataPointNum} \|\geodesic^\diffeo_{\VectorC, \Vector^{\dataPointNum}} (t_{\sumIndC}) - \geodesic^\diffeo_{\Vector^1, \Vector^{\dataPointNum}} (t_{\sumIndC}) \|_2.
    \label{eq:app-geodesic-variation-error}
\end{equation}


\end{document}